\documentclass[12pt]{article}

\usepackage{amsmath,amssymb,amsthm,geometry,nicefrac,enumerate,mathtools,hyperref,comment,bbm,mathrsfs,mathrsfs,
color,verbatim,amsfonts,graphicx,bm,dsfont}

\usepackage{datetime}
\usepackage[inline]{enumitem}
\usepackage[capitalise,sort]{cleveref}

\crefname{enumi}{item}{items}
\crefname{equation}{}{}

\usepackage{cite}

\usepackage[english]{babel}
\usepackage[utf8]{inputenc}

\newcommand{\with}{\curvearrowleft}

\theoremstyle{plain}

\newtheorem{lemma}{Lemma}[section]

\newtheorem{theorem}[lemma]{Theorem}
\newtheorem{definition}[lemma]{Definition}
\newtheorem{corollary}[lemma]{Corollary}
\newtheorem{proposition}[lemma]{Proposition}

\numberwithin{equation}{section}
\hypersetup{colorlinks=true, linkcolor = {blue}}

\usepackage{mleftright}
\usepackage{mathtools}
\DeclarePairedDelimiter{\pr}{(}{)}

\newcommand{\pa}[1]{\left({#1}\right)}
\newcommand{\pb}[1]{\left[{#1}\right]}
\newcommand{\pc}[1]{\left\{{#1}\right\}}
\newcommand{\pabs}[1]{\left\lvert{#1}\right\rvert}

\newcommand{\R}{\mathbb{R}}
\newcommand{\N}{\mathbb{N}}

\newcommand{\smallsum}{\textstyle\sum}

\providecommand{\B}{\mathcal{B}}

\newcommand{\norm}[1]{\cfadd{def:p-norm}\lVert #1 \rVert}

\newcommand{\normmm}[1]{{\left\vert\kern-0.25ex\left\vert\kern-0.25ex\left\vert #1
    \right\vert\kern-0.25ex\right\vert\kern-0.25ex\right\vert}}

\newcommand{\abs}[1]{\lvert #1\rvert}
\newcommand{\Abs}[1]{\big| #1 \big|}

\newcommand{\eps}{\varepsilon}

\newcommand{\mf}[1]{\mathfrak{#1}}

\newcommand{\cL}{\mathcal{L}}

\newcommand{\fB}{\mathfrak{B}}

\newcommand{\fL}{\mathfrak{L}}

\newcommand{\fR}{\mathfrak{R}}

\newcommand{\fW}{\mathfrak{W}}

\newcommand{\fa}{\mathfrak{a}}
\newcommand{\fb}{\mathfrak{b}}

\newcommand{\ff}{\mathfrak{f}}
\newcommand{\fg}{\mathfrak{g}}

\newcommand{\fl}{\mathfrak{l}}
\newcommand{\fm}{\mathfrak{m}}
\newcommand{\fn}{\mathfrak{n}}

\newcommand{\bfN}{\mathbf{N}}

\newcommand{\bfP}{\mathbf{P}}

\newcommand{\scrR}{\mathscr{R}}

\usepackage{xparse}
\usepackage{etoolbox}

\NewDocumentEnvironment{cproof}{m}
{\begin{proof}[Proof of \cref{#1}]}%
{\noindent The proof of \cref{#1} is thus complete.
\end{proof}}

\NewDocumentEnvironment{cproof2}{m}
{\begin{proof}[Proof of \cref{#1}]}%
{\noindent This completes the proof of \cref{#1}.
\end{proof}}

\ExplSyntaxOn

\NewDocumentCommand{\enum}{ O{;} m o }
 {
  \my_enum:nnn { #1 } { #2 } { #3 }
 }

\seq_new:N \l__my_enum_seq
\tl_new:N \l__my_enum_item_tl
\prop_const_from_keyval:Nn \l__verbs
{
show = shows ,
imply = implies ,
demonstrate = demonstrates ,
prove = proves ,
establish = establishes ,
ensure = ensures ,
assure = assures
}
\int_new:N \l__number_of_args

\cs_new_protected:Nn \my_enum:nnn
 {
  \seq_set_split:Nnn \l__my_enum_seq { #1 } { #2 }
  \seq_remove_all:Nn \l__my_enum_seq {}
  \int_set_eq:NN \l__number_of_args { \seq_count:N \l__my_enum_seq }
  \seq_use:Nnnn \l__my_enum_seq { ~and~ } { ,~ } { ,~and~ }
  \IfNoValueTF{#3}{}{
    \space
    \int_compare:nNnTF{ \l__number_of_args } < {2}{ \prop_item:Nn \l__verbs {#3} }{ #3 }
  }
 }

\seq_new:N \g_cflist_loaded
\seq_new:N \g_cflist_pending

\NewDocumentCommand{\cfadd}{ m }
{
  \seq_if_in:NnF \g_cflist_loaded { #1 } {
    \seq_if_in:NnF \g_cflist_pending { #1 } {
      \seq_gput_right:Nn \g_cflist_pending { #1 }
    }
  }
}

\NewDocumentCommand{\cfconsiderloaded}{ m }{
  \seq_gput_right:Nn \g_cflist_loaded {#1}
}

\NewDocumentCommand{\cfremove}{ m }
{
  \seq_gremove_all:Nn \g_cflist_pending { #1 }
}

\NewDocumentCommand{\cfload}{ o }
{
  \seq_if_empty:NTF \g_cflist_pending {\unskip} {
    (cf.\ \cref{\seq_use:Nn \g_cflist_pending {,}})\IfValueTF{#1}{#1~}{\unskip}
    \seq_gconcat:NNN \g_cflist_loaded \g_cflist_loaded \g_cflist_pending
    \seq_gclear:N \g_cflist_pending
  }
}

\NewDocumentCommand{\cfclear} {} {
  \seq_gclear:N \g_cflist_loaded
  \seq_gclear:N \g_cflist_pending
}

\NewDocumentCommand{\cfout}{ o }
{
  \seq_if_empty:NTF \g_cflist_pending {\unskip} {
    (cf.\ \cref{\seq_use:Nn \g_cflist_pending {,}})\IfValueTF{#1}{#1~}{\unskip}
    \seq_gclear:N \g_cflist_pending
  }
}

\NewDocumentCommand{\ifnocf} { m } {
  \seq_if_empty:NT \g_cflist_pending { #1 }
}

\ExplSyntaxOff

\ExplSyntaxOn

\bool_new:N \g_noteobserve

\NewDocumentCommand{\nobs}{}{
  \bool_if:nTF { \g_noteobserve } {
    \bool_gset_false:N \g_noteobserve
    note~
  } {
    \bool_gset_true:N \g_noteobserve
    observe~
  }
}

\NewDocumentCommand{\Nobs}{}{
  \bool_if:nTF { \g_noteobserve } {
    \bool_gset_false:N \g_noteobserve
    Note~
  } {
    \bool_gset_true:N \g_noteobserve
    Observe~
  }
}

\ExplSyntaxOff

\NewDocumentEnvironment {athm} {m m} {%
\begin{#1}\label{#2}\global\def\loc{#2}%
}{%
\end{#1}%
}

\NewDocumentEnvironment {adef} {m} {%
\begin{definition}\label{#1}\global\def\loc{#1}%
}{%
\end{definition}%
}

\NewDocumentEnvironment{aproof} {} {%
\begin{proof}[Proof~of~\cref{\loc}]%
}{%
\end{proof}%
}

\ExplSyntaxOff

\newcommand{\idMatrix}{\cfadd{def:identityMatrix}\operatorname{I}}

\newcommand{\pmat}[1]{\begin{pmatrix}#1 \end{pmatrix}}
\newcommand{\asinfnorm}[1]{\cfadd{def:p-norm}\left\lVert #1 \right\rVert_{\infty}}
\newcommand{\infnorm}[1]{\cfadd{def:p-norm}\norm{#1}_{\infty}}

\newcommand{\ANNs}{\cfadd{def:ANN}\mathbf{N}}
\newcommand{\paramANN}{\cfadd{def:ANN}\mathcal{P}}
\newcommand{\lengthANN}{\cfadd{def:ANN}\mathcal{L}}
\newcommand{\inDimANN}{\cfadd{def:ANN}\mathcal{I}}
\newcommand{\outDimANN}{\cfadd{def:ANN}\mathcal{O}}
\newcommand{\dims}{\cfadd{def:ANN}\mathcal{D}}
\newcommand{\hiddenLength}{\cfadd{def:ANN}\mathcal{H}}
\newcommand{\dimANNlevel}{\cfadd{def:ANN}\mathbb{D}}
\newcommand{\network}{\cfadd{def:ANN2}\text{neural network }}

\newcommand{\functionANN}{\cfadd{def:ANN}\mathcal{R}}

\newcommand{\vectorNN}{\cfadd{def:vectorANN}\mathcal{T}}

\newcommand{\sizenetwork}[1]{{\left \vert \kern -0.25ex \left \vert \kern -0.25ex \left \vert #1 \right \vert \kern -0.25ex \right \vert \kern -0.25ex \right \vert}}

\newcommand{\scalarMultANN}[2]{\cfadd{def:ANNscalar}#1\circledast#2}

\newcommand{\parallelizationSpecial}{\cfadd{def:simpleParallelization}\mathbf{P}}

\newcommand{\idRelu}{\cfadd{def:ReLu:identity}\mathfrak{I}}

\newcommand{\compANN}[2]{\cfadd{def:ANNcomposition}{#1 \bullet #2}}
\newcommand{\powANN}[2]{\cfadd{def:iteratedANNcomposition}#1^{\bullet #2}}
\newcommand{\compANNbullet}{\cfadd{def:ANNcomposition} \bullet}

\newcommand{\Rect}{\cfadd{def:Rectifier}\mathfrak R}

\newcommand{\AffineANN}{\cfadd{def:ANN:affine}\mathbf{A}}

\newcommand{\sumANN}{\cfadd{def:ANN:sum}\mathfrak{S}}

\newcommand{\expconst}{\sigma}

\newcommand{\constantR}{R}
\newcommand{\constantfrakc}{\mathfrak{c}}
\newcommand{\constantfrakC}{\mathfrak{C}}
\newcommand{\domainD}{D}

\newcommand{\shiftedabsnetwork}{\mathbb{J}}

\newcommand{\idPowerrr}{\mathscr{I}}

\begin{document}

\title{Lower bounds for artificial neural network \\
approximations:
A proof that shallow neural networks \\
fail to overcome the curse of dimensionality
}

\author{
Philipp Grohs$^1$,
Shokhrukh Ibragimov$^2$,\\
Arnulf Jentzen$^{3, 4}$,
and
Sarah Koppensteiner$^5$
\bigskip
\\
\small{$^1$ Faculty of Mathematics and Research Platform Data Science, University of Vienna,}
\vspace{-0.1cm}\\
\small{Vienna, Austria, e-mail: \texttt{philipp.grohs@univie.ac.at}}
\smallskip
\\
\small{$^2$ Faculty of Mathematics and Computer Science, University of M{\"u}nster,}
\vspace{-0.1cm}\\
\small{M{\"u}nster, Germany, e-mail: \texttt{sibragim@uni-muenster.de}}
\smallskip
\\
\small{$^3$ Faculty of Mathematics and Computer Science, University of M{\"u}nster,}
\vspace{-0.1cm}\\
\small{M{\"u}nster, Germany, e-mail: \texttt{ajentzen@uni-muenster.de}}
\smallskip
\\
\small{$^4$ School of Data Science and Shenzhen Research Institute of Big Data,}
\vspace{-0.1cm}\\
\small{The Chinese University of Hong Kong, Shenzhen, China, e-mail: \texttt{ajentzen@cuhk.edu.cn}}
\smallskip
\\
\small{$^5$ Faculty of Mathematics, University of Vienna,}
\vspace{-0.1cm}\\
\small{Vienna, Austria, e-mail: \texttt{sarah.koppensteiner@univie.ac.at}}
}

\date{\today}

\maketitle

\begin{abstract}
Artificial neural networks (ANNs) have become a very powerful tool in the approximation of high-dimensional functions. Especially, deep ANNs, consisting of a large number of hidden layers, have been very successfully used in a series of practical relevant computational problems involving high-dimensional input data ranging from classification tasks in supervised learning to optimal decision problems in reinforcement learning. There are also a number of mathematical results in the scientific literature which study the approximation capacities of ANNs in the context of high-dimensional target functions. In particular, there are a series of mathematical results in the scientific literature which show that sufficiently deep ANNs have the capacity to overcome the curse of dimensionality in the approximation of certain target function classes in the sense that the number of parameters of the approximating ANNs grows at most polynomially in the dimension $d \in \N$ of the target functions under considerations. In the proofs of several of such high-dimensional approximation results it is crucial that the involved ANNs are sufficiently deep and consist a sufficiently large number of hidden layers which grows in the dimension of the considered target functions. It is the topic of this work to look a bit more detailed to the deepness of the involved ANNs in the approximation of high-dimensional target functions. In particular, the main result of this work proves that there exists a concretely specified sequence of functions which can be approximated without the curse of dimensionality by sufficiently deep ANNs but which cannot be approximated without the curse of dimensionality if the involved ANNs are shallow or not deep enough.
\end{abstract}

\tableofcontents
\newpage
\section{Introduction}

Artificial neural networks (ANNs) have become a very powerful tool in the approximation of high-dimensional functions. Especially, deep ANNs, consisting of a large number of hidden layers, have been very successfully used in a series of practical relevant computational problems involving high-dimensional input data ranging from classification tasks in supervised learning to optimal decision problems in reinforcement learning. 

There are also a large number of mathematical results in the scientific literature which study the approximation capacities of ANNs; see, e.g., Cybenko~\cite{MR1015670}, Funahashi~\cite{FUNAHASHI1989183}, Hornik et al.~\cite{HORNIK1991251,HORNIK1989359}, Leshno et al.~\cite{LESHNO1993861}, Guliyev~\&~Ismailov~\cite{GULIYEV2018296}, Elbr{\"a}chter et al.~\cite{Dmytro2020DNNapproximationtheory}, and the references mentioned therein. Moreover, in the recent years a series of articles have appeared in the scientific literature which study the approximation capacities of ANNs in the context of high-dimensional target functions. In particular, the results in such articles show that deep ANNs have the capacity to overcome the curse of dimensionality in the approximation of certain target function classes in the sense that the number of parameters of the approximating ANNs grows at most polynomially in the dimension $d \in \N$ of the target functions under considerations. For example, we refer to Elbr{\"a}chter et al.~\cite{ElbraechterSchwab2018}, Jentzen et al.~\cite{jentzen2018proofarxiv1809}, Gonon et al.~\cite{GononGrohsEtAl2019,gonon2020deep}, Grohs et al.~\cite{GrohsHerrmann2020,GrohsWurstemberger2018,Salimova2019MonteCarlo}, Kutyniok et al.~\cite{kutyniok2019theoretical}, Reisinger~\&~Zhang~\cite{reisinger2019rectified}, Beneventano et al.~\cite{beneventano2020highdimensional}, Berner et al.~\cite{BernerGrohsJentzen2018}, Hornung et al.~\cite{HornungJentzenSalimova2020}, Hutzenthaler et al.~\cite{hutzenthaler2020proof}, and the overview articles Beck et al.~\cite{Kuckuck2020overview} and E et al.~\cite{Jiequn2020AlgorithmsPDEs} for such high-dimensional ANN approximation results in the numerical approximation of solutions of PDEs and we refer to Barron~\cite{Barron1993UniversalApproximation,Barron1992,Barron1994}, Jones~\cite{Jones1992}, Girosi~\&~Anzellotti~\cite{GirAnz1993}, Donahue et al.~\cite{DonGurvDarSont1997}, Gurvits~\&~Koiran~\cite{GurvKoir1997}, K{\r{u}}rkov{\'a} et al.~\cite{KurSang2008,KurKaiKre1997,KurSang2002,Kurkova2008}, Kainen et al.~\cite{KaiKurSang2009,KaiKurSang2012}, Klusowski~\&~Barron~\cite{KlusBar2018}, Li et al.~\cite{LiTangYu}, and Cheridito et al.~\cite{cheridito2019efficient} for such high-dimensional ANN approximation results in the numerical approximation of certain specific target function classes independent of solutions of PDEs (cf., e.g., also Maiorov~\&~Pinkus~\cite{MaioPinkus1999}, Pinkus~\cite{Pinkus1999}, Guliyev~\&~Ismailov~\cite{GuliIsm2018a}, Petersen~\&~Voigtlaender~\cite{PeteVoigt2018}, and B\"{o}lcskei et al.~\cite{BolcGrohsKutynPete2019} for related results). In the proofs of several of the above named high-dimensional approximation results it is crucial that the involved ANNs are sufficiently deep and consist a sufficiently large number of hidden layers which grows in the dimension of the considered target functions. 

It is the key topic of this work to look a bit more detailed to the deepness of the involved ANNs in the approximation of high-dimensional target functions. More specifically, \cref{thm:main_result} in \cref{sec:lower_and_upper_bounds_for_number_of_parameters_in_NN_approximations} below, which is the main result of this work, proves that there exists a concretely specified sequence of high-dimensional functions which can be approximated without the curse of dimensionality by sufficiently deep ANNs but which cannot be approximated without the curse of dimensionality if the involved ANNs are shallow or not deep enough. In the scientific literature related ANN approximation results can also be found in Daniely~\cite{pmlr-v65-daniely17a}, Eldan~\&~Shamir~\cite{pmlr-v49-eldan16}, and Safran~\&~Shamir~\cite{pmlr-v70-safran17a}. One of the differences between the results in the above named references and the results in this work is, roughly speaking, that the considered target functions in the above named references can be approximated by ANNs with two hidden layers without the curse of dimensionality but not with ANNs with one hidden layer while in this work the considered target functions can only be approximated without the curse of dimensionality if the number of the hidden layers of the approximating ANN grows like the dimensions of the target functions. 

To illustrate the findings of this work in more detail, we now present in the following result, \cref{thm:introduction} below, a special case of \cref{thm:main_result}. Below \cref{thm:introduction} we also add some explanatory comments regarding the mathematical objects appearing in \cref{thm:introduction} and regarding the statement of \cref{thm:introduction}. 
 
\begin{samepage}
\cfclear
\begin{theorem}\label{thm:introduction}
Let $\varphi \colon (\cup_{d \in \N}\R^d) \to \R$ and $\Rect \colon (\cup_{d \in \N} \R^d) \allowbreak \to (\cup_{d \in \N} \R^d)$ satisfy for all $d \in \N$, $x = (x_1, \ldots, x_d) \in \R^d$ that $\varphi(x)=(2 \pi)^{\nicefrac{-d}{2}}\exp(- \frac{1}{2} (\smallsum_{j=1}^d \abs{x_j}^2))$ and $\Rect(x) = (\max\{x_1, 0\}, \allowbreak \ldots, \max\{x_d, 0\})$, let $\ANNs = \cup_{L \in \N} \cup_{ l_0, l_1, \ldots, l_L \in \N } (\times_{k = 1}^L (\R^{l_k \times l_{k-1}} \times \R^{l_k}))$, and let $\functionANN \colon \ANNs \to \allowbreak (\cup_{k, l \in \N} \, \allowbreak C(\R^k, \allowbreak \R^l))$, $\hiddenLength \colon \ANNs \to \N_0$, $\paramANN \allowbreak \colon \allowbreak \ANNs \allowbreak \to \N$, and $\sizenetwork{\cdot} \colon \ANNs \allowbreak \to \R$ \allowbreak satisfy \allowbreak for all $L \in \N$, $l_0, \allowbreak l_1, \ldots, l_L \allowbreak \in \N$, $v_0 \allowbreak \in \R^{l_0}, v_1 \in \R^{l_1}, \ldots, \allowbreak v_{L} \in \R^{l_{L}}$, \allowbreak $\Phi = ((W_1, B_1),\allowbreak \ldots, \allowbreak (W_L,\allowbreak B_L)) = \allowbreak (((W_{1, i, j})_{(i, j) \in \{1, \ldots, l_1\} \times \{1, \ldots, l_{0}\}}, \allowbreak (B_{1, i})_{i \in \{1, \ldots, l_{1}\}}), \allowbreak \ldots, \allowbreak ((W_{L, i, j})_{(i, j) \in \{1, \ldots, l_L\} \times \{1, \ldots, l_{L - 1}\}}, \allowbreak (B_{L, i})_{i \in \{1, \ldots, l_{L}\}})) \in \allowbreak (\times_{k = 1}^L \allowbreak (\R^{l_k \times l_{k-1}} \times \allowbreak \R^{l_k}))$ \allowbreak with \allowbreak $\forall \, k \in \{1, 2, \ldots, L\} \colon v_k = \Rect(W_k v_{k-1} + B_k)$ that \allowbreak $\functionANN(\Phi) \in \allowbreak C(\R^{l_0},\allowbreak \R^{l_L})$, $(\functionANN(\Phi)) (v_0) = W_L v_{L-1} + B_L$, $\hiddenLength(\Phi)=L-1$, $\paramANN(\Phi) = \sum_{k = 1}^L l_k (l_{k-1} + 1)$, and $\sizenetwork{\Phi} = \max_{1 \leq n \leq L} \max_{1 \le i \le l_n} \max_{1 \le j \le l_{n-1}} \allowbreak \max\{\abs{W_{n, i, j}}, \abs{B_{n, i}}\}$. Then there exist continuously differentiable $\ff_d \colon \R^d \to \R$, $d \in \N$, such that for all $\delta \in (0,1]$, $\eps \in (0, \nicefrac{1}{2}]$ there exists $\constantfrakC \in (0, \infty)$ such that
\begin{enumerate}[label=(\roman *)]
\item
\label{item1:thm:introduction} it holds for all $\constantfrakc \in [\constantfrakC, \infty)$, $d \in \N$ that
\begin{equation}\label{eqn:item1:thm:introduction}
\min \! \pc{p \in \N \colon
\pb{
\begin{gathered}
\exists \, \Phi \in \ANNs \colon \, p = \paramANN(\Phi), \, \sizenetwork{\Phi} \le \constantfrakc d^{\constantfrakc}, \\ 
d \le \hiddenLength(\Phi) \le \constantfrakc d, \, \functionANN(\Phi) \in C(\R^d, \R), \\
[\smallint\nolimits_{\R^d}\abs{\pr{\functionANN\pr{\Phi}}(x)- \ff_d(x)}^2\varphi(x)\,dx]^{\nicefrac{1}{2}} \le \eps
\end{gathered}}} \le \constantfrakc d^3
\end{equation}

and
\item
\label{item2:thm:introduction} it holds for all $\constantfrakc \in [\constantfrakC, \infty)$, $d \in \N$ that
\begin{equation}\label{eqn:item2:thm:introduction}
\min \! \pc{p \in \N \colon
\pb{
\begin{gathered}
\exists \, \Phi \in \ANNs \colon \, p = \paramANN(\Phi), \, \sizenetwork{\Phi} \le \constantfrakc d^{\constantfrakc}, \\ 
\hiddenLength(\Phi)\leq \constantfrakc d^{1-\delta}, \, \functionANN(\Phi) \in C(\R^d, \R), \\
[\smallint\nolimits_{\R^d}\abs{\pr{\functionANN\pr{\Phi}}(x)- \ff_d(x)}^2\varphi(x)\,dx]^{\nicefrac{1}{2}} \le \eps
\end{gathered}}} \ge (1 + {\constantfrakc}^{-3})^{(d^{\delta})}.
\end{equation}
\end{enumerate}
\end{theorem}
\end{samepage}

\cref{thm:introduction} above is an immediate consequence of \cref{cor:main_result} in Subsection~\ref{subsec:ann_approximation_without_specifying_the_target_function} below. \cref{cor:main_result}, in turn, follows from \cref{thm:main_result} in Subsection~\ref{subsec:ann_approximation_specifying_the_target_function} below, which is the main result of the article. In the following we provide some explanatory comments regarding the statement of \cref{thm:introduction} and regarding the mathematical objects appearing in \cref{thm:introduction}.

In \cref{thm:introduction} we measure the error between the target function and the realization of the approximating ANN in the $L^2$-sense on the whole $\R^d$, $d \in \N$, with respect to standard normal distribution. In particular, we \nobs that the function $\varphi \colon (\cup_{d \in \N} \R^d) \to \R$ in \cref{thm:introduction} appears in the $L^2$-errors in \cref{item1:thm:introduction,item2:thm:introduction} in \cref{thm:introduction} and describes the densities of the standard normal distribution. More formally, \nobs that for all $d \in \N$ it holds that the function $\R^d \ni x \mapsto \varphi(x) = (2 \pi)^{\nicefrac{-d}{2}} \exp(-\frac{1}{2} (\smallsum_{j=1}^d \abs{x_j}^2)) \in \R$ is nothing else but the density of the $d$-dimensional standard normal distribution.

\cref{thm:introduction} is an approximation result for ANNs with the rectifier function as the activation function and the function $\Rect \colon (\cup_{d \in \N} \R^d) \to (\cup_{d \in \N} \R^d)$ in \cref{thm:introduction} describes multidimensional versions of the rectifier function. More specifically, \nobs that for all $d \in \N$ it holds that the function $\R^d \ni x \mapsto \Rect(x) = ( \max\{ x_1, 0 \} , \ldots, \max\{x_d, 0\}) \in \R^d$ is the $d$-dimensional version of the rectifier activation function $\R \ni x \mapsto \max\{x, 0\} \in \R$.

The set $\ANNs = \cup_{L \in \N} \cup_{ l_0, l_1, \ldots, l_L \in \N } (\times_{k = 1}^L (\R^{l_k \times l_{k-1}} \times \R^{l_k}))$ in \cref{thm:introduction} represents the set of all ANNs and the function $\functionANN \colon \ANNs \to (\cup_{k, l \in \N} C(\R^k, \R^l))$ in \cref{thm:introduction} assigns to each ANN in $\ANNs$ its realization function. More formally, \nobs that for every ANN $\Phi \in \ANNs$ it holds that the function $\functionANN(\Phi) \in ( \cup_{k, l \in \N} C(\R^k, \R^l))$ is the realization function associated to the ANN $\Phi$. 

The function $\hiddenLength \colon \ANNs \to \N_0$ in \cref{thm:introduction} describes the number of hidden layers of the considered ANN, the function $\paramANN \colon \ANNs \to \N$ in \cref{thm:introduction} counts the number of parameters (the number of weights and biases) used to describe the considered ANN, 
and the function $\sizenetwork{\cdot} \colon \ANNs \to \R$ in \cref{thm:introduction} specifies the size of the absolute values of the parameters of the considered ANN. More specificially, \nobs that for every ANN $\Phi \in \bfN$ it holds that $\hiddenLength(\Phi)$ is the number of hidden layers of the ANN $\Phi$, that $\paramANN(\Phi)$ is the number of real parameters used to describe the ANN $\Phi$, and that $\sizenetwork{\Phi}$ is the maximum of the absolute values of the real parameters used to describe the ANN $\Phi$. 

Roughly speaking, \cref{thm:introduction} asserts that there exists a sequence of continuously differentiable target functions $f_d \colon \R^d \to \R$, $d \in \N$, such that for every arbitrarily small prescribed approximation accuracy $\eps \in (0,\nicefrac{1}{2}]$ it holds that the class of all sufficiently deep ANNs can approximate the target functions $f_d \colon \R^d \to \R$, $d \in \N$, without the curse of dimensionality (with the number of ANN parameters growing at most cubically in the dimension $d \in \N$; see \cref{eqn:item1:thm:introduction} in \cref{item1:thm:introduction} in \cref{thm:introduction}) and that the class of all shallow ANNs can only approximate the target functions $f_d \colon \R^d \to \R$, $d \in \N$, with the curse of dimensionality (with the number of ANN parameters growing at least exponentially in the dimension $d \in \N$; 
see \cref{eqn:item2:thm:introduction} in \cref{item2:thm:introduction} in \cref{thm:introduction}). In that sense \cref{thm:introduction} shows for a specific class of target functions that deep ANNs can overcome the curse of dimensionality but shallow ANNs fail to do so.

The remainder of this article is organized as follows. In \cref{sec:basics_on_anns} we briefly recall a few general concepts and results from the scientific literature to describe and operate on ANNs. In \cref{sec:upper_bounds_for_weighted_Gaussian_tails} we establish suitable upper bounds for certain weighted tails of standard normal distributions. In \cref{sec:lower_bounds_for_number_of_parameters_in_ANN_approximations} we use the upper bounds for certain weighted tails of standard normal distributions from \cref{sec:upper_bounds_for_weighted_Gaussian_tails} to establish appropriate lower bounds for the number of parameters of ANNs that approximate certain high-dimensional target functions. In \cref{sec:upper_bounds_for_number_of_parameters_in_ANN_approximations} we establish suitable upper bounds for the number of parameters of ANNs that approximate such high-dimensional 
target functions. In \cref{sec:lower_and_upper_bounds_for_number_of_parameters_in_NN_approximations} we combine the lower bounds from \cref{sec:lower_bounds_for_number_of_parameters_in_ANN_approximations} with the upper bounds from \cref{sec:upper_bounds_for_number_of_parameters_in_ANN_approximations} to establish in \cref{thm:main_result} the main ANN approximation result of this work. 
\cref{thm:introduction} above is a direct consequence of \cref{cor:main_result} in \cref{sec:lower_and_upper_bounds_for_number_of_parameters_in_NN_approximations}, which, in turn, follows from \cref{thm:main_result} in \cref{sec:lower_and_upper_bounds_for_number_of_parameters_in_NN_approximations}.

\section{Basics on artificial neural networks (ANNs)}\label{sec:basics_on_anns}
In this section we briefly recall a few general concepts and results from the scientific literature to describe and operate on ANNs. All the notions and the results in this section are well-known in the scientific literature. In particular, regarding \cref{def:ANN} we refer, e.g., to \cite[Definitions 2.1 and 2.3]{Zimmermann2019spacetime}, regarding \cref{def:ANNcomposition} we refer, e.g., to \cite[Definition 2.5]{Zimmermann2019spacetime}, regarding \cref{def:identityMatrix} we refer, e.g., to \cite[Definition 2.10]{Zimmermann2019spacetime}, regarding \cref{def:iteratedANNcomposition} we refer, e.g., to \cite[Definition 2.11]{Zimmermann2019spacetime}, regarding \cref{def:simpleParallelization} we refer, e.g., to \cite[Definition 2.17]{Zimmermann2019spacetime}, regarding \cref{def:ReLu:identity} we refer, e.g., to \cite[Definition 3.15]{Salimova2019MonteCarlo}, regarding \cref{def:ANN:affine} we refer, e.g., to \cite[Definitions 3.7 and 3.10]{Salimova2019MonteCarlo}, regarding \cref{def:ANNscalar} we refer, e.g., to \cite[Definition 3.13]{Salimova2019MonteCarlo}, regarding \cref{def:ANN:sum} we refer, e.g., to \cite[Definition 3.17]{Salimova2019MonteCarlo}, and regarding \cref{def:vectorANN} we refer, e.g., to \cite[Definition 2.11]{Benno2019Fullerroranalysis}. Moreover, note that
\cref{prop:ANNcomposition_elementary_properties} is, e.g., proved as \cite[Proposition 2.6]{Zimmermann2019spacetime}, note that \cref{lemma:associativity_of_ANN_compositions} is, e.g., proved as \cite[Lemma 2.8]{Zimmermann2019spacetime}, note that \cref{Lemma:PropertiesOfANNenlargementGeometry} is, e.g., proved as \cite[Lemma 2.13]{Zimmermann2019spacetime}, note that \cref{Lemma:PropertiesOfParallelizationEqualLength} is, e.g., proved as \cite[Proposition 2.19]{Zimmermann2019spacetime}, note that \cref{Lemma:PropertiesOfParallelizationEqualLengthDims} is, e.g., proved as \cite[Proposition 2.20]{Zimmermann2019spacetime}, note that \cref{lem:Relu:identity} is, e.g., proved as \cite[Lemma 3.16]{Salimova2019MonteCarlo}, note that \cref{lem:ANNscalar} is, e.g., proved as \cite[Lemma 3.14]{Salimova2019MonteCarlo}, and note that \cref{lem:def:ANNsum} is, e.g., proved as \cite[Lemma 3.18]{Salimova2019MonteCarlo}. The proof of \cref{lem:ANN:affine} is clear and therefore is omitted.

\subsection{Structured description of ANNs}
\label{subsec:stcructured_description_anns}

\begin{definition}\label{def:Rectifier}
We denote by $\fR \colon (\cup_{d \in \N} \R^d) \to (\cup_{d \in \N} \R^d)$ the function which satisfies for all $d \in \N$, $x = (x_1, x_2, \ldots, x_d) \in \R^d$ that $\Rect(x) = (\max\{x_1, 0\}, \max\{x_2, 0\}, \ldots, \max\{x_d, 0\})$.
\end{definition}

\cfclear
\begin{definition}\label{def:ANN}
\cfconsiderloaded{def:ANN} We denote by $\mathbf{N}$ the set given by
\begin{equation}\label{eq:defANN}
\textstyle \ANNs \textstyle = \bigcup_{L \in \N} \bigcup_{ l_0, l_1, \ldots, l_L \in \N } \big( \bigtimes_{k = 1}^{L} (\R^{l_k \times l_{k-1}} \times \R^{l_k})\big)
\end{equation}

\noindent
and we denote by $\functionANN \colon \ANNs \to (\cup_{k, l \in \N}\,C(\R^k,\R^l))$, $\paramANN \colon \ANNs \to \N$, $\lengthANN \colon \ANNs \to \N$, $\inDimANN \colon \ANNs \to \N$, $\outDimANN \colon \ANNs \to \N$, $\hiddenLength \colon \ANNs \to \N_0$, $\dims \colon \ANNs \to (\cup_{ L = 2 }^\infty \N^{L})$, and  $\mathbb{D}_n \colon \ANNs \to \N_0$, $n \in \N_0$, the functions which satisfy for all $L \in\N$, $l_0, l_1,\ldots, l_L \allowbreak \in \N$, $\Phi = ((W_1, B_1),(W_2, B_2),\allowbreak \ldots, (W_{L},\allowbreak B_{L})) \in \allowbreak ( \times_{k = 1}^{L} \allowbreak(\R^{l_k \times l_{k-1}} \times \R^{l_k}))$, $v_0 \in \R^{l_0}, v_1 \in \R^{l_1}, \ldots, \allowbreak v_{L} \in \R^{l_{L}}$, $n \in \N_0$ with $ \forall \, k \in \{1, 2, \ldots, L\} \colon v_k = \Rect(W_k v_{k-1} + B_k)$ that $\functionANN(\Phi) \in C(\R^{l_0},\R^{l_L})$, $(\functionANN(\Phi)) (v_0) = W_{L} v_{L - 1} + B_{L}$, $\paramANN(\Phi) = \sum_{k = 1}^{L} l_k(l_{k-1} + 1)$, $\lengthANN(\Phi)=L$,  $\inDimANN(\Phi)=l_0$, $\outDimANN(\Phi)=l_L$, $\hiddenLength(\Phi)=L-1$, $\dims(\Phi)= (l_0, l_1, \ldots, l_L)$, and
\begin{align}
\label{def:ANN:eq1}
\begin{split}
\mathbb{D}_n (\Phi) =
\begin{cases}
l_n &\colon n \leq L \\
0 &\colon n > L\ifnocf.
\end{cases}
\end{split}
\end{align}
\cfload[.]
\end{definition}

\cfclear
\begin{definition}[Neural network]\label{def:ANN2}
We say that $\Phi$ is a neural network if and only if it holds that $\Phi \in \ANNs$ \cfload.
\end{definition}

\subsection{Compositions of ANNs}
\label{subsec:compositions_anns}

\cfclear
\begin{definition}[Compositions of ANNs]
\label{def:ANNcomposition}
\cfconsiderloaded{def:ANNcomposition} We denote by $\compANN{(\cdot)}{(\cdot)} \colon \allowbreak \{(\Phi_1,\Phi_2)\allowbreak \in \ANNs \times \ANNs \colon \inDimANN(\Phi_1)=\outDimANN(\Phi_2)\}\allowbreak \to \ANNs$ the function which satisfies for all $L,\fL \in \N$, $l_0,l_1,\ldots, l_L, \fl_0,\fl_1,\ldots, \fl_\fL \in \N$, $ \Phi_1 = ((W_1, B_1),(W_2, B_2),\allowbreak \ldots, (W_L,\allowbreak B_L)) \in \allowbreak (\times_{k = 1}^L\allowbreak(\R^{l_k \times l_{k-1}} \times \R^{l_k}))$, $ \Phi_2 = ((\fW_1, \fB_1),\allowbreak(\fW_2, \fB_2),\allowbreak \ldots, (\fW_\fL,\allowbreak \fB_\fL)) \in  \allowbreak (\times_{k = 1}^\fL \allowbreak (\R^{\fl_k \times \fl_{k-1}} \times \R^{\fl_k}))$ with $l_0=\inDimANN(\Phi_1)=\outDimANN(\Phi_2)=\fl_{\fL}$ that
\begin{equation}\label{ANNoperations:Composition}
\begin{split}
&\compANN{\Phi_1}{\Phi_2} = \\&
\begin{cases}
\begin{array}{r}
\big((\fW_1, \fB_1),(\fW_2, \fB_2),\ldots, (\fW_{\fL-1},\allowbreak \fB_{\fL-1}),
(W_1 \fW_{\fL}, W_1 \fB_{\fL}+B_{1}),\\ (W_2, B_2), (W_3, B_3),\ldots,(W_{L},\allowbreak B_{L})\big)
\end{array}
& \colon L>1<\fL \\[3ex]
\big((W_1 \fW_{1}, W_1 \fB_1+B_{1}), (W_2, B_2), (W_3, B_3),\ldots,(W_{L},\allowbreak B_{L}) \big)
& \colon L>1=\fL\\[1ex]
\big((\fW_1, \fB_1),(\fW_2, \fB_2),\allowbreak \ldots, (\fW_{\fL-1},\allowbreak \fB_{\fL-1}),(W_1 \fW_{\fL}, W_1 \fB_{\fL}+B_{1}) \big)
& \colon L=1<\fL  \\[1ex]
\bigl((W_1 \fW_{1}, W_1 \fB_1+B_{1})\bigr)
& \colon L=1=\fL
\end{cases}\ifnocf.
\end{split}
\end{equation}
\cfload[.]
\end{definition}





\cfclear
\begin{proposition}\label{prop:ANNcomposition_elementary_properties}
Let $\Phi_1,\Phi_2\in\ANNs$ satisfy $\inDimANN(\Phi_1)=\outDimANN(\Phi_2)$ \cfload. Then
\begin{enumerate}[label=(\roman *)]
\item \label{PropertiesOfCompositions:Dims} it holds that
\begin{equation}
\dims(\compANN{\Phi_1}{\Phi_2}) = (\dimANNlevel_0(\Phi_2), \dimANNlevel_1(\Phi_2), \ldots, \dimANNlevel_{\hiddenLength(\Phi_2)}(\Phi_2), \dimANNlevel_1(\Phi_1), \dimANNlevel_2(\Phi_1), \ldots, \dimANNlevel_{\lengthANN(\Phi_1)}(\Phi_1)),
\end{equation}

\item \label{PropertiesOfCompositions:HiddenLength} it holds that $\hiddenLength(\compANN{\Phi_1}{\Phi_2})=\hiddenLength(\Phi_1)+\hiddenLength(\Phi_2)$,

\item \label{PropertiesOfCompositions:Realization1} it holds that $\functionANN(\compANN{\Phi_1}{\Phi_2})\in C(\R^{\inDimANN(\Phi_2)},\R^{\outDimANN(\Phi_1)})$, and

\item \label{PropertiesOfCompositions:Realization2} it holds that $\functionANN(\compANN{\Phi_1}{\Phi_2})=[\functionANN(\Phi_1)]\circ [\functionANN(\Phi_2)]$\ifnocf.
\end{enumerate}
\cfout[.]
\end{proposition}

\cfclear
\begin{lemma}\label{lemma:associativity_of_ANN_compositions}
Let $\Phi_1, \Phi_2, \Phi_3 \in \ANNs$ satisfy $\inDimANN(\Phi_1) = \outDimANN(\Phi_2)$ and $\inDimANN(\Phi_2) = \outDimANN(\Phi_3)$\cfload. Then $(\Phi_1 \compANNbullet \Phi_2) \compANNbullet \Phi_3 = \Phi_1 \compANNbullet (\Phi_2 \compANNbullet \Phi_3)$\cfout.
\end{lemma}

\subsection{Powers of ANNs}
\label{subsec:powers_anns}

\begin{definition}\label{def:identityMatrix}
Let $n \in \N$. Then we denote by $\idMatrix_n \in \R^{n \times n}$ the identity matrix in $\R^{n \times n}$.
\end{definition}

\cfclear
\begin{definition}\label{def:iteratedANNcomposition}
\cfconsiderloaded{def:iteratedANNcomposition}
We denote by $\powANN{(\cdot)}{n}\colon \{\Phi \in \ANNs \colon \inDimANN(\Phi)=\outDimANN(\Phi)\}\allowbreak\to\ANNs$, $n\in\N_0$, the functions which satisfy for all $n \in \N_0$, $\Phi\in\ANNs$ with $\inDimANN(\Phi)=\outDimANN(\Phi)$ that
\begin{equation}\label{iteratedANNcomposition:equation}
\begin{split}
\powANN{\Phi}{n}=
\begin{cases} \big(\idMatrix_{\outDimANN(\Phi)},(0,0,\dots, 0)\big)\in\R^{\outDimANN(\Phi)\times \outDimANN(\Phi)}\times \R^{\outDimANN(\Phi)} &: n=0 \\
\,\compANN{\Phi}{(\powANN{\Phi}{(n-1)})} &: n \in \N
\end{cases}
\end{split}
\end{equation}
\cfload.
\end{definition}


\cfclear
\begin{lemma}\label{Lemma:PropertiesOfANNenlargementGeometry}
Let  $d, \mathfrak{i} \in \N$, $\Psi \in \ANNs$ satisfy $\dims(\Psi)=(d,\mathfrak{i},d)$ \cfload. Then it holds for all $n \in \N_0$ that
$\hiddenLength(\powANN{\Psi}{n}) = n$,  $\dims(\powANN{\Psi}{n})\in \N^{n+2}$, and
\begin{equation}\label{BulletPower:Dimensions}
\dims(\powANN{\Psi}{n}) =
\begin{cases}
(d,d) &: n=0\\
(d,\mathfrak{i},\mathfrak{i},\dots,\mathfrak{i},d) &: n\in\N
\end{cases}\ifnocf.
\end{equation}
\cfout[.]
\end{lemma}




\subsection{Parallelizations of ANNs}
\label{subsec:parallelizations_anns}

\cfclear
\begin{definition}[Parallelization of ANNs with the same length]\label{def:simpleParallelization}
\cfconsiderloaded{def:simpleParallelization} Let $n\in\N$. Then we denote by
\begin{equation}
\parallelizationSpecial_{n} \colon \pc{(\Phi_1,\Phi_2,\dots, \Phi_n)\in \ANNs^n \colon \cL(\Phi_1) = \cL(\Phi_2) = \ldots = \cL(\Phi_n) } \to \ANNs
\end{equation}

\noindent
the function which satisfies for all  $L\in\N$,
$ (l_{1,0},l_{1,1},\dots, l_{1,L}), (l_{2,0},l_{2,1},\dots, l_{2,L}),\dots,\allowbreak (l_{n,0},\allowbreak l_{n,1},\allowbreak\dots, l_{n,L})\in\N^{L+1}$,
$\Phi_1=((W_{1,1}, B_{1,1}),(W_{1,2}, B_{1,2}),\allowbreak \ldots, (W_{1,L},\allowbreak B_{1,L}))\in (\times_{k = 1}^L\allowbreak(\R^{l_{1,k} \times l_{1,k-1}} \times \R^{l_{1,k}}))$, $\Phi_2 = ((W_{2,1}, B_{2,1}),\allowbreak(W_{2,2}, B_{2,2}),\allowbreak \ldots, (W_{2,L},\allowbreak B_{2,L})) \in (\times_{k = 1}^L\allowbreak(\R^{l_{2,k} \times l_{2,k-1}} \times \R^{l_{2,k}}))$, \dots, $\allowbreak \Phi_n \allowbreak = \allowbreak ((W_{n,1},\allowbreak B_{n,1}), \allowbreak (W_{n,2}, B_{n,2}), \allowbreak \ldots, (W_{n,L}, \allowbreak B_{n,L})) \in (\times_{k = 1}^L (\R^{l_{n,k} \times l_{n,k-1}} \times \R^{l_{n,k}}))$
that
\begin{equation}\label{parallelisationSameLengthDef}
\begin{split}
\bfP_{n}(\Phi_1,\Phi_2,\dots,\Phi_n)&=
\left(
\pa{\begin{pmatrix}
W_{1,1}& 0& 0& \cdots& 0\\
0& W_{2,1}& 0&\cdots& 0\\
0& 0& W_{3,1}&\cdots& 0\\
\vdots& \vdots&\vdots& \ddots& \vdots\\
0& 0& 0&\cdots& W_{n,1}
\end{pmatrix} ,\begin{pmatrix}B_{1,1}\\B_{2,1}\\B_{3,1}\\\vdots\\ B_{n,1}\end{pmatrix}},\right.
\\&\quad
\pa{\begin{pmatrix}
W_{1,2}& 0& 0& \cdots& 0\\
0& W_{2,2}& 0&\cdots& 0\\
0& 0& W_{3,2}&\cdots& 0\\
\vdots& \vdots&\vdots& \ddots& \vdots\\
0& 0& 0&\cdots& W_{n,2}
\end{pmatrix} ,\begin{pmatrix}B_{1,2}\\B_{2,2}\\B_{3,2}\\\vdots\\ B_{n,2}\end{pmatrix}}
,\dots,
\\&\quad\left.
\pa{\begin{pmatrix}
W_{1,L}& 0& 0& \cdots& 0\\
0& W_{2,L}& 0&\cdots& 0\\
0& 0& W_{3,L}&\cdots& 0\\
\vdots& \vdots&\vdots& \ddots& \vdots\\
0& 0& 0&\cdots& W_{n,L}
\end{pmatrix} ,\begin{pmatrix}B_{1,L}\\B_{2,L}\\B_{3,L}\\\vdots\\ B_{n,L}\end{pmatrix}}\right)
\end{split}
\end{equation}
\cfload.
\end{definition}

\cfclear
\begin{proposition}\label{Lemma:PropertiesOfParallelizationEqualLength}
Let $n\in\N$, $\Phi=(\Phi_1,\Phi_2,\allowbreak\dots,\allowbreak \Phi_n)\in\ANNs^n$ satisfy $\lengthANN(\Phi_1)=\lengthANN(\Phi_2)=\ldots=\lengthANN(\Phi_n)$ \cfload. Then
 \begin{enumerate}[label=(\roman{*})]
\item
\label{PropertiesOfParallelizationEqualLength:ItemOne} it holds that $\functionANN(\parallelizationSpecial_{n}(\Phi))\in C(\R^{[\sum_{j=1}^n \inDimANN(\Phi_j)]},\R^{[\sum_{j=1}^n \outDimANN(\Phi_j)]})$ and

\item
\label{PropertiesOfParallelizationEqualLength:ItemTwo} it holds for all    $x_1\in\R^{\inDimANN(\Phi_1)},x_2\in\R^{\inDimANN(\Phi_2)},\dots, x_n\in\R^{\inDimANN(\Phi_n)}$ that
\begin{equation}\label{PropertiesOfParallelizationEqualLengthFunction}
\begin{split}
&\big( \functionANN\big(\parallelizationSpecial_{n}(\Phi)\big) \big) (x_1,x_2,\dots, x_n)
\\&=\big((\functionANN(\Phi_1))(x_1), (\functionANN(\Phi_2))(x_2),\dots,
(\functionANN(\Phi_n))(x_n) \big)\in \R^{[\sum_{j=1}^n \outDimANN(\Phi_j)]}
\end{split}
\end{equation}
\end{enumerate}
\cfout[.]
\end{proposition}

\cfclear
\begin{proposition}\label{Lemma:PropertiesOfParallelizationEqualLengthDims}
Let $n \in \N$, $\Phi_1,\Phi_2,\allowbreak\dots,\allowbreak \Phi_n\in\ANNs$ satisfy $\lengthANN(\Phi_1)= \lengthANN(\Phi_2)=\ldots =\lengthANN(\Phi_n)$ \cfload. Then
\begin{equation}
\dims\big(\parallelizationSpecial_{n}(\Phi_1, \Phi_2, \allowbreak \dots, \allowbreak \Phi_n)\big) = \big(\smallsum_{j=1}^n \dimANNlevel_0(\Phi_j), \smallsum_{j=1}^n \dimANNlevel_1(\Phi_j),\dots, \smallsum_{j=1}^n \dimANNlevel_L(\Phi_j)\big)\ifnocf.
\end{equation}
\cfout[.]
\end{proposition}

\cfclear
\begin{definition}\label{def:ReLu:identity}
\cfconsiderloaded{def:ReLu:identity} We denote by $\idRelu = (\idRelu_d)_{d \in \N} \colon \N \to \ANNs$ the function which satisfies for all $d \in \N$ that
\begin{equation}\label{eq:def:id:1}
\idRelu_1 = \pa{\pa{\begin{pmatrix}
1 \\
-1
\end{pmatrix}, \begin{pmatrix}
0 \\
0
\end{pmatrix}}, \begin{pmatrix}
\begin{pmatrix}
1 & -1
\end{pmatrix}, 0
\end{pmatrix}} \in \pa{\pa{\R^{2 \times 1} \times \R^2} \times \pa{\R^{1 \times 2} \times \R^1}}
\end{equation}

\noindent
and
\begin{equation}\label{eq:def:id:2}
\idRelu _d = \parallelizationSpecial_d(\idRelu_1, \idRelu_1, \ldots, \idRelu_1)
\end{equation}
\cfload.
\end{definition}


\cfclear
\begin{lemma}\label{lem:Relu:identity}
Let $d \in \N$. Then
\begin{enumerate}[label=(\roman{*})]
\item
\label{item:lem:Relu:dims} it holds that $\dims(\idRelu_d) = (d, 2d, d) \in \N^3$,

\item
\label{item:lem:Relu:cont} it holds that $\functionANN(\idRelu_d) \in C(\R^d, \R^d)$, and

\item
\label{item:lem:Relu:real} it holds for all $x \in \R^d$ that $(\functionANN(\idRelu_d))(x) = x$\ifnocf.
\end{enumerate}
\cfout[.]
\end{lemma}




\subsection{Linear transformations as ANNs}
\label{subsec:linear_transformation_anns}

\cfclear
\begin{definition}[Affine linear transformation NN]\label{def:ANN:affine}
\cfconsiderloaded{def:ANN:affine} Let $m,n\in\N,$ $W \in \R^{m\times n},$ $B \in\R^m$. Then we denote by
$\AffineANN_{W, B} \in (\R^{m\times n} \times \R^m) \subseteq \ANNs$ the \network
given by $\AffineANN_{W, B} = (W, B)$ \cfload.
\end{definition}

\cfclear
\begin{lemma}\label{lem:ANN:affine}
Let $m,n\in\N$, $W \in \R^{m\times n}$, $B \in\R^m$. Then
\begin{enumerate}[label=(\roman *)]
\item
\label{lem:ANN:affine:item1} it holds that $\dims(\AffineANN_{W, B}) = (n,m) \in \N^2$,

\item
\label{lem:ANN:affine:item2} it holds that $\functionANN(\AffineANN_{W, B}) \in C(\R^n,\R^m)$,
and

\item
\label{lem:ANN:affine:item3} it holds for all $x \in \R^n$ that $(\functionANN(\AffineANN_{W, B}))(x) = Wx + B$
\end{enumerate}
\cfout[.]
\end{lemma}




\subsection{Scalar multiplications of ANNs}
\label{subsec:scalar_mult_anns}

\cfclear
\begin{definition}[Scalar multiplications of ANNs]\label{def:ANNscalar}
\cfconsiderloaded{def:ANNscalar} We denote by $\scalarMultANN{(\cdot)}{(\cdot)} \colon \R \times \ANNs \to \ANNs$ the function which satisfies for all $\lambda \in \R$, $\Phi \in \ANNs$ that $\scalarMultANN{\lambda}{\Phi} = \compANN{\AffineANN_{\lambda \idMatrix_{\outDimANN(\Phi)},0}}{\Phi}$ \cfload.
\end{definition}


\cfclear
\begin{lemma}\label{lem:ANNscalar}
Let $\lambda \in \R$, $\Phi \in \ANNs$ \cfload. Then
\begin{enumerate}[label=(\roman{*})]
\item\label{item:ANN:scalar:1} it holds that $\dims(\scalarMultANN\lambda\Phi ) = \dims(\Phi)$,

\item\label{item:ANN:scalar:2} it holds that $\functionANN(\scalarMultANN{\lambda}{\Phi}) \in C(\R^{\inDimANN(\Phi)}, \R^{\outDimANN(\Phi)})$, and

\item\label{item:ANN:scalar:3} it holds for all $x \in \R^{\inDimANN(\Phi)}$ that $(\functionANN(\scalarMultANN{\lambda}{\Phi}))(x) = \lambda \big( (\functionANN(\Phi))(x) \big)$\ifnocf.
\end{enumerate}
\cfout[.]
\end{lemma}




\subsection{Sums of ANNs}
\label{subsec:sums_anns}

\cfclear
\begin{definition}\label{def:ANN:sum}
\cfconsiderloaded{def:ANN:sum} Let $m, n \in \N$. Then we denote by $\sumANN_{m, n} \in (\R^{m \times (mn)} \times \R^m)$ the \network \allowbreak given by $\sumANN_{m, n} = \AffineANN_{(\idMatrix_m \,\,\,  \idMatrix_m \,\,\, \ldots \,\,\, \idMatrix_m), \, 0}$ \cfload.
\end{definition}

\cfclear
\begin{lemma}\label{lem:def:ANNsum}
Let $m, n \in \N$. Then
\begin{enumerate}[label=(\roman{*})]
\item\label{item:ANNsum:2} it holds that $\dims (\sumANN_{m, n}) = (mn, m) \in \N^2$,
\item\label{item:ANNsum:3} it holds that $\functionANN(\sumANN_{m, n}) \in C(\R^{mn}, \R^m)$, and
\item\label{item:ANNsum:4} it holds for all $x_1, x_2, \ldots, x_n \in \R^{m}$ that $(\functionANN(\sumANN_{m, n})) (x_1, x_2, \ldots, x_n) = \smallsum_{k=1}^n x_k$\ifnocf.
\end{enumerate}
\cfout[.]
\end{lemma}





\subsection{On the connection to the vectorized description of ANNs}
\label{subsec:connection_to_vectorized_anns}

\begin{definition}[$p$-norm]
\label{def:p-norm}
We denote by
$ \norm{\cdot}_p \colon \bigl(\bigcup_{d=1}^{\infty} \R^d \bigr) \to [0, \infty)$, $p \in [1, \infty]$, the functions which satisfy for all $p \in [1, \infty)$, $d \in \N$, $\theta = (\theta_1, \theta_2, \ldots, \theta_d) \in \R^d$ that $\norm{\theta}_p = \bigl( \smallsum_{i=1}^d \lvert \theta_i \rvert^p \bigr)^{ \!\! \nicefrac{1}{p} }$ and $\lVert \theta \rVert_\infty = \max_{ i \in \{ 1, 2, \ldots, d \} } \lvert \theta_i \rvert$.
\end{definition}

\cfclear
\begin{definition}\label{def:vectorANN}
\cfconsiderloaded{def:vectorANN} We denote by $\vectorNN \colon \ANNs \to \big(\! \bigcup_{d\in\N} \R^d \big)$ the function which satisfies for all $ L,d \in \N $, $ l_0, l_1, \ldots, l_L \in \N $, $ \Phi = \pa{ (W_1, B_1), (W_2, B_2), \ldots, (W_L, B_L) } \in (\times_{m = 1}^L (\R^{l_{m} \times l_{m-1}} \times \R^{l_{m}}))$, $\theta=(\theta_1,\theta_2,\dots,\theta_d)\in\R^d$, $k\in\{1,2,\dots,L\}$ with $\vectorNN(\Phi)=\theta$
that
\begin{equation}\label{eq:translate}
\begin{split}
d &= \paramANN(\Phi), \qquad B_k = \!
\begin{pmatrix}
\theta_{(\sum_{i=1}^{k-1}l_i(l_{i-1}+1))+l_kl_{k-1}+1}\\
\theta_{(\sum_{i=1}^{k-1}l_i(l_{i-1}+1))+l_kl_{k-1}+2}\\
\theta_{(\sum_{i=1}^{k-1}l_i(l_{i-1}+1))+l_kl_{k-1}+3}\\
\vdots\\
\theta_{(\sum_{i=1}^{k-1}l_i(l_{i-1}+1))+l_kl_{k-1}+l_k}
\end{pmatrix}\!,
\qquad \text{and} \\
W_k &=\!
\begin{pmatrix}
\theta_{(\sum_{i=1}^{k-1}l_i(l_{i-1}+1))+1} & \theta_{(\sum_{i=1}^{k-1}l_i(l_{i-1}+1))+2} & \cdots & \theta_{(\sum_{i=1}^{k-1}l_i(l_{i-1}+1))+l_{k-1}} \\
\theta_{(\sum_{i=1}^{k-1}l_i(l_{i-1}+1))+l_{k-1}+1} & \theta_{(\sum_{i=1}^{k-1}l_i(l_{i-1}+1))+l_{k-1}+2} & \cdots & \theta_{(\sum_{i=1}^{k-1}l_i(l_{i-1}+1))+2l_{k-1}} \\
\theta_{(\sum_{i=1}^{k-1}l_i(l_{i-1}+1))+2l_{k-1}+1} & \theta_{(\sum_{i=1}^{k-1}l_i(l_{i-1}+1))+2l_{k-1}+2} & \cdots & \theta_{(\sum_{i=1}^{k-1}l_i(l_{i-1}+1))+3l_{k-1}} \\
\vdots & \vdots & \ddots & \vdots \\
\theta_{(\sum_{i=1}^{k-1}l_i(l_{i-1}+1))+(l_k-1)l_{k-1}+1} & \theta_{(\sum_{i=1}^{k-1}l_i(l_{i-1}+1))+(l_k-1)l_{k-1}+2} & \cdots & \theta_{(\sum_{i=1}^{k-1}l_i(l_{i-1}+1))+l_kl_{k-1}}
\end{pmatrix}
\end{split}
\end{equation}
\cfload.
\end{definition}


\cfclear
\begin{lemma}\label{lem:composition_infnorm}
Let $L, \mathfrak{L} \in \N$, $l_0, l_1, \ldots, l_L, \mathfrak{l}_0, \mathfrak{l}_1, \ldots, \mathfrak{l}_\mathfrak{L} \in \N$, $\Phi_1 = ((W_1, B_1),(W_2, B_2),\allowbreak \ldots, (W_L, \allowbreak B_L)) \in \allowbreak (\times_{k = 1}^L\allowbreak(\R^{l_k \times l_{k-1}} \times \R^{l_k}))$, $\Phi_2 = ((\mathfrak{W}_1, \mathfrak{B}_1),\allowbreak(\mathfrak{W}_2, \mathfrak{B}_2),\allowbreak \ldots, (\mathfrak{W}_\mathfrak{L},\allowbreak \mathfrak{B}_\mathfrak{L})) \in \allowbreak (\times_{k = 1}^\mathfrak{L} \allowbreak(\R^{\mathfrak{l}_k \times \mathfrak{l}_{k-1}} \times \R^{\mathfrak{l}_k}))$ \cfload. Then
\begin{equation}\label{eq:composition_infnorm}
\norm{\vectorNN(\compANN{\Phi_1}{\Phi_2})}_{\infty} \leq \max\bigl\{
\norm{\vectorNN(\Phi_1)}_{\infty}, \norm{\vectorNN(\Phi_2)}_{\infty}, \norm{\vectorNN\bigl(((W_1\mf W_{\mf L},W_1\mf B_{\mf L}+B_1))\bigr)}_{\infty} \bigr\}
\end{equation}
\cfout[.]
\end{lemma}
\begin{cproof}{lem:composition_infnorm} Observe that \cref{ANNoperations:Composition} and \cref{eq:translate} establish \cref{eq:composition_infnorm}.
\end{cproof}




\section{Upper bounds for weighted Gaussian tails}
\label{sec:upper_bounds_for_weighted_Gaussian_tails}

In this section we establish in \cref{lemma:sharp_bounds_specific_integrals} in Subsection~\ref{subsec:bounds_weighted_Gaussian_tails} below 
suitable upper bounds for certain weighted tails of standard normal distributions. Our proof of \cref{lemma:sharp_bounds_specific_integrals} employs the Gaussian segment type estimate in \cref{lemma:Gaussiantail_sharp_bounds} in Subsection~\ref{subsec:bounds_Gaussian_tails} below, the elementary integration formula for certain radial symmetric functions in \cref{lemma:transformation_Rd_to_R} in Subsection~\ref{subsec:bounds_Gaussian_tails}, and the Gaussian tail estimate in \cref{cor:Gaussiantail_alpha=2d} in Subsection~\ref{subsec:bounds_Gaussian_tails}.

\cref{lemma:transformation_Rd_to_R} is a direct consequence of the integral transformation theorem and only for completeness we include in Subsection~\ref{subsec:bounds_Gaussian_tails} the detailed proof for \cref{lemma:transformation_Rd_to_R}. Our proof of \cref{lemma:Gaussiantail_sharp_bounds} uses \cref{lemma:transformation_Rd_to_R} and the elementary estimates for the Gamma function which we present in \cref{cor:Gamma&Stirling} in Subsection~\ref{subsec:bounds_evaluations_Gamma_function} below. Our proof of \cref{cor:Gamma&Stirling} employs the well-known representation result for the Gamma function in \cref{cor:GammaBeta} in Subsection~\ref{subsec:bounds_evaluations_Gamma_function} and the elementary estimates for factorials in \cref{cor:Stirling} in Subsection~\ref{subsec:bounds_evaluations_Gamma_function}. Our proof of \cref{cor:Stirling}, in turn, uses the well-known Stirling inequalities which we recall in \cref{lemma:Stirling} in Subsection~\ref{subsec:bounds_evaluations_Gamma_function}. \cref{lemma:Stirling} is, e.g., proved in Robbins~\cite{MR69328}. The equality in \cref{eqn:lemma:Stirling:Wallis_formula} in the proof of \cref{lemma:Stirling} is also referred to as Wallis's formula in the scientific literature. Our proof of \cref{cor:GammaBeta} employs well-known functional equations for the Gamma and the Beta function which we briefly recall in \cref{lemma:GammaBeta} in Subsection~\ref{subsec:bounds_evaluations_Gamma_function}. \cref{lemma:GammaBeta} is, e.g., proved in Egan~\cite{10.2307/3605337}. Our proof of \cref{cor:Gaussiantail_alpha=2d} uses the well-known Bernoulli inequality which we recall in \cref{lemma:Bernoulli} in Subsection~\ref{subsec:bounds_Gaussian_tails} and the elementary Gaussian tail estimates in \cref{lemma:Gaussiantail_R1} and \cref{cor:Gaussiantail_Rd} in Subsection~\ref{subsec:bounds_Gaussian_tails}. Only for completeness we include in this section also the detailed proofs for \cref{lemma:GammaBeta}, \cref{cor:GammaBeta}, \cref{lemma:Stirling}, and \cref{lemma:Bernoulli}.



\subsection{Lower and upper bounds for evaluations of the Gamma function}
\label{subsec:bounds_evaluations_Gamma_function}



\begin{lemma}\label{lemma:GammaBeta}
Let $\Gamma \colon (0, \infty) \to (0, \infty)$ and $\mathbb{B} \colon (0, \infty)^2 \to (0, \infty)$ satisfy for all $x, y \in (0, \infty)$ that $\Gamma (x) = \int_0^{\infty} t^{x - 1} e^{-t} \, dt$ and $\mathbb{B}(x, y) = \int_0^1 t^{x-1} (1-t)^{y-1} \, dt$. Then
\begin{enumerate}[label=(\roman *)]
\item
\label{item:lemma:GammaBeta1} it holds for all $x \in (0, \infty)$ that $\Gamma (x + 1) = x \, \Gamma (x)$,

\item
\label{item:lemma:GammaBeta2} it holds that $\Gamma(\nicefrac{1}{2})=\sqrt{\pi}$, and

\item
\label{item:lemma:GammaBeta3} it holds for all $x, y \in (0, \infty)$ that $\mathbb{B}(x, y) = \mathbb{B}(y, x) = \frac{\Gamma(x) \Gamma(y)}{\Gamma(x+y)}$.
\end{enumerate}
\end{lemma}
\begin{cproof}{lemma:GammaBeta}
Throughout this proof let $\Phi \colon (0,\infty) \times (0,1) \to (0,\infty)^2$ satisfy for all $u \in (0,\infty)$, $v \in (0,1)$ that $\Phi(u, v) = (u (1 - v), uv)$ and let $f_{x, y} \colon (0,\infty)^2 \to (0,\infty)$, $x, y \in (0, \infty)$, satisfy for all $x, y, s, t \in (0,\infty)$ that $f_{x, y}(s,t)=s^{(x-1)} \, t^{(y-1)} \, e^{-(s+t)}$. \Nobs that the integration by parts formula assures that for all $x \in (0,\infty)$ it holds that
\begin{equation}
\begin{split}
\Gamma(x + 1) & = \int_0^{\infty} t^{((x + 1) - 1)} \, e^{-t} \, dt = - \int_0^{\infty} t^{x} \! \left[ -e^{- t}\right] dt \\
& = - \! \left( \! \left[ t^x e^{-t}\right]^{t = \infty}_{t = 0} - x \!\pb{\int_0^{\infty} t^{(x - 1)} \, e^{-t} \, dt} \right) \! = x \! \pb{\int_0^{\infty} t^{(x - 1)} \,e^{-t} \, dt} \! = x \, \Gamma(x).
\end{split}
\end{equation}

\noindent
This establishes \cref{item:lemma:GammaBeta1}. Next \nobs that the integral transformation theorem shows that
\begin{equation}
\begin{split}
\Gamma\!\pa{\frac{1}{2}} \! & = \int_0^{\infty} t^{\nicefrac{-1}{2}} e^{-t} \, dt = \int_0^{\infty} t^{-1} e^{-t^2} 2t \, dt = 2 \! \pb{\int_0^{\infty} e^{-t^2} \, dt} \\
& = \frac{2}{\sqrt{2}} \! \pb{\int_0^{\infty} e^{-\frac{t^2}{2}} \, dt} \! = \sqrt{\pi} \! \pb{\int_{\R} (2 \pi)^{\nicefrac{-1}{2}} e^{-\frac{t^2}{2}} \, dt} \! = \sqrt{\pi}.
\end{split}
\end{equation}

\noindent
This establishes \cref{item:lemma:GammaBeta2}. Moreover, \nobs that the integral transformation theorem ensures that for all $x, y \in (0, \infty)$ it holds that
\begin{equation}
\begin{split}
\mathbb{B} (x, y) & = \int_0^1 t^{(x - 1)} \, (1 - t)^{(y - 1)} \, dt = \int_1^{\infty} \! \left[\tfrac{1}{t} \right]^{(x - 1)} \left[1 - \tfrac{1}{t} \right]^{(y - 1)} \! \tfrac{1}{t^2} \, dt \\
& = \int_1^{\infty} t^{(- x - 1)} \! \left[\tfrac{t - 1}{ t}\right]^{(y - 1)} dt = \int_1^{\infty} t^{(- x - y)} (t - 1)^{(y - 1)} \, dt \\
& = \int_0^{\infty} (t + 1)^{(- x - y)} t^{(y - 1)} \, dt = \int_0^{\infty} \frac{t^{(y - 1)}}{(t + 1)^{(x + y)}} \, dt.
\end{split}
\end{equation}

\noindent
In addition, \nobs that Fubini's theorem shows that for all $x, y \in (0, \infty)$ it holds that
\begin{equation}\label{eq:pre.transform}
\begin{split}
\Gamma(x) \Gamma(y) & = \! \left[\int_0^{\infty} t^{(x - 1)} \, e^{- t} \, dt \right] \! \left[\int_0^{\infty} t^{(y - 1)} \, e^{- t} \, dt \right] = \! \left[\int_0^{\infty} s^{(x - 1)} \, e^{- s} \, ds \right] \! \left[\int_0^{\infty} t^{(y - 1)} \, e^{- t} \, dt \right] \\
& = \int_0^{\infty} \int_0^{\infty} s^{(x - 1)} \, t^{(y - 1)} \, e^{-(s + t)} \, dt \, ds = \int_{(0,\infty)^2} f_{x, y}(s,t) \, d(s,t).
\end{split}
\end{equation}

\noindent
Furthermore, \nobs that for all $u \in (0, \infty)$, $v \in (0,1)$ it holds that
\begin{equation}
\Phi'(u,v) =
\begin{pmatrix}
1-v & -u \\
v & u
\end{pmatrix}\!.
\end{equation}

\noindent
Hence, we obtain that for all $u \in (0, \infty)$, $v \in (0,1)$ it holds that
\begin{equation}
\operatorname{det}(\Phi'(u,v)) = (1-v)u - v(-u) = u-vu+vu = u \in (0,\infty).
\end{equation}

\noindent
Combining this with \eqref{eq:pre.transform} and the integral transformation theorem shows that for all $x, y \in (0, \infty)$ it holds that
\begin{equation}
\begin{split}
\Gamma(x) \Gamma(y) &= \int_{(0,\infty)\times(0,1)} f_{x, y}(\Phi(u,v)) \, |\!\operatorname{det}(\Phi'(u,v))| \, d(u,v) \\
& = \int_0^\infty \int_0^1 (u(1-v))^{(x-1)} \, (uv)^{(y-1)} \, e^{-(u(1-v)+uv)} \, u \, dv \, du \\
& = \int_0^\infty \int_0^1 u^{(x+y-1)} \, e^{-u} \, v^{(y-1)} \, (1-v)^{(x-1)} \, dv \, du \\
& = \! \left[\int_0^\infty u^{(x+y-1)} \, e^{-u} \, du \right] \! \left[\int_0^1 v^{(y-1)} \, (1-v)^{(x-1)} \, dv \right] \\
& = \! \left[\int_0^\infty u^{(x+y-1)} \, e^{-u} \, du \right] \! \left[\int_0^1 (1-v)^{(y-1)} \, v^{(x-1)} \, dv \right] \\
& = \Gamma(x+y) \, \mathbb{B} (x, y).
\end{split}
\end{equation}

\noindent
This establishes \cref{item:lemma:GammaBeta3}.
\end{cproof}


\begin{corollary}\label{cor:GammaBeta}
Let $\Gamma \colon (0, \infty) \to (0, \infty)$ satisfy for all $x \in (0, \infty)$ that $\Gamma (x) = \int_0^{\infty} t^{x - 1} e^{-t} \, dt$. Then
\begin{enumerate}[label=(\roman *)]
\item
\label{item:cor:GammaBeta1} it holds that $\Gamma(1) = 1$ and
\item
\label{item:cor:GammaBeta2} it holds for all $d \in \N$ that
\begin{equation}
\Gamma\!\pa{\frac{d}{2}}=
\begin{cases}
\big(\frac{d}{2}-1\big)! & \colon \frac{d}{2} \in \N \\[1ex]
\frac{(d-1)!\sqrt{\pi}}{(\frac{d-1}{2})! \, 2^{d-1}} & \colon \frac{d}{2} \notin \N.
\end{cases}
\end{equation}
\end{enumerate}
\end{corollary}
\begin{cproof}{cor:GammaBeta}
\Nobs that the assumption that for all $x \in (0, \infty)$ it holds that $\Gamma(x) = \int_0^{\infty} t^{x-1} e^{-t} \, dt$ ensures that
\begin{equation}
\Gamma(1) = \int_0^{\infty} t^{0} e^{-t} \, dt = \int_0^{\infty} e^{-t} \, dt = [-e^{-t}]_{t = 0}^{t = \infty} = 1.
\end{equation}

\noindent
This establishes \cref{item:cor:GammaBeta1}. Next \nobs that \cref{lemma:GammaBeta}, induction, and \cref{item:cor:GammaBeta1} assure that for all $\fl, \fm \in \N$ with $\fl = 2\fm$ it holds that
\begin{equation}\label{eqn:cor:GammaBeta_even}
\Gamma\!\pa{\frac{\fl}{2}}\!=\Gamma(\fm)=(\fm-1)!=\!\pa{\frac{\fl}{2}-1}!.
\end{equation}

\noindent
Moreover, \nobs that \cref{lemma:GammaBeta} and induction show that for all $\fl, \fm \in \N$ with $\fl = 2\fm-1$ it holds that
\begin{equation}
\begin{split}
\Gamma\!\pa{\frac{\fl}{2}}\! & = \Gamma\!\pa{\fm-\frac{1}{2}}\! = \!\pa{\fm-\frac{3}{2}}\!\Gamma\!\pa{\fm-\frac{3}{2}}\! = \ldots = \!\pa{\fm-\frac{3}{2}}\!\pa{\fm-\frac{5}{2}}\!\cdots\frac{1}{2} \, \Gamma\!\pa{\frac{1}{2}}\! \\
& = \frac{(2\fm-3)!!\sqrt{\pi}}{2^{\fm-1}} = \frac{(2\fm-3)!!(2\fm-2)!!\sqrt{\pi}}{2^{\fm-1}(2\fm-2)!!} = \frac{(2\fm-2)!\sqrt{\pi}}{4^{\fm-1}(\fm-1)!} = \frac{(\fl-1)!\sqrt{\pi}}{2^{\fl-1} \! \pa{\frac{\fl-1}{2}}!}.
\end{split}
\end{equation}

\noindent
Combining this with \cref{eqn:cor:GammaBeta_even} establishes \cref{item:cor:GammaBeta2}.
\end{cproof}

\begin{lemma}\label{lemma:Stirling}
Let $\fn \in \N$. Then
\begin{equation}
\sqrt{2 \pi \fn} \! \pb{\frac{\fn}{e}}^\fn \! e^{\frac{1}{12\fn+1}} < \fn! < \sqrt{2 \pi \fn} \! \pb{\frac{\fn}{e}}^\fn \! e^{\frac{1}{12\fn}}.
\end{equation}
\end{lemma}
\begin{cproof}{lemma:Stirling}
Throughout this proof let $a = (a_n)_{n \in \N} \colon \N \to \R$, $b = (b_n)_{n \in \N} \colon \N \to \R$, $c = (c_n)_{n \in \N} \colon \N \to \R$, $r = (r_{n})_{n \in \N} \colon \N \to [0, \infty]$, and $S = (S_{n})_{n \in \N} \colon \N \to \R$ satisfy for all $n \in \N$ that 
\begin{equation}\label{eqn:lemma:Stirling:a_n_b_n_c_n}
a_n \textstyle = \int_{n}^{n+1} \ln (x) \, dx, \quad b_n = \frac{1}{2} [\ln (n+1) - \ln (n)], \quad c_n = a_n - \frac{1}{2} [\ln (n+1) + \ln (n)],
\end{equation}

\noindent
$r_n = \sum_{k=n}^{\infty} \abs{c_k}$, and $S_{n} = \ln (n!)$, let $C \in [0, \infty]$ satisfy $C = r_1 = \sum_{k=1}^{\infty} \abs{c_k}$, and let $I \colon \N_0 \to \R$ satisfy for all $n \in \N_0$ that $I(n) = \int_0^{\pi} [\sin (x)]^n \, dx$. \Nobs that \cref{eqn:lemma:Stirling:a_n_b_n_c_n} ensures that for all $n \in \N$ it holds that $\ln (n+1) = a_n + b_n - c_n$ and $S_n = \ln (n!) = \textstyle \smallsum_{k=1}^{n-1} \ln (k+1)$. The fact that for all $x \in (0, \infty)$ it holds that $[x \ln (x) - x]'= \ln (x)$ therefore shows that for all $n \in \N$ it holds that
\begin{equation}\label{eqn:lemma:Stirling_S_n_explicit}
\begin{split}
S_n & = \smallsum_{k=1}^{n-1} (a_k + b_k - c_k) = \smallint\nolimits_{1}^{n} \ln (x) \, dx + \tfrac{1}{2} \ln (n) - \smallsum_{k=1}^{n-1} c_k \\
& = [x \ln (x) - x]_{x = 1}^{x = n} + \tfrac{1}{2} \ln (n) - \smallsum_{k=1}^{n-1} c_k = \! \pb{n+\tfrac{1}{2}} \! \ln (n) - n + 1 - \smallsum_{k=1}^{n-1} c_k.
\end{split}
\end{equation}

\noindent
Next \nobs that the fact that for all $x \in (0, \infty)$ it holds that $[x \ln (x) - x]'= \ln (x)$ assures that for all $n \in \N$ it holds that
\begin{equation}\label{eqn:lemma:Stirling_c_n_explicit}
\begin{split}
c_n & = \! \pa{\smallint\nolimits_{n}^{n+1} \ln (x) \, dx} \! - \tfrac{1}{2} [\ln (n+1) + \ln (n)] = [x \ln (x) - x]_{x = n}^{x = n + 1} - \tfrac{1}{2} [\ln (n+1) + \ln (n)] \\
& = \! \pb{n + \tfrac{1}{2}} \! [\ln (n+1) - \ln (n)] - 1 = \! \pb{n + \tfrac{1}{2}} \! \ln\!\pa{1 + \tfrac{1}{n}} \! - 1.
\end{split}
\end{equation}

\noindent
Moreover, \nobs that the fact that for all $x \in (-1, 1)$ it holds that $\ln (1-x) = - \! \pb{\sum_{n=1}^{\infty} \frac{x^n}{n}}$ shows that for all $x \in (-1, 1)$ it holds that
\begin{equation}
\ln\!\pa{\tfrac{1+x}{1-x}} = \ln (1+x) - \ln (1-x) = \! \pb{- \smallsum_{n=1}^{\infty} \tfrac{(-x)^n}{n}} - \! \pb{- \smallsum_{n=1}^{\infty} \tfrac{x^n}{n}} \! = 2 \! \pb{\smallsum_{n=1}^{\infty} \tfrac{x^{2n-1}}{2n-1}}\!.
\end{equation}

\noindent
This implies that for all $n \in \N$ it holds that
\begin{equation}
\begin{split}
\pb{n+\tfrac{1}{2}} \! \ln\!\pa{1+\tfrac{1}{n}} & = \! \pb{n+\tfrac{1}{2}} \! \ln\!\pa{\tfrac{1+(2n+1)^{-1}}{1-(2n+1)^{-1}}} = (2n+1) \! \pb{\smallsum_{p=1}^{\infty} \tfrac{(2n+1)^{1-2p}}{2p-1}} \\
& = \smallsum_{p=1}^{\infty} \tfrac{(2n+1)^{2-2p}}{2p-1} = 1 + \tfrac{1}{3(2n+1)^2} + \tfrac{1}{5(2n+1)^4} + \tfrac{1}{7(2n+1)^6}+\ldots.
\end{split}
\end{equation}

\noindent
Combining this with \cref{eqn:lemma:Stirling_c_n_explicit} ensures that for all $n \in \N$ it holds that
\begin{equation}\label{eqn:lemma:Stirling_c_n_explicit_series}
0 < c_n = \smallsum_{p=2}^{\infty} \tfrac{(2n+1)^{2-2p}}{2p-1} = \smallsum_{p=1}^{\infty} \tfrac{(2n+1)^{-2p}}{2p+1} = \tfrac{1}{3(2n+1)^2} + \tfrac{1}{5(2n+1)^4} + \tfrac{1}{7(2n+1)^6} + \ldots.
\end{equation}

\noindent
The fact that for all $x \in (-1, 1)$ it holds that $\sum_{n=0}^{\infty} x^n = (1-x)^{-1}$ therefore assures that for all $n \in \N$ it holds that
\begin{equation}
\begin{split}
0 < c_n & = \smallsum_{p=1}^{\infty} \tfrac{(2n+1)^{-2p}}{2p+1} < \tfrac{1}{3} \! \pb{\smallsum_{p=1}^{\infty} (2n+1)^{-2p}} = \tfrac{1}{3(2n+1)^2} \! \pb{\smallsum_{p=0}^{\infty} \! \pb{(2n+1)^{-2}}^{p}} \\
& = \tfrac{1}{3(2n+1)^2} (1-(2n+1)^{-2})^{-1} = \tfrac{1}{3} [(2n+1)^2-1]^{-1} = \tfrac{1}{3} [2n(2n+2)]^{-1} \\
& = \tfrac{1}{12} [n(n+1)]^{-1} = \tfrac{1}{12} \big[\tfrac{1}{n(n+1)}\big] =  \tfrac{1}{12} \! \pb{\tfrac{1}{n} - \tfrac{1}{n+1}}\!.
\end{split}
\end{equation}

\noindent
This implies that for all $n \in \N$ it holds that
\begin{equation}\label{eqn:lemma:Stirling_r_n_upper_bound}
r_n = \smallsum_{k=n}^{\infty} \abs{c_k} = \smallsum_{k=n}^{\infty} c_k < \tfrac{1}{12} \! \pa{\sum_{k=n}^{\infty} \! \pb{\tfrac{1}{k} - \tfrac{1}{k+1}}} \! = \tfrac{1}{12n}.
\end{equation}

\noindent
Next \nobs that the fact that for all $x \in (-1, 1)$ it holds that $\sum_{n=0}^{\infty} x^n = (1-x)^{-1}$, the fact that for all $x \in (1, \infty)$ it holds that $3^x > 2x+1$, and \cref{eqn:lemma:Stirling_c_n_explicit_series} show that for all $n \in \N$ it holds that
\begin{equation}\label{eqn:lemma:Stirling_c_n_lower_bound}
\begin{split}
c_n & = \smallsum_{p=1}^{\infty} \tfrac{(2n+1)^{-2p}}{2p+1} > \smallsum_{p=1}^{\infty} \tfrac{(2n+1)^{-2p}}{3^p} = \tfrac{1}{3(2n+1)^2} \! \pb{\smallsum_{p=0}^{\infty} [3^{-1}(2n+1)^{-2}]^p} \\
& = \tfrac{1}{3(2n+1)^2} (1-3^{-1}(2n+1)^{-2})^{-1} = \tfrac{1}{3(2n+1)^2 - 1} \ge \tfrac{1}{12} \! \pb{\tfrac{1}{n+(12)^{-1}} - \tfrac{1}{n+1+(12)^{-1}}}\!.
\end{split}
\end{equation}

\noindent
This ensures that for all $n \in \N$ it holds that 
\begin{equation}\label{eqn:lemma:Stirling_r_n_lower_bound}
r_n = \smallsum_{k=n}^{\infty} \abs{c_k} = \smallsum_{k=n}^{\infty} c_k > \tfrac{1}{12} \big(\sum_{k=n}^{\infty} \! \big[\tfrac{1}{k+(12)^{-1}} - \tfrac{1}{k+1+(12)^{-1}}\big]\big) \! = \tfrac{1}{12n+1}.
\end{equation}

\noindent
In addition, \nobs that the fact that for all $n \in \N$ it holds that $c_n > 0$ and \cref{eqn:lemma:Stirling_r_n_upper_bound} ensure that for all $n \in \N$ it holds that $0 < r_n \le r_1 = C < \infty$. Combining this with \cref{eqn:lemma:Stirling_S_n_explicit} shows that for all $n \in \N$ it holds that
\begin{equation}
\begin{split}
\ln (n!) = S_n & = \! \pb{n+\tfrac{1}{2}} \! \ln (n) - n + 1 - \smallsum_{k=1}^{n-1} c_k = \! \pb{n+\tfrac{1}{2}} \! \ln (n) - n + 1 - C + r_n.
\end{split}
\end{equation}

\noindent
Therefore, we obtain that for all $n \in \N$ it holds that
\begin{equation}
n! = \sqrt{n} \! \pb{\frac{n}{e}}^n e^{r_n} e^{1-C}.
\end{equation}

\noindent
Combining this with \cref{eqn:lemma:Stirling_r_n_lower_bound,eqn:lemma:Stirling_r_n_upper_bound} ensures that for all $n \in \N$ it holds that
\begin{equation}\label{eqn:lemma:Stirling_with_const}
\sqrt{n} \! \pb{\frac{n}{e}}^n e^{\frac{1}{12n+1}} \, e^{1-C} < n! < \sqrt{n} \! \pb{\frac{n}{e}}^n e^{\frac{1}{12n}} \, e^{1-C}.
\end{equation}

\noindent
This implies that for all $n \in \N$ it holds that
\begin{equation}
e^{\frac{1}{12n+1}} < n! \! \pb{\frac{e}{n}}^n n^{\nicefrac{-1}{2}} \, e^{C-1} < e^{\frac{1}{12n}}.
\end{equation}

\noindent
Hence, we obtain that
\begin{equation}\label{eqn:lemma:Stirling_squeeze}
\lim_{n \to \infty} \! \pb{n! \! \pb{\frac{e}{n}}^n n^{\nicefrac{-1}{2}} \, e^{C-1}} \! = e^0 = 1.
\end{equation}

\noindent
Moreover, \nobs that the integration by parts formula and the chain rule ensure that for all $n \in \N \cap [2, \infty)$ it holds that
\begin{equation}
\begin{split}
I(n) & = \textstyle \int_0^{\pi} [\sin (x)]^n \, dx = \int_0^{\pi} [\sin (x)]^{n-1} [\sin (x)] \, dx = - \int_0^{\pi} [\sin (x)]^{n-1} [\frac{d}{d x} \cos (x)] \, dx \\[1ex]
& = \textstyle -[(\sin (x))^{n-1} \cos (x)]_{x=0}^{x=\pi} + \int_0^{\pi} [\frac{d}{dx} ((\sin (x))^{n-1})] [\cos (x)] \, dx \\[1ex]
& \textstyle = (n-1) \int_0^{\pi} [\sin (x)]^{n-2} [\cos (x)]^2 \, dx = \textstyle (n-1) \int_0^{\pi} [\sin (x)]^{n-2} [1 - (\sin (x))^2] \, dx  \\[1ex] 
& = (n-1) [I(n-2) - I(n)].
\end{split}
\end{equation}

\noindent
This implies that for all $n \in \N \cap [2, \infty)$ it holds that
\begin{equation}\label{eqn:lemma:Stirling_I(n)_recurrence}
I(n) = \! \pb{\tfrac{n-1}{n}} \! I(n-2).
\end{equation}

\noindent
Combining this with the fact that $I(0) = \textstyle \int_0^{\pi} \, dx = \pi$ assures that for all $n \in \N$ it holds that
\begin{equation}\label{eqn:lemma:Stirling_I(n)_recurrence_even}
I(2n) = \! \pb{\tfrac{2n-1}{2n}} \! I(2n-2) = \ldots = \! \pb{\tfrac{2n-1}{2n}} \! \pb{\tfrac{2n-3}{2n-2}} \! \cdots \! \pb{\tfrac{1}{2}} \! I(0) = \pi \! \pb{\prod_{k=1}^{n} \! \pb{\tfrac{2k-1}{2k}}}\!.
\end{equation}

\noindent
Next \nobs that the fact that $I(1) = \int_0^{\pi} \sin (x) \, dx = [- \cos (x)]_{x=0}^{x=\pi} = 2$ and \cref{eqn:lemma:Stirling_I(n)_recurrence} demonstrate that for all $n \in \N$ it holds that
\begin{equation}\label{eqn:lemma:Stirling_I(n)_recurrence_odd}
\begin{split}
I(2n+1) & = \! \pb{\tfrac{2n}{2n+1}} \! I(2n-1) = \ldots = \! \pb{\tfrac{2n}{2n+1}} \! \pb{\tfrac{2n-2}{2n-1}} \! \cdots \! \pb{\tfrac{2}{3}} \! I(1) = 2 \! \pb{\prod_{k=1}^{n} \! \pb{\tfrac{2k}{2k+1}}}\!.
\end{split}
\end{equation}

\noindent
Furthermore, \nobs that the fact that for all $n \in \N$, $x \in (0, \pi)$ it holds that $0 < [\sin (x)]^{2n+1} \le [\sin (x)]^{2n} \le [\sin (x)]^{2n-1}$ ensures that for all $n \in \N$ it holds that $0 < I(2n+1) \le I(2n) \le I(2n-1)$. This and \cref{eqn:lemma:Stirling_I(n)_recurrence} imply that for all $n \in \N$ it holds that
\begin{equation}
1 \le \tfrac{I(2n)}{I(2n+1)} \le \tfrac{I(2n-1)}{I(2n+1)} = \tfrac{2n+1}{2n} = 1 + \tfrac{1}{2n}.
\end{equation}

\noindent
Combining this with \cref{eqn:lemma:Stirling_I(n)_recurrence_even,eqn:lemma:Stirling_I(n)_recurrence_odd} shows that
\begin{equation}
1 = \lim_{n \to \infty} \! \pb{\tfrac{I(2n)}{I(2n+1)}} \! = \lim_{n \to \infty} \! \pb{\tfrac{\pi}{2} \prod_{k=1}^{n} \! \pb{\tfrac{(2k-1)(2k+1)}{(2k)^2}}} \! = \tfrac{\pi}{2} \! \pb{\lim_{n \to \infty} \! \pb{\tfrac{[(2n-1)!!]^2 (2n+1)}{[(2n)!!]^2}}}\!.
\end{equation}

\noindent
This demonstrates that
\begin{equation}\label{eqn:lemma:Stirling:Wallis_formula}
\lim_{n \to \infty} \! \pb{\tfrac{(2n)!!}{(2n-1)!! \sqrt{2n}}} \! = \lim_{n \to \infty} \! \pb{\tfrac{(2n)!!}{(2n-1)!! \sqrt{2n+1}}} \! = \! \pb{\tfrac{\pi}{2}}^{\nicefrac{1}{2}}\!.
\end{equation}

\noindent
Combining this with \cref{eqn:lemma:Stirling_squeeze} establishes that
\begin{equation}
\begin{split}
\pb{\tfrac{\pi}{2}}^{\nicefrac{1}{2}}\! & = \lim_{n \to \infty} \! \pb{\tfrac{(2n)!!}{(2n-1)!! \sqrt{2n}}} \! = \lim_{n \to \infty} \! \pb{\tfrac{[(2n)!!]^2}{(2n)! \sqrt{2n}}} \! = \lim_{n \to \infty} \! \pb{\tfrac{[2^n (n!)]^2}{(2n)! \sqrt{2n}}} \! = \lim_{n \to \infty} \! \pb{\tfrac{4^n [n!]^2}{(2n)! \sqrt{2n}}} \\
& = \lim_{n \to \infty} \! \pb{\tfrac{\pb{n! \pb{\frac{e}{n}}^n n^{\nicefrac{-1}{2}} \, e^{C-1}}^2}{(2n)! \pb{\frac{e}{2n}}^{2n} (2n)^{\nicefrac{-1}{2}} \, e^{C-1}}} \lim_{n \to \infty}\! \pb{\tfrac{4^n \pb{\pb{\frac{n}{e}}^n n^{\nicefrac{1}{2}} \, e^{1-C}}^2}{\sqrt{2n} \pb{\frac{2n}{e}}^{2n} (2n)^{\nicefrac{1}{2}} \, e^{1-C}}} \\
& = \lim_{n \to \infty}\! \pb{\tfrac{4^n \pb{\pb{\frac{n}{e}}^n n^{\nicefrac{1}{2}} \, e^{1-C}}^2}{\sqrt{2n} \pb{\frac{2n}{e}}^{2n} (2n)^{\nicefrac{1}{2}} \, e^{1-C}}} \! = \lim_{n \to \infty}\! \pb{\tfrac{4^n n^{2n} e^{-2n} n [e^{1-C}] [e^{1-C}]}{\sqrt{2n} \, (2n)^{2n} e^{-2n} (2n)^{\nicefrac{1}{2}} \, e^{1-C}}} \\
& = \lim_{n \to \infty}\! \pb{\tfrac{2^{2n} n^{2n} n \, e^{1-C}}{(2n) 2^{2n} n^{2n}}} \! = \lim_{n \to \infty} \! \pb{\tfrac{e^{1-C}}{2}} \! = \tfrac{e^{1-C}}{2}.
\end{split}
\end{equation}

\noindent
Hence, we obtain that $e^{1-C} = \sqrt{2 \pi}$. This and \cref{eqn:lemma:Stirling_with_const} show that for all $n \in \N$ it holds that
\begin{equation}
\sqrt{2 \pi n} \! \pb{\frac{n}{e}}^n e^{\frac{1}{12n+1}} < n! < \sqrt{2 \pi n} \! \pb{\frac{n}{e}}^n e^{\frac{1}{12n}}.
\end{equation}
\end{cproof}

\begin{corollary}\label{cor:Stirling}
Let $\fm \in \N \cap [2, \infty)$. Then
\begin{enumerate}[label=(\roman *)]
\item
\label{item:cor:Stirling1} it holds that
\begin{equation}
\sqrt{2 \pi(\fm-1)} \! \pb{\frac{\fm-1}{e}}^{\fm-1} \! \le (\fm-1)! \le  \sqrt{3 \pi(\fm-1)} \! \pb{\frac{\fm-1}{e}}^{\fm-1}
\end{equation}

and
\item
\label{item:cor:Stirling2} it holds that
\begin{equation}
\sqrt{\pi} \! \pb{\frac{\fm-1}{e}}^{\fm-1} \! \le \frac{(2\fm-2)!\sqrt{\pi}}{4^{\fm-1}(\fm-1)!} \le \sqrt{2 \pi} \! \pb{\frac{\fm-1}{e}}^{\fm-1}\!.
\end{equation}
\end{enumerate}
\end{corollary}
\begin{cproof}{cor:Stirling}
\Nobs that \cref{lemma:Stirling} (applied with $\fn \with \fm-1$ in the notation of \cref{lemma:Stirling}) implies that
\begin{equation}\label{eqn:cor:Stirling_m-1}
\sqrt{2 \pi (\fm-1)} \! \pb{\frac{\fm-1}{e}}^{\fm-1} \! e^{\frac{1}{12\fm-11}} \le (\fm-1)! \le \sqrt{2 \pi (\fm-1)} \! \pb{\frac{\fm-1}{e}}^{\fm-1} \! e^{\frac{1}{12\fm-12}}.
\end{equation}

\noindent
The fact that $e \le \! \big(\frac{3}{2}\big)^{(6\fm - 6)}$ therefore assures that
\begin{equation}
\begin{split}
& \sqrt{2 \pi (\fm-1)} \! \pb{\frac{\fm-1}{e}}^{\fm-1} \! \le \sqrt{2 \pi (\fm-1)} \! \pb{\frac{\fm-1}{e}}^{\fm-1} \! e^{\frac{1}{12\fm-11}} \le (\fm-1)! \\
& \le \sqrt{2 \pi (\fm-1)} \! \pb{\frac{\fm-1}{e}}^{\fm-1} \! e^{\frac{1}{12\fm-12}} \le \sqrt{3 \pi (\fm-1)} \! \pb{\frac{\fm-1}{e}}^{\fm-1}\!.
\end{split}
\end{equation}

\noindent
This establishes \cref{item:cor:Stirling1}. Moreover, \nobs that \cref{lemma:Stirling} (applied with $\fn \with 2\fm-2$ in the notation of \cref{lemma:Stirling}) ensures that
\begin{equation}\label{eqn:cor:Stirling_2m-2}
\sqrt{2 \pi (2\fm-2)} \! \pb{\frac{2\fm-2}{e}}^{2\fm-2} \! e^{\frac{1}{24\fm-23}} \le (2\fm-2)! \le \sqrt{2 \pi (2\fm-2)} \! \pb{\frac{2\fm-2}{e}}^{2\fm-2} \! e^{\frac{1}{24\fm-24}}.
\end{equation}

\noindent
Combining this with \cref{eqn:cor:Stirling_m-1} demonstrates that
\begin{equation}
\sqrt{2 \pi} \! \pb{\frac{\fm-1}{e}}^{\fm-1} e^{\frac{-12\fm+11}{(24\fm-23)(12\fm-12)}} \! \le \frac{(2\fm-2)!\sqrt{\pi}}{4^{\fm-1}(\fm-1)!} \le \sqrt{2\pi} \! \pb{\frac{\fm-1}{e}}^{\fm-1} \! e^{\frac{-12\fm+13}{(24\fm-24)(12\fm-11)}}.
\end{equation}

\noindent
The fact that $e^{(11-12 \fm)} \ge 2^{(24\fm - 23)(6 - 6\fm)}$ and the fact that for all $x \in [1, \infty)$ it holds that $x^{(-12m+13)} \le 1$ hence ensure that
\begin{equation}
\begin{split}
& \sqrt{\pi} \! \pb{\frac{\fm-1}{e}}^{\fm-1} \! \le \sqrt{2 \pi} \! \pb{\frac{\fm-1}{e}}^{\fm-1} e^{\frac{-12\fm+11}{(24\fm-23)(12\fm-12)}} \! \le \frac{(2\fm-2)!\sqrt{\pi}}{4^{\fm-1}(\fm-1)!} \\
& \le \sqrt{2\pi} \! \pb{\frac{\fm-1}{e}}^{\fm-1} \! e^{\frac{-12\fm+13}{(24\fm-24)(12\fm-11)}} \le \sqrt{2\pi} \! \pb{\frac{\fm-1}{e}}^{\fm-1}\!.
\end{split}
\end{equation}

\noindent
This establishes \cref{item:cor:Stirling2}.
\end{cproof}





\begin{corollary}\label{cor:Gamma&Stirling}
Let $\Gamma \colon (0, \infty) \to (0, \infty)$ satisfy for all $x \in (0, \infty)$ that $\Gamma (x) = \int_0^{\infty} t^{x - 1} e^{-t} \, dt$. Then
\begin{enumerate}[label=(\roman *)]
\item
\label{item:cor:Gamma&Stirling1} it holds for all $\fm \in \N \cap [2, \infty)$ that
\begin{equation}
\sqrt{2 \pi(\fm-1)} \! \pb{\frac{\fm-1}{e}}^{\fm-1} \! \le \Gamma(\fm) \le  \sqrt{3 \pi(\fm-1)} \! \pb{\frac{\fm-1}{e}}^{\fm-1}
\end{equation}

and
\item
\label{item:cor:Gamma&Stirling2} it holds for all $\fm \in \N \cap [2, \infty)$ that
\begin{equation}
\sqrt{\pi} \! \pb{\frac{\fm-1}{e}}^{\fm-1} \! \le \Gamma\!\pa{\fm-\frac{1}{2}} \! \le \sqrt{2 \pi} \! \pb{\frac{\fm-1}{e}}^{\fm-1}\!.
\end{equation}
\end{enumerate}
\end{corollary}
\begin{cproof}{cor:Gamma&Stirling}
\Nobs that \cref{cor:Stirling} and \cref{item:cor:GammaBeta2} in \cref{cor:GammaBeta} establish \cref{item:cor:Gamma&Stirling1,item:cor:Gamma&Stirling2}.
\end{cproof}



\subsection{Lower and upper bounds for Gaussian tails}
\label{subsec:bounds_Gaussian_tails}


\begin{lemma}\label{lemma:Gaussiantail_R1}
Let $\expconst, s \in (0, \infty)$. Then
\begin{enumerate}[label=(\roman *)]
\item
\label{item:lemma:Gaussiantail_R1_1} it holds that $\int_0^{\infty} e^{-\expconst x^2} \, dx = \frac{\sqrt{\pi}}{2\sqrt{\expconst}}$,

\item
\label{item:lemma:Gaussiantail_R1_2} it holds that
\begin{equation}
\int_s^{\infty} e^{-\expconst x^2} \, dx \le \! \pb{\frac{\sqrt{\pi}}{2\sqrt{\expconst}}} \! e^{-\expconst s^2},
\end{equation}

and
\item
\label{item:lemma:Gaussiantail_R1_3} it holds that
\begin{equation}
\int_0^s e^{-\expconst x^2} \, dx \ge \! \pb{\frac{\sqrt{\pi}}{2\sqrt{\expconst}}} \! \big(1 - e^{-\expconst s^2}\big).
\end{equation}
\end{enumerate}
\end{lemma}
\begin{cproof}{lemma:Gaussiantail_R1}
\Nobs that the integral transformation theorem shows that
\begin{equation}\label{eqn:lemma:Gaussiantail_R1_1}
\int_0^{\infty} e^{- \expconst x^2} \, dx = \frac{1}{\sqrt{2 \expconst}} \int_0^{\infty} \exp(-\tfrac{x^2}{2}) \, dx = \frac{\sqrt{\pi}}{ \sqrt{\expconst}} \int_0^{\infty} \frac{\exp(-\frac{x^2}{2})}{\sqrt{2 \pi}} \, dx = \frac{\sqrt{\pi}}{2 \sqrt{\expconst}}.
\end{equation}

\noindent
This establishes \cref{item:lemma:Gaussiantail_R1_1}. Next \nobs that the integral transformation theorem and \cref{eqn:lemma:Gaussiantail_R1_1} ensure that
\begin{equation}
\label{eqn:lemma:Gaussiantail_R1_2}
\begin{split}
\int_s^{\infty} e^{-\expconst x^2} \, dx & = \int_0^{\infty} e^{-\expconst (x+s)^2} \, dx = \int_0^{\infty} \! \pa{e^{-\expconst x^2 - 2 \expconst s x - \expconst s^2}} \! dx \\
& = e^{-\expconst s^2} \! \pb{\int_0^{\infty} e^{-\expconst x^2 - 2\expconst s x} \, dx} \! \le e^{-\expconst s^2} \! \pb{\int_0^{\infty} e^{-\expconst x^2} \, dx} \! = \! \pb{\frac{\sqrt{\pi}}{2\sqrt{\expconst}}} \! e^{-\expconst s^2}.
\end{split}
\end{equation}

\noindent
This establishes \cref{item:lemma:Gaussiantail_R1_2}. Next we combine \eqref{eqn:lemma:Gaussiantail_R1_1} and \eqref{eqn:lemma:Gaussiantail_R1_2} to obtain that
\begin{equation}
\int_0^s e^{-\expconst x^2} \, dx = \int_0^{\infty} e^{-\expconst x^2} \, dx - \int_s^{\infty} e^{-\expconst x^2} \, dx = \frac{\sqrt{\pi}}{2\sqrt{\expconst}} - \int_s^{\infty} e^{-\expconst x^2} \, dx \ge \! \pb{\frac{\sqrt{\pi}}{2\sqrt{\expconst}}} \! \big(1-e^{-\expconst s^2}\big).
\end{equation}

\noindent
This establishes \cref{item:lemma:Gaussiantail_R1_3}.
\end{cproof}


\cfclear
\begin{corollary}\label{cor:Gaussiantail_Rd}
Let $d \in \N$, $\expconst, s \in (0, \infty)$ \cfload. Then
\begin{enumerate}[label=(\roman *)]
\item
\label{item:cor:Gaussiantail_Rd_1} it holds that
\begin{equation}
\int_{\{y \in \R^d \colon \norm{y}_2 \le s\}} e^{-\expconst \norm{x}_2^2} \, dx \ge \!\pb{\frac{\pi}{\expconst}}^{\nicefrac{d}{2}}\!\pb{1-e^{\nicefrac{-\expconst s^2}{d}}}^d
\end{equation}

and
\item
\label{item:cor:Gaussiantail_Rd_2} it holds that
\begin{equation}
\int_{\{y \in \R^d \colon \norm{y}_2 \ge s\}} e^{-\expconst \norm{x}_2^2} \, dx \le \!\pb{\frac{\pi}{\expconst}}^{\nicefrac{d}{2}}\!\pa{1 - \! \pb{1-e^{\nicefrac{-\expconst s^2}{d}}}^d} \!\ifnocf.
\end{equation}
\end{enumerate}
\cfout.
\end{corollary}
\begin{cproof}{cor:Gaussiantail_Rd}
\Nobs that \cref{item:lemma:Gaussiantail_R1_1} in \cref{lemma:Gaussiantail_R1} implies that
\begin{equation}
\label{eqn:cor:Gaussiantail_Rd_1}
\begin{split}
\int_{\R^d} e^{-\expconst \norm{x}_2^2} \, dx & = \int_{\R} \int_{\R} \ldots \int_{\R} e^{-\expconst \pr{\abs{x_1}^2+\abs{x_2}^2+\ldots+\abs{x_d}^2}} \, dx_d \ldots dx_2 dx_1 \\
& = \! \pb{\int_{\R} e^{-\expconst x^2} \, dx}^d \! = \! \pb{2 \int_0^{\infty} e^{-\expconst x^2} \, dx}^d \! = \!\pb{\frac{\pi}{\expconst}}^{\nicefrac{d}{2}}\ifnocf. 
\end{split}
\end{equation}

\noindent
\cfload. Next \nobs that \cref{item:lemma:Gaussiantail_R1_3} in \cref{lemma:Gaussiantail_R1} (applied with $\expconst \with \expconst$, $s \with d^{\nicefrac{-1}{2}} s$ in the notation of \cref{lemma:Gaussiantail_R1}) and the fact that
\begin{equation}
\{y = (y_1, y_2, \ldots, y_d) \in \R^d \colon (\forall \, j \in \{1, 2, \ldots, d\} \colon \abs{y_j} \le d^{\nicefrac{-1}{2}} s)\} \subseteq \{y \in \R^d \colon \norm{y}_2 \le s\}
\end{equation}

\noindent
ensure that
\begin{equation}
\begin{split}
\int_{\{y \in \R^d \colon \norm{y}_2 \le s \}} e^{-\expconst \norm{x}_2^2} \, dx & \ge \prod_{j=1}^d \! \pb{\int_{-d^{\nicefrac{-1}{2}} s}^{d^{\nicefrac{-1}{2}} s} e^{-\expconst \abs{x_j}^2} \, dx_j} \\
& = \pb{2 \int_{0}^{d^{\nicefrac{-1}{2}} s} e^{-\expconst \abs{x}^2} \, dx}^d \! \ge \! \pb{\frac{\pi}{\expconst}}^{\nicefrac{d}{2}} \! \pb{1-e^{\nicefrac{-\expconst s^2}{d}}}^d \!. 
\end{split}
\end{equation}

\noindent
Combining this with \eqref{eqn:cor:Gaussiantail_Rd_1} establishes \cref{item:cor:Gaussiantail_Rd_1,item:cor:Gaussiantail_Rd_2}.
\end{cproof}

\begin{lemma}\label{lemma:Bernoulli}
Let $\alpha \in \R \backslash (0, 1)$. Then it holds for all $x \in (-1, \infty)$ that $(1+x)^{\alpha} \ge 1+\alpha x$.
\end{lemma}
\begin{cproof}{lemma:Bernoulli}
Throughout this proof let $f \colon (-1, \infty) \to \R$ satisfy for all $x \in (-1, \infty)$ that $f(x) = (1+x)^{\alpha} - 1 - \alpha x$. \Nobs that the chain rule ensures that for all $x \in (-1, \infty)$ it holds that
\begin{equation}\label{eqn:lemma:Bernoulli_derivative}
f'(x) = \alpha (1 + x)^{\alpha - 1} - \alpha = \alpha [(1 + x)^{\alpha - 1} - 1].
\end{equation}

\noindent
The assumption that $\alpha \in \R \backslash (0,1)$ hence ensures that for all $x \in (-1, 0]$ it holds that $f'(x) = \alpha [(1 + x)^{\alpha - 1} - 1] \leq 0$. This implies that the function $(-1,0] \ni x \mapsto f(x) \in \R$ is non-increasing. Hence, we obtain that for all $x \in (-1,0]$ it holds that $f(x) \geq f(0) = 0$. Next \nobs that \cref{eqn:lemma:Bernoulli_derivative} and the assumption that $\alpha \in \R \backslash (0,1)$ demonstrate that for all $x \in [0,\infty)$ it holds that $f'(x) = \alpha [(1 + x)^{\alpha - 1} - 1] \geq 0$. This ensures that the function $[0,\infty) \ni x \mapsto f(x) \in \R$ is non-decreasing. Therefore, we obtain that for all $x \in [0,\infty)$ it holds that $f(x) \geq f(0) = 0$.
\end{cproof}

\cfclear
\begin{corollary}\label{cor:Gaussiantail_alpha=2d}
Let $d \in \N$, $\expconst, s \in (0, \infty)$. Then
\begin{equation}
\int_{\{y \in \R^d \colon \norm{y}_2 \ge s\}} \!\pb{\frac{\expconst}{\pi}}^{\nicefrac{d}{2}} \! e^{- \expconst \norm{x}_2^2} \, dx \le d e^{\nicefrac{- \expconst s^2}{d}} \ifnocf.
\end{equation}
\cfout[.]
\end{corollary}
\begin{cproof}{cor:Gaussiantail_alpha=2d}
\Nobs that the fact that $-e^{\nicefrac{- \expconst s^2}{d}} \in (-1, \infty)$ and \cref{lemma:Bernoulli} (applied with $\alpha \with d$ in the notation of \cref{lemma:Bernoulli}) ensure that $\big(1-e^{\nicefrac{- \expconst s^2}{d}}\big)^d \ge 1-d e^{\nicefrac{- \expconst s^2}{d}}$. Combining this with \cref{item:cor:Gaussiantail_Rd_2} in \cref{cor:Gaussiantail_Rd} (applied with $d \with d$, $\expconst \with \expconst$, $s \with s$ in the notation of \cref{cor:Gaussiantail_Rd}) implies that
\begin{equation}
\int_{\{y \in \R^d \colon \norm{y}_2 \ge s\}\!} \!\pb{\frac{\expconst}{\pi}}^{\nicefrac{d}{2}}\! e^{- \expconst \norm{x}_2^2} \, dx \le 1 - \! \pb{1-e^{\nicefrac{-\expconst s^2}{d}}}^d \! \le 1 - \! \pb{1-de^{\nicefrac{-\expconst s^2}{d}}} \! = d e^{\nicefrac{- \expconst s^2}{d}}\ifnocf.
\end{equation}
\cfload.
\end{cproof}

\cfclear
\begin{lemma}\label{lemma:transformation_Rd_to_R}
Let $d \in \N$, $\expconst \in (0, \infty)$, $\alpha, s \in [0, \infty)$ and let $\Gamma \colon (0, \infty) \to (0, \infty)$ satisfy for all $x \in (0, \infty)$ that $\Gamma(x) = \int_0^{\infty} t^{x-1} e^{-t} \, dt$ \cfload. Then
\begin{equation}
\int_{\{y \in \R^d \colon \norm{y}_2 \ge s\}} \! \pb{\frac{\expconst}{\pi}}^{\nicefrac{d}{2}} \norm{x}_2^{\alpha} \, e^{- \expconst \norm{x}_2^2} \, dx = \frac{2 {\expconst}^{\nicefrac{d}{2}}}{\Gamma\!\pa{\frac{d}{2}}} \! \pb{\int_s^{\infty} e^{- \expconst r^2} r^{\alpha + d - 1} \, dr} \! \ifnocf.
\end{equation}
\cfout[.]
\end{lemma}
\begin{cproof}{lemma:transformation_Rd_to_R}
Throughout this proof assume w.l.o.g. $d>1$, let $\mathbb{B} \colon (0, \infty)^2 \to (0, \infty)$ satisfy for all $x, y \in (0, \infty)$ that $\mathbb{B} (x, y) = \int_0^1 t^{x-1} (1-t)^{y-1} \, dt$, let $\domainD \subseteq \R^{d-1}$ satisfy
\begin{equation}
\domainD = 
\begin{cases}
(0, 2\pi) & \colon d=2 \\
(0, \pi)^{d-2} \times (0, 2\pi) & \colon d>2,
\end{cases}
\end{equation}

\noindent
and let $\Psi \colon (0, \infty) \times \domainD \to \R^d$ satisfy for all $r \in (0, \infty)$, $\varphi = (\varphi_1, \varphi_2, \ldots, \varphi_{d-1}) \in \domainD$ that 
\begin{align}\label{eqn:lemma:transformation_Rd_to_R:Psi}
& \Psi (r, \varphi) \textstyle = \\ 
& \textstyle \bigg(r \cos (\varphi_1) \! \pb{\prod\limits_{k=1}^{0} \sin (\varphi_k)}\!, r \cos (\varphi_2) \! \pb{\prod\limits_{k=1}^{1} \sin (\varphi_k)}\!, \ldots,
r \cos (\varphi_{d-1}) \! \pb{\prod\limits_{k=1}^{d-2} \sin (\varphi_{k})}\!, r \!\pb{\prod\limits_{k=1}^{d-1} \sin(\varphi_k)}\!\bigg). \nonumber
\end{align}

\noindent
\Nobs that \cref{eqn:lemma:transformation_Rd_to_R:Psi} shows that for all $r \in (0, \infty)$, $\varphi = (\varphi_1, \varphi_2, \ldots, \varphi_{d-1}) \in \domainD$ it holds that $\norm{\Psi(r, \varphi)}_2 = r$ and
\begin{equation}
\abs{\operatorname{det}(\Psi'(r, \varphi))} = r^{d-1} \! \pb{\prod_{k=1}^{d-2} [\sin(\varphi_k)]^{d-k-1}}\!\ifnocf.
\end{equation}

\noindent
\cfload[.]The integral transformation theorem hence ensures that
\begin{equation}\label{eqn:lemma:transformation_Rd_to_R:transformation}
\begin{split}
& \int_{\{y \in \R^d \colon \norm{y}_2 \ge s\}} \norm{x}_2^{\alpha} \, e^{- \expconst \norm{x}_2^2} \, dx \\
& = \int_{\R^d} \norm{x}_2^{\alpha} \, e^{- \expconst \norm{x}_2^2} \, \mathbbm{1}_{[s, \infty)} (\norm{x}_2) \, dx \\
& = \int_0^{\infty} \! \int_{\domainD} \norm{\Psi(r, \varphi)}_2^{\alpha} \, e^{- \expconst \norm{\Psi(r, \varphi)}_2^2} \, \abs{\operatorname{det}(\Psi'(r, \varphi))} \, \mathbbm{1}_{[s, \infty)}(\norm{\Psi(r, \varphi)}_2) \, d\varphi \, dr \\
& = \int_s^{\infty} \! \int_{\domainD} \norm{\Psi(r, \varphi)}_2^{\alpha} \, e^{- \expconst \norm{\Psi(r, \varphi)}_2^2} \, \abs{\operatorname{det}(\Psi'(r, \varphi))} \, d\varphi \, dr \\
& = \int_s^{\infty} \! \int_{\domainD} e^{- \expconst r^2} r^{\alpha} \, \abs{\operatorname{det}(\Psi'(r, \varphi))} \, d\varphi \, dr \\
& = 2 \pi \!\pb{\prod_{k=1}^{d-2} \int_0^{\pi} [\sin (x)]^k \, dx} \! \pb{\int_s^{\infty} e^{- \expconst r^2} r^{\alpha + d - 1} \, dr}\!\ifnocf.
\end{split}
\end{equation}

\noindent
\cfload[.]Next \nobs that the chain rule assures that for all $x \in (0, 1)$ it holds that
\begin{equation}
1 = \tfrac{d}{dx} \big( \operatorname{id}_{\R} (x) \big) = \tfrac{d}{dx} \big( \sin(\arcsin(x)) \big) = \cos (\arcsin (x)) [\arcsin'(x)] = [\arcsin'(x)] \sqrt{1 - x^2}.
\end{equation}

\noindent
This implies that for all $x \in (0, 1)$ it holds that
\begin{equation}
\arcsin'(x) = (1 - x^2)^{\nicefrac{-1}{2}}.
\end{equation}

\noindent
The integral transformation theorem and \cref{lemma:GammaBeta} hence show that for all $k \in \N$ it holds that
\begin{equation}
\begin{split}
\int_0^{\pi} [\sin (x)]^k \, dx & = 2 \int_0^{\nicefrac{\pi}{2}} [\sin (x)]^k \, dx = 2 \int_0^1 \! \pb{\frac{x^k}{(1-x^2)^{\nicefrac{1}{2}}}} dx = \int_0^1 x^{\frac{k-1}{2}} (1-x)^{-\frac{1}{2}} \, dx \\
& = \mathbb{B}\!\pa{\frac{k+1}{2}, \frac{1}{2}} = \frac{\Gamma\!\pa{\frac{k+1}{2}} \Gamma\!\pa{\frac{1}{2}}}{\Gamma\!\pa{\frac{k+2}{2}}} = \! \frac{\Gamma\!\pa{\frac{k+1}{2}} \! \sqrt{\pi}}{\Gamma\!\pa{\frac{k+2}{2}}}.
\end{split}
\end{equation}

\noindent
Combining this with \cref{eqn:lemma:transformation_Rd_to_R:transformation} and \cref{item:cor:GammaBeta1} in \cref{cor:GammaBeta} demonstrates that
\begin{equation}
\begin{split}
\int_{\{y \in \R^d \colon \norm{y}_2 \ge s\}} \norm{x}_2^{\alpha} \, e^{- \expconst \norm{x}_2^2} \, dx & = 2 \pi \!\pb{\prod_{k=1}^{d-2} \int_0^{\pi} [\sin (x)]^k \, dx} \! \int_s^{\infty} e^{- \expconst r^2} r^{\alpha + d - 1} \, dr \\
& = 2 \pi \!\pb{\prod_{k=1}^{d-2} \! \pb{\frac{\Gamma\!\pa{\frac{k+1}{2}} \! \sqrt{\pi}}{\Gamma\!\pa{\frac{k+2}{2}}}}} \! \int_s^{\infty} e^{- \expconst r^2} r^{\alpha + d - 1} \, dr \\
& = \frac{2 \pi^{\nicefrac{d}{2}}}{\Gamma\!\pa{\frac{d}{2}}} \! \pb{\int_s^{\infty} e^{- \expconst r^2} r^{\alpha + d - 1} \, dr}\!.
\end{split}
\end{equation}

\noindent
Therefore, we obtain that
\begin{equation}
\int_{\{y \in \R^d \colon \norm{y}_2 \ge s\}} \! \pb{\frac{\expconst}{\pi}}^{\nicefrac{d}{2}} \norm{x}_2^{\alpha} \, e^{- \expconst \norm{x}_2^2} \, dx = \frac{2 {\expconst}^{\nicefrac{d}{2}}}{\Gamma\!\pa{\frac{d}{2}}} \! \pb{\int_s^{\infty} e^{- \expconst r^2} r^{\alpha + d - 1} \, dr}\!.
\end{equation}
\end{cproof}

\cfclear
\begin{lemma}\label{lemma:Gaussiantail_sharp_bounds}
Let $d \in \N \cap [3, \infty)$, $\beta, \expconst \in (0, \infty)$ \cfload. Then
\begin{equation}
\int_{\pc{y \in \R^d \colon \frac{\sqrt{d(1+\beta)}}{\sqrt{2 \expconst}} \le \norm{y}_2 \le \frac{d \sqrt{1+\beta}}{\sqrt{2 \expconst}}}} \! \pb{\frac{\expconst}{\pi}}^{\nicefrac{d}{2}} \! e^{- \expconst \norm{x}_2^2} \, dx \le d \! \pb{\frac{1+\beta}{e^{\beta}}}^{\nicefrac{d}{2}}\!\ifnocf.
\end{equation}
\cfout[.]
\end{lemma}
\begin{cproof}{lemma:Gaussiantail_sharp_bounds}
Throughout this proof let $\Gamma \colon (0, \infty) \to (0, \infty)$ satisfy for all $x \in (0, \infty)$ that $\Gamma(x) = \int_0^{\infty} t^{x-1} e^{-t} \, dt$. \Nobs that \cref{lemma:transformation_Rd_to_R} (applied with $d \with d$, $\expconst \with \expconst$, $\alpha \with 0$, $s \with (2 \expconst)^{\nicefrac{-1}{2}} (d (1+\beta))^{\nicefrac{1}{2}}$ in the notation of \cref{lemma:transformation_Rd_to_R}) implies that
\begin{equation}\label{eqn:lemma:Gaussiantail_sharp_bounds_sqrt_d}
\int_{\pc{y \in \R^d \colon \norm{y}_2 \ge \frac{\sqrt{d(1+\beta)}}{\sqrt{2 \expconst}}}} \! \pb{\frac{\expconst}{\pi}}^{\nicefrac{d}{2}} \! e^{- \expconst \norm{x}_2^2} \, dx = \frac{2 \expconst^{\nicefrac{d}{2}}}{\Gamma\!\pa{\frac{d}{2}}} \int_{\frac{\sqrt{d(1+\beta)}}{\sqrt{2 \expconst}}}^{\infty}  e^{- \expconst r^2} r^{d-1}  \, dr\ifnocf.
\end{equation}

\noindent
\cfload[.]Next \nobs that \cref{lemma:transformation_Rd_to_R} (applied with $d \with d$, $\expconst \with \expconst$, $\alpha \with 0$, $s \with (2 \expconst)^{\nicefrac{-1}{2}} d (1+\beta)^{\nicefrac{1}{2}}$ in the notation of \cref{lemma:transformation_Rd_to_R}) shows that
\begin{equation}
\int_{\pc{y \in \R^d \colon \norm{y}_2 \ge \frac{d \sqrt{1+\beta}}{\sqrt{2 \expconst}}}} \! \pb{\frac{\expconst}{\pi}}^{\nicefrac{d}{2}} \! e^{- \expconst \norm{x}_2^2} \, dx = \frac{2 \expconst^{\nicefrac{d}{2}}}{\Gamma\!\pa{\frac{d}{2}}} \int_{\frac{d \sqrt{1+\beta}}{\sqrt{2 \expconst}}}^{\infty}  e^{- \expconst r^2} r^{d-1}  \, dr.
\end{equation}

\noindent
Combining this with \cref{eqn:lemma:Gaussiantail_sharp_bounds_sqrt_d} ensures that
\begin{equation}\label{eqn:lemma:Gaussiantail_sharp_bounds_substitution_1+delta}
\begin{split}
& \int_{\pc{y \in \R^d \colon \frac{\sqrt{d(1+\beta)}}{\sqrt{2 \expconst}} \le \norm{y}_2 \le \frac{d \sqrt{1+\beta}}{\sqrt{2 \expconst}}}} \! \pb{\frac{\expconst}{\pi}}^{\nicefrac{d}{2}} \! e^{- \expconst \norm{x}_2^2} \, dx \\
& = \int_{\pc{y \in \R^d \colon \norm{y}_2 \ge \frac{\sqrt{d(1+\beta)}}{\sqrt{2 \expconst}}}} \! \pb{\frac{\expconst}{\pi}}^{\nicefrac{d}{2}} \! e^{- \expconst \norm{x}_2^2} \, dx - \int_{\pc{y \in \R^d \colon \norm{y}_2 \ge \frac{d \sqrt{1+\beta}}{\sqrt{2 \expconst}}}} \! \pb{\frac{\expconst}{\pi}}^{\nicefrac{d}{2}} \! e^{- \expconst \norm{x}_2^2} \, dx \\
& = \frac{2 \expconst^{\nicefrac{d}{2}}}{\Gamma\!\pa{\frac{d}{2}}} \int_{\frac{\sqrt{d(1+\beta)}}{\sqrt{2 \expconst}}}^{\infty}  e^{- \expconst r^2} r^{d-1}  \, dr - \frac{2 \expconst^{\nicefrac{d}{2}}}{\Gamma\!\pa{\frac{d}{2}}} \int_{\frac{d \sqrt{1+\beta}}{\sqrt{2 \expconst}}}^{\infty}  e^{- \expconst r^2} r^{d-1}  \, dr \\
& = \frac{2 \expconst^{\nicefrac{d}{2}}}{\Gamma\!\pa{\frac{d}{2}}} \int_{\frac{\sqrt{d(1+\beta)}}{\sqrt{2 \expconst}}}^{\frac{d \sqrt{1+\beta}}{\sqrt{2 \expconst}}}  e^{- \expconst r^2} r^{d-1}  \, dr.
\end{split}
\end{equation}
Next \nobs that the chain rule ensures that for all $x \in [(2 \expconst)^{\nicefrac{-1}{2}} d^{\nicefrac{1}{2}}, \infty)$ it holds that
\begin{equation}\label{eqn:decreasing}
\big[e^{- \expconst x^2} x^{d-1}\big]' \! = e^{- \expconst x^2} x^{d-2} (d-1 - 2 \expconst x^2) \le e^{- \expconst x^2} x^{d-2} (d-1 - d) < 0.
\end{equation}

\noindent
This ensures that the function $[(2 \expconst)^{\nicefrac{-1}{2}} d^{\nicefrac{1}{2}}, \infty) \ni x \mapsto e^{- \expconst x^2} x^{d-1} \in \R$ is strictly decreasing. Hence, we obtain that
\begin{equation}\label{eqn:sharp:item2:upper_bound}
\begin{split}
& \frac{2 \expconst^{\nicefrac{d}{2}}}{\Gamma\!\pa{\frac{d}{2}}} \int_{\frac{\sqrt{d(1+\beta)}}{\sqrt{2 \expconst}}}^{\frac{d \sqrt{1+\beta}}{\sqrt{2 \expconst}}}  e^{- \expconst r^2} r^{d-1}  \, dr \le \frac{2 \expconst^{\nicefrac{d}{2}}}{\Gamma\!\pa{\frac{d}{2}}} \int_{\frac{\sqrt{d(1+\beta)}}{\sqrt{2 \expconst}}}^{\frac{d \sqrt{1+\beta}}{\sqrt{2 \expconst}}} \! \pb{e^{-\frac{d(1+\beta)}{2}}} \! \pb{\frac{d(1+\beta)}{2 \expconst}}^{\frac{d-1}{2}} \! dr \\
& = \frac{2 \expconst^{\nicefrac{d}{2}}}{\Gamma\!\pa{\frac{d}{2}}} \! \pb{e^{-\frac{d(1+\beta)}{2}}} \! \pb{\frac{d(1+\beta)}{2 \expconst}}^{\frac{d-1}{2}} \! \pb{\frac{(d-\sqrt{d})\sqrt{1+\beta}}{\sqrt{2 \expconst}}} \\
& \le \frac{2 \expconst^{\nicefrac{d}{2}}}{\Gamma\!\pa{\frac{d}{2}}} \! \pb{e^{-\frac{d(1+\beta)}{2}}} \! \pb{\frac{d(1+\beta)}{2 \expconst}}^{\frac{d-1}{2}} \! \pb{d \! \pb{\frac{1+\beta}{2\expconst}}^{\nicefrac{1}{2}}} = \frac{2}{\Gamma\!\pa{\frac{d}{2}}} \! \pb{e^{-\frac{d(1+\beta)}{2}}} \! \pb{d^{\frac{d+1}{2}}} \! \pb{\frac{1+\beta}{2}}^{\nicefrac{d}{2}}\!.
\end{split}
\end{equation}

\noindent
Next \nobs that \cref{item:cor:Gamma&Stirling1} in \cref{cor:Gamma&Stirling} and the fact that for all $m \in \N \cap [2, \infty)$ it holds that
\begin{equation}
\pb{\pb{1+\frac{1}{m-1}}^{m-1}}^{\frac{2m-1}{2m-2}} \! \le e \! \pb{\pb{1+\frac{1}{m-1}}^{m-1}}^{\frac{1}{2m-2}} \! = e \! \pb{1+\frac{1}{m-1}}^{\nicefrac{1}{2}} \!  \le e \big[2^{\nicefrac{1}{2}}\big] \le \frac{3e}{2}
\end{equation}

\noindent
assure that for all $k, m \in \N$ with $k = 2m \ge 4$ it holds that
\begin{equation}\label{eqn:sharp:item2:upper_bound_even}
\begin{split}
& \frac{2}{\Gamma\!\pa{\frac{k}{2}}} \! \pb{e^{- \frac{k(1+\beta)}{2}}} \! \pb{k^{\frac{k+1}{2}}} \! \pb{\frac{1+\beta}{2}}^{\nicefrac{k}{2}} = \frac{2}{\Gamma\!\pa{m}} \! \pb{e^{-m(1+\beta)}} \! \pb{(2m)^{m+\frac{1}{2}}} \! \pb{\frac{1+\beta}{2}}^m \\
& \le \frac{2}{\sqrt{2\pi (m-1)}} \! \pb{\frac{e}{m-1}}^{m-1} \! \pb{e^{-m(1+\beta)}} \! \pb{(2m)^{m+\frac{1}{2}}} \! \pb{\frac{1+\beta}{2}}^m \\
& = \frac{2m}{\sqrt{\pi}} \! \pb{e^{-1-m\beta}} \! \pb{\pb{1+\frac{1}{m-1}}^{m-1}}^{\frac{2m-1}{2m-2}} \! (1+\beta)^{m} \le \frac{2m}{\sqrt{\pi}} \! \pb{e^{-1-m\beta}} \! \pb{\frac{3e}{2}} \! (1+\beta)^{m} \\
& = \frac{3k}{2 \sqrt{\pi}} \! \pb{\frac{1+\beta}{e^{\beta}}}^{\nicefrac{k}{2}} \le k \! \pb{\frac{1+\beta}{e^{\beta}}}^{\nicefrac{k}{2}}\!.
\end{split}
\end{equation}

\noindent
Next \nobs that \cref{item:cor:Gamma&Stirling2} in \cref{cor:Gamma&Stirling} and the fact that for all $m \in \N \cap [2, \infty)$ it holds that $(1+(2m-2)^{-1})^{2m-2} \le e$ show that for all $k, m \in \N$ with $k = 2m - 1 \ge 3$ it holds that
\begin{equation}\label{eqn:sharp:item2:upper_bound_odd}
\begin{split}
& \pb{\frac{2}{\Gamma\!\pa{\frac{k}{2}}}} \! \pb{e^{- \frac{k(1+\beta)}{2}} k^{\frac{k+1}{2}}} \! \pb{\frac{1+\beta}{2}}^{\nicefrac{k}{2}} \!  = \! \pb{\frac{2}{\Gamma\!\pa{m-\frac{1}{2}}}} \! \pb{e^{(-m+\frac{1}{2})(1+\beta)} (2m-1)^m} \! \pb{\frac{1+\beta}{2}}^{m-\frac{1}{2}} \\
& \le \frac{2}{\sqrt{\pi}} \! \pb{\frac{e}{m-1}}^{m-1} \! \pb{e^{(-m+\frac{1}{2})(1+\beta)} (2m-1)^m} \! \pb{\frac{1+\beta}{2}}^{m-\frac{1}{2}} \\
& = \! \pb{\frac{2}{\pi}}^{\nicefrac{1}{2}} \! \pb{(2m-1) e^{-m\beta +\frac{\beta}{2}-\frac{1}{2}}} \! \pb{\pb{1+\frac{1}{2m-2}}^{2m-2}}^{\nicefrac{1}{2}} \! (1+\beta)^{m-\frac{1}{2}} \\
& \le \! \pb{\frac{2}{\pi}}^{\nicefrac{1}{2}} \! (2m-1) \, e^{-m\beta +\frac{\beta}{2}-\frac{1}{2}} \, e^{\frac{1}{2}} (1+\beta)^{m-\frac{1}{2}} = \! \pb{\frac{2}{\pi}}^{\nicefrac{1}{2}} \! k \! \pb{\frac{1+\beta}{e^{\beta}}}^{\nicefrac{k}{2}} \le k \! \pb{\frac{1+\beta}{e^{\beta}}}^{\nicefrac{k}{2}}\!.
\end{split}
\end{equation}

\noindent
Combining this with \cref{eqn:sharp:item2:upper_bound,eqn:sharp:item2:upper_bound_even} assures that
\begin{equation}\label{eqn:sharp:item2:upper_bound_partial}
\frac{2 \expconst^{\nicefrac{d}{2}}}{\Gamma\!\pa{\frac{d}{2}}} \int_{\frac{\sqrt{d(1+\beta)}}{\sqrt{2 \expconst}}}^{\frac{d \sqrt{1+\beta}}{\sqrt{2 \expconst}}}  e^{- \expconst r^2} r^{d-1}  \, dr \le d \! \pb{\frac{1+\beta}{e^{\beta}}}^{\nicefrac{d}{2}}\!.
\end{equation}

\noindent
This and \eqref{eqn:lemma:Gaussiantail_sharp_bounds_substitution_1+delta} imply that
\begin{equation}
\int_{\pc{y \in \R^d \colon \frac{\sqrt{d(1+\beta)}}{\sqrt{2 \expconst}} \le \norm{y}_2 \le \frac{d \sqrt{1+\beta}}{\sqrt{2 \expconst}}}} \! \pb{\frac{\expconst}{\pi}}^{\nicefrac{d}{2}} \! e^{- \expconst \norm{x}_2^2} \, dx \le d \! \pb{\frac{1+\beta}{e^{\beta}}}^{\nicefrac{d}{2}}\!.
\end{equation}
\end{cproof}



\subsection{Upper bounds for weighted Gaussian tails}
\label{subsec:bounds_weighted_Gaussian_tails}


\cfclear
\begin{lemma}\label{lemma:sharp_bounds_specific_integrals}
Let $d \in \N \cap [3, \infty)$, $\beta, \expconst \in (0, \infty)$, $k \in \N_0$ and let $\Gamma \colon (0, \infty) \to (0, \infty)$ satisfy for all $x \in (0, \infty)$ that $\Gamma(x) = \int_0^{\infty} t^{x-1} e^{-t} \, dt$. Then
\begin{equation}
\begin{split}
& \int_{\pc{y \in \R^d \colon \norm{y}_2 \ge \frac{\sqrt{d(1+\beta)}}{\sqrt{2 \expconst}}}} \! \pb{\frac{\expconst}{\pi}}^{\nicefrac{d}{2}} \! \norm{x}_2^k \, e^{- \expconst \norm{x}_2^2} \, dx \\
& \le d^{1+k} \! \pb{\frac{1+\beta}{2 \expconst}}^{\nicefrac{k}{2}} \! \pb{\frac{1+\beta}{e^{\beta}}}^{\nicefrac{d}{2}} + \frac{\Gamma\!\pa{\frac{d+k}{2}}}{\Gamma\!\pa{\frac{d}{2}}} \pb{\frac{d+k}{\expconst^{\nicefrac{k}{2}}}} e^{- \frac{d^2(1+\beta)}{2(d+k)}}\ifnocf.
\end{split}
\end{equation}
\cfout[.]
\end{lemma}
\begin{cproof}{lemma:sharp_bounds_specific_integrals}
\Nobs that \cref{lemma:Gaussiantail_sharp_bounds} (applied with $d \with d$, $\beta \with \beta$, $\expconst \with \expconst$ in the notation of \cref{lemma:Gaussiantail_sharp_bounds}) ensures that
\begin{equation}\label{eqn:lemma:sharp_bounds_specific_integrals_split_1}
\begin{split}
& \int_{\pc{y \in \R^d \colon \frac{\sqrt{d(1+\beta)}}{\sqrt{2 \expconst}} \le \norm{y}_2 \le \frac{d \sqrt{1+\beta}}{\sqrt{2 \expconst}}}} \! \pb{\frac{\expconst}{\pi}}^{\nicefrac{d}{2}} \! \norm{x}_2^k \, e^{- \expconst \norm{x}_2^2} \, dx \\
& \le d^k \! \pb{\frac{1+\beta}{2 \expconst}}^{\nicefrac{k}{2}} \! \int_{\pc{y \in \R^d \colon \frac{\sqrt{d(1+\beta)}}{\sqrt{2 \expconst}} \le \norm{y}_2 \le \frac{d \sqrt{1+\beta}}{\sqrt{2 \expconst}}}} \! \pb{\frac{\expconst}{\pi}}^{\nicefrac{d}{2}} \! e^{- \expconst \norm{x}_2^2} \, dx \\
& \le d^{1+k} \! \pb{\frac{1+\beta}{2 \expconst}}^{\nicefrac{k}{2}} \! \pb{\frac{1+\beta}{e^{\beta}}}^{\nicefrac{d}{2}}\ifnocf.
\end{split}
\end{equation}

\noindent
\cfload[.]Moreover, \nobs that \cref{lemma:transformation_Rd_to_R} (applied with $d \with d+k$, $\expconst \with \expconst$, $\alpha \with 0$, $s \with (2 \expconst)^{\nicefrac{-1}{2}} d(1+\beta)^{\nicefrac{1}{2}}$ in the notation of \cref{lemma:transformation_Rd_to_R}) and \cref{cor:Gaussiantail_alpha=2d} (applied with $d \with d+k$, $\expconst \with \expconst$, $s \with (2 \expconst)^{\nicefrac{-1}{2}} d (1+\beta)^{\nicefrac{1}{2}}$ in the notation of \cref{cor:Gaussiantail_alpha=2d}) assure that
\begin{equation}
\frac{2 \expconst^{\frac{d+k}{2}}}{\Gamma\!\pa{\frac{d+k}{2}}} \int_{\frac{d \sqrt{1+\beta}}{\sqrt{2 \expconst}}}^{\infty}  e^{- \expconst r^2} r^{d+k-1} \, dr = \int_{\pc{y \in \R^{d+k} \colon \norm{y}_2 \ge \frac{d \sqrt{1+\beta}}{\sqrt{2 \expconst}}}\!} \! \pb{\frac{\expconst}{\pi}}^{\frac{d+k}{2}} \! e^{- \expconst \norm{x}_2^2} \, dx \le (d+k) e^{-\frac{d^2(1+\beta)}{2(d+k)}}.
\end{equation}

\noindent
\cref{lemma:transformation_Rd_to_R} (applied with $d \with d$, $\expconst \with \expconst$, $\alpha \with k$, $s \with (2 \expconst)^{\nicefrac{-1}{2}} d(1+\beta)^{\nicefrac{1}{2}}$ in the notation of \cref{lemma:transformation_Rd_to_R}) hence shows that
\begin{equation}
\begin{split}
& \int_{\pc{y \in \R^d \colon \norm{y}_2 \ge \frac{d \sqrt{1+\beta}}{\sqrt{2 \expconst}}}} \! \pb{\frac{\expconst}{\pi}}^{\nicefrac{d}{2}} \! \norm{x}_2^k \, e^{- \expconst \norm{x}_2^2} \, dx = \frac{2 \expconst^{\nicefrac{d}{2}}}{\Gamma\!\pa{\frac{d}{2}}} \int_{\frac{d \sqrt{1+\beta}}{\sqrt{2 \expconst}}}^{\infty}  e^{- \expconst r^2} r^{d+k-1} \, dr \\
& = \frac{\Gamma\!\pa{\frac{d+k}{2}}}{\Gamma\!\pa{\frac{d}{2}} \expconst^{\nicefrac{k}{2}}} \! \pb{\frac{2 \expconst^{\frac{d+k}{2}}}{\Gamma\!\pa{\frac{d+k}{2}}} \int_{\frac{d \sqrt{1+\beta}}{\sqrt{2 \expconst}}}^{\infty}  e^{- \expconst r^2} r^{d+k-1} \, dr} \le \frac{\Gamma\!\pa{\frac{d+k}{2}}}{\Gamma\!\pa{\frac{d}{2}}} \pb{\frac{d+k}{\expconst^{\nicefrac{k}{2}}}} e^{- \frac{d^2(1+\beta)}{2(d+k)}}.
\end{split}
\end{equation}

\noindent
Combining this with \cref{eqn:lemma:sharp_bounds_specific_integrals_split_1} demonstrates that
\begin{equation}
\begin{split}
& \int_{\pc{y \in \R^d \colon \norm{y}_2 \ge \frac{\sqrt{d(1+\beta)}}{\sqrt{2 \expconst}}}} \! \pb{\frac{\expconst}{\pi}}^{\nicefrac{d}{2}} \! \norm{x}_2^k \, e^{- \expconst \norm{x}_2^2} \, dx \\
& = \int_{\pc{y \in \R^d \colon \frac{\sqrt{d(1+\beta)}}{\sqrt{2 \expconst}} \le \norm{y}_2 \le \frac{d \sqrt{1+\beta}}{\sqrt{2 \expconst}}}} \! \pb{\frac{\expconst}{\pi}}^{\nicefrac{d}{2}} \! \norm{x}_2^k \, e^{- \expconst \norm{x}_2^2} \, dx \\
& \quad + \int_{\pc{y \in \R^d \colon \norm{y}_2 \ge \frac{d \sqrt{1+\beta}}{\sqrt{2 \expconst}}}} \! \pb{\frac{\expconst}{\pi}}^{\nicefrac{d}{2}} \! \norm{x}_2^k \, e^{- \expconst \norm{x}_2^2} \, dx \\
& \le d^{1+k} \! \pb{\frac{1+\beta}{2 \expconst}}^{\nicefrac{k}{2}} \! \pb{\frac{1+\beta}{e^{\beta}}}^{\nicefrac{d}{2}}\! + \frac{\Gamma\!\pa{\frac{d+k}{2}}}{\Gamma\!\pa{\frac{d}{2}}} \pb{\frac{d+k}{\expconst^{\nicefrac{k}{2}}}} e^{- \frac{d^2(1+\beta)}{2(d+k)}}.
\end{split}
\end{equation}
\end{cproof}




\section{Lower bounds for the number of ANN parameters in the approximation of high-dimensional functions}\label{sec:lower_bounds_for_number_of_parameters_in_ANN_approximations}

In this section we employ the upper bounds for certain weighted tails of standard normal distributions from \cref{sec:upper_bounds_for_weighted_Gaussian_tails} to establish in \cref{cor:main_2} in Subsection~\ref{subsec:lower_bounds_anns_for_specific_class} below suitable lower bounds for the number of parameters of appropriate ANNs that approximate certain high-dimensional target functions. Our proof of \cref{cor:main_2} employs appropriate lower bounds for the product of the number of ANN parameters and the maximum of the absolute values of the ANN parameters which we establish in \cref{cor:main} in Subsection~\ref{subsec:lower_bounds_anns_for_specific_class} below. Our proof of \cref{cor:main}, in turn, employs the lower bounds for general ANNs in \cref{thm:main1}. Our proof of \cref{thm:main1} uses the elementary lower bounds for normalized $L^2$-scalar products in \cref{lemma:CauchySchwarz_like} in Subsection~\ref{subsec:connection_distance&scalar_prod} below as well as the upper bounds for $L^2$-scalar products involving realizations of ANNs in \cref{lemma:upper_bounds_for_the_scalar_product} in Subsection~\ref{subsec:upper_bounds_scalar_prod} below. Our proof of \cref{lemma:upper_bounds_for_the_scalar_product} employs the priori estimates for realizations of ANNs in \cref{lemma:on_Realization_of_DNN}, \cref{cor:on_Realization_of_DNN}, and \cref{lemma:on_the_norm_of_realization} in Subsection~\ref{subsec:upper_bounds_for_realizations_anns} below. Our proofs of \cref{lemma:on_Realization_of_DNN} and \cref{cor:on_Realization_of_DNN} use the well-known matrix norm estimates in \cref{lemma:on_matrix_norm_inequalities} below. Only for the sake of completeness we include in this section also the detailed proofs for \cref{lemma:on_matrix_norm_inequalities} and \cref{lemma:CauchySchwarz_like}.

\subsection{Upper bounds for realizations of ANNs}
\label{subsec:upper_bounds_for_realizations_anns}


\cfclear
\begin{lemma}\label{lemma:on_matrix_norm_inequalities}
Let $m, n \in \N$, $A = (A_{i, j})_{(i, j) \in \{1, 2, \ldots, m\} \times \{1, 2, \ldots, n\}} \in \R^{m \times n}$, $B = (B_1, B_2, \ldots, B_m) \in \R^m$, $x \in \R^n$ \cfload. Then
\begin{enumerate}[label=(\roman *)]
\item
\label{item1:lemma:on_matrix_norm_inequalities} it holds that
\begin{equation}
\norm{Ax + B}_{\infty} \le \sqrt{n} \Big[\max_{i \in \{1, 2, \ldots, m\}} \max_{j \in \{1, 2, \ldots, n\}} \abs{A_{i, j}}\Big] \norm{x}_2 +  \norm{B}_{\infty}
\end{equation}

and
\item
\label{item2:lemma:on_matrix_norm_inequalities} it holds that
\begin{equation}
\norm{Ax + B}_{\infty} \le n \Big[\max_{i \in \{1, 2, \ldots, m\}} \max_{j \in \{1, 2, \ldots, n\}} \abs{A_{i, j}}\Big] \norm{x}_{\infty} + \norm{B}_{\infty}\ifnocf.
\end{equation}
\end{enumerate}
\cfout[.]
\end{lemma}
\begin{cproof}{lemma:on_matrix_norm_inequalities}
Throughout this proof let $\alpha \in \R$ satisfy $\alpha = \max_{i \in \{1, 2, \ldots, m\}} \allowbreak \max_{j \in \{1, 2, \ldots, n\}} \allowbreak \abs{A_{i, j}}$ and let $\beta \in \R$ satisfy $\beta = \norm{B}_{\infty}$ \cfload[.]\Nobs that the triangle inequality and the fact that for all $v = (v_1, v_2, \ldots, v_n) \in \R^n$ it holds that $\sum_{j=1}^n \abs{v_j} \le \sqrt{n} \norm{v}_2$ ensure that for all $v = (v_1, v_2, \ldots, v_n) \in \R^n$ it holds that
\begin{equation}
\begin{split}
\norm{Av + B}_{\infty} & = \max_{i \in \{1, 2, \ldots, m\}} \! \bigg|B_i + \sum_{j=1}^n A_{i, j} v_j \bigg| \le \max_{i \in \{1, 2, \ldots, m\}} \! \bigg(\abs{B_i} + \sum_{j=1}^n \abs{A_{i, j} v_j} \bigg) \\
& \le \beta + \alpha \sum_{j=1}^n \abs{v_j} \le \beta + \alpha \sqrt{n} \norm{v}_2.
\end{split}
\end{equation}

\noindent
This establishes \cref{item1:lemma:on_matrix_norm_inequalities}. Moreover, \nobs that the fact that for all $v \in \R^n$ it holds that $\sqrt{n} \norm{v}_2 \le n \infnorm{v}$ and \cref{item1:lemma:on_matrix_norm_inequalities} demonstrate that for all $v \in \R^n$ it holds that
\begin{equation}
\norm{Av + B}_{\infty} \le \beta + \alpha \sqrt{n} \norm{v}_2 \le \beta + \alpha n \infnorm{v}.
\end{equation}

\noindent
This establishes \cref{item2:lemma:on_matrix_norm_inequalities}.
\end{cproof}

\cfclear
\begin{lemma}\label{lemma:on_Realization_of_DNN}
Let $L \in \N \cap [2, \infty)$, $l_0, l_1, \ldots, l_L \in \N$, $\Phi = ((W_1, B_1), (W_2, B_2), \ldots, \allowbreak (W_L, B_L)) \in \big(\! \bigtimes_{k=1}^{L} (\R^{l_k \times l_{k-1}} \times \R^{l_k})\big)$, $x_0 \in \R^{l_0}, x_1 \in \R^{l_1}, \ldots, x_L \in \R^{l_L}$ satisfy for all $k \in \{1, 2, \ldots, L\}$ that $x_k = \Rect (W_k x_{k-1}+B_k)$ \cfload. Then
\begin{enumerate}[label=(\roman *)]
\item\label{item:lemma:on_Realization_of_DNN_1}
it holds for all $k \in \{1, 2, \ldots, L\}$, $j \in \pc{1, 2, \ldots, k}$ that
\begin{equation}\label{eqn:lemma:on_Realization_of_DNN_1}
\norm{x_k}_{\infty} \le l_{k-1}l_{k-2} \cdots l_{k-j} \pa{\max\!\pc{1, \norm{\vectorNN(\Phi)}_{\infty}}}^{j} \! (\norm{x_{k-j}}_{\infty}+j)
\end{equation}

and
\item\label{item:lemma:on_Realization_of_DNN_2}
it holds that
\begin{equation}
\norm{(\functionANN(\Phi))(x_0)}_{\infty} \le l_{L-1} l_{L-2} \cdots l_1 \! \pa{\max\!\pc{1, \norm{\vectorNN(\Phi)}_{\infty}}}^{L-1} \! (\norm{x_1}_{\infty}+L-1)\ifnocf.
\end{equation}

\end{enumerate}
\cfout[.]
\end{lemma}
\begin{cproof}{lemma:on_Realization_of_DNN}
Throughout this proof let $\alpha = \max \! \pc{1, \norm{\vectorNN(\Phi)}_{\infty}}$ \cfload. \Nobs that the fact that for all $x \in \R$ it holds that $\abs{\max\{x, 0\}} \le \abs{x}$ and \cref{item2:lemma:on_matrix_norm_inequalities} in \cref{lemma:on_matrix_norm_inequalities} (applied for every $k \in \{1, 2, \ldots, L\}$ with $m \with l_k$, $n \with l_{k-1}$, $A \with W_k$, $B \with B_k$, $x \with x_{k-1}$ in the notation of \cref{lemma:on_matrix_norm_inequalities}) imply that for all $k \in \{1, 2, \ldots, L\}$ it holds that
\begin{equation}\label{eqn:lemma:on_Realization_of_DNN_induction_1}
\begin{split}
\norm{x_k}_{\infty} & = \norm{\Rect (W_k x_{k-1}+B_k)}_{\infty} \le \norm{W_k x_{k-1}+B_k}_{\infty} \\
& \le \alpha \, l_{k-1}\norm{x_{k-1}}_{\infty}+\alpha \le \alpha \, l_{k-1} (\norm{x_{k-1}}_{\infty}+1).
\end{split}
\end{equation}

\noindent
This demonstrates that for all $k \in \{2, 3, \ldots, L\}$, $i \in \pc{1, 2, \ldots, k-1}$ with $\norm{x_k}_{\infty} \le l_{k-1} \allowbreak l_{k-2} \allowbreak \cdots \allowbreak l_{k-i} \allowbreak \alpha^{i}(\norm{x_{k-i}}_{\infty}+i)$ it holds that
\begin{equation}
\begin{split}
\norm{x_k}_{\infty} & \le l_{k-1}l_{k-2} \cdots l_{k-i} \alpha^{i}(\norm{x_{k-i}}_{\infty}+i) \\
& \le l_{k-1}l_{k-2} \cdots l_{k-i} \alpha^{i}(\alpha l_{k-i-1}(\norm{x_{k-i-1}}_{\infty} + 1)+i) \\
& \le l_{k-1}l_{k-2} \cdots l_{k-i} l_{k-i-1} \alpha^{i+1}(\norm{x_{k-i-1}}_{\infty}+i+1).
\end{split}
\end{equation}

\noindent
This, \cref{eqn:lemma:on_Realization_of_DNN_induction_1}, and induction show that for all $k \in \{1, 2, \ldots, L\}$, $j \in \pc{1, 2, \ldots, k}$ it holds that
\begin{equation}
\norm{x_k}_{\infty} \le l_{k-1}l_{k-2} \cdots l_{k-j} \pa{\max\!\pc{1, \norm{\vectorNN(\Phi)}_{\infty}}}^{j}\!(\norm{x_{k-j}}_{\infty}+j).
\end{equation}

\noindent
This establishes \cref{item:lemma:on_Realization_of_DNN_1}. Next \nobs that \cref{item2:lemma:on_matrix_norm_inequalities} in \cref{lemma:on_matrix_norm_inequalities} (applied with $m \with l_L$, $n \with l_{L-1}$, $A \with W_L$, $B \with B_L$, $x \with x_{L-1}$ in the notation of \cref{lemma:on_matrix_norm_inequalities}) ensures that
\begin{equation}
\begin{split}
\norm{(\functionANN(\Phi))(x_0)}_{\infty} = \norm{W_L x_{L-1} + B_L}_{\infty} \le \alpha l_{L-1} \norm{x_{L-1}}_{\infty}+\alpha \le \alpha l_{L-1} (\norm{x_{L-1}}_{\infty}+1)\ifnocf.
\end{split}
\end{equation}

\noindent
\cfload[.]This and \cref{item:lemma:on_Realization_of_DNN_1} demonstrate that
\begin{equation}
\begin{split}
\norm{(\functionANN(\Phi))(x_0)}_{\infty} & \le \alpha l_{L-1} (\norm{x_{L-1}}_{\infty}+1) \\
& \le \alpha l_{L-1} ([l_{L-2} l_{L-3} \cdots l_1 \alpha^{L-2} (\norm{x_1}_{\infty}+L-2)]+1) \\
& \le l_{L-1}l_{L-2} \cdots l_1 \alpha^{L-1} (\norm{x_1}_{\infty}+L-1).
\end{split}
\end{equation}

\noindent
This establishes \cref{item:lemma:on_Realization_of_DNN_2}.
\end{cproof}


\cfclear
\begin{corollary}\label{cor:on_Realization_of_DNN}
It holds for all $\Phi \in \ANNs$, $x \in \R^{\inDimANN(\Phi)}$ that
\begin{equation}\label{eqn:cor:on_Realization_of_DNN_0}
\norm{(\functionANN(\Phi))(x)}_{\infty} \le \! \pb{\frac{\paramANN(\Phi) \max\!\pc{1, \norm{\vectorNN(\Phi)}_{\infty}}}{2 \lengthANN(\Phi)}}^{\lengthANN(\Phi)} \! (\norm{x}_2+\lengthANN(\Phi))\ifnocf.
\end{equation}
\cfout[.]
\end{corollary}
\begin{cproof}{cor:on_Realization_of_DNN}
Throughout this proof let $L \in \N$, $l_0, l_1, \ldots, l_L \in \N$,  $\Phi =((W_1, B_1), \allowbreak (W_2, B_2), \ldots, (W_L, B_L)) \in (\bigtimes_{k=1}^{L} (\R^{l_k \times l_{k-1}} \times \R^{l_k}))$, $\alpha = \max \! \pc{1, \norm{\vectorNN(\Phi)}_{\infty}}$, $x_0 \in \R^{l_0}$, $x_1 \in \R^{l_1}$ satisfy $x_1 = \Rect (W_1 x_0 +B_1)$ \cfload. \Nobs that \cref{item1:lemma:on_matrix_norm_inequalities} in \cref{lemma:on_matrix_norm_inequalities} (applied with $m \with l_1$, $n \with l_0$, $A \with W_1$, $B \with B_1$, $x \with x_0$ in the notation of \cref{lemma:on_matrix_norm_inequalities}) ensures that
\begin{equation}\label{eqn:cor:on_Realization_of_DNN:lemma_apply}
\norm{W_1 x_0 + B_1}_{\infty} \le \alpha \sqrt{l_0} \norm{x_0}_2 + \alpha \le \alpha \sqrt{l_0}(\norm{x_0}_2 + 1).
\end{equation}

\noindent
In the following we distinguish between the case $\lengthANN(\Phi) = 1$ and the case $\lengthANN(\Phi) > 1$. We first prove \cref{eqn:cor:on_Realization_of_DNN_0} in the case $\lengthANN(\Phi) = 1$. \Nobs that \cref{eqn:cor:on_Realization_of_DNN:lemma_apply} demonstrates that
\begin{equation}
\begin{split}
\norm{(\functionANN(\Phi))(x_0)}_{\infty} & = \norm{W_1 x_0 + B_1}_{\infty} \le \alpha \sqrt{l_0} (\norm{x_0}_2 + 1) \\
& \le \frac{(l_0 + 1) \alpha}{2} \, (\norm{x_0}_2 + 1) \le \frac{l_1 (l_0 + 1) \alpha}{2} \, (\norm{x_0}_2 + 1) \\
&  = \! \pb{\frac{\paramANN(\Phi) \max\!\pc{1, \norm{\vectorNN(\Phi)}_{\infty}}}{2 \lengthANN(\Phi)}}^{\lengthANN(\Phi)} \! (\norm{x_0}_2+\lengthANN(\Phi)).
\end{split}
\end{equation}

\noindent
This proves \cref{eqn:cor:on_Realization_of_DNN_0} in case $\lengthANN(\Phi) = 1$. We now prove \cref{eqn:cor:on_Realization_of_DNN_0} in the case $\lengthANN(\Phi) > 1$. \Nobs that \cref{eqn:cor:on_Realization_of_DNN:lemma_apply} and the fact that for all $x \in \R$ it holds that $\abs{\max\{x, 0\}} \le \abs{x}$ show that
\begin{equation}
\norm{x_1}_{\infty} = \norm{\Rect (W_1 x_0 +B_1)}_{\infty} \le \norm{W_1 x_0 +B_1}_{\infty} \le \alpha \sqrt{l_0} (\norm{x_0}_2 + 1).
\end{equation}

\noindent
This and \cref{item:lemma:on_Realization_of_DNN_2} in \cref{lemma:on_Realization_of_DNN} (applied with $L \with L$, $l_0 \with l_0$, $l_1 \with l_1$, \ldots, $l_L \with l_L$, $\Phi \with \Phi$, $x_0 \with x_0$, $x_1 \with x_1$ in the notation of \cref{lemma:on_Realization_of_DNN}) ensure that
\begin{equation}\label{eqn:cor:on_Realization_of_DNN_2}
\begin{split}
\norm{(\functionANN(\Phi))(x_0)}_{\infty} & \le l_{L-1}l_{L-2} \cdots l_1 \alpha^{L-1} (\norm{x_1}_{\infty}+L-1) \\
& \le l_{L-1}l_{L-2} \cdots l_1 \alpha^{L-1} (\alpha \sqrt{l_0} (\norm{x_0}_2 + 1) +L-1) \\
& \le l_{L-1}l_{L-2} \cdots l_1 \sqrt{l_0} \, \alpha^L (\norm{x_0}_2 + L).
\end{split}
\end{equation}

\noindent
In the next step \nobs that the inequality of arithmetic and geometric means assures that
\begin{equation}
\begin{split}
\paramANN(\Phi) & = \sum_{k=1}^{L} l_k (l_{k-1} + 1) = l_1 + l_2 + \ldots + l_L + l_0 l_1 +l_1 l_2 + \ldots +l_{L-1} l_L \\
& \ge 2L \! \pb{(l_1 l_2 \cdots l_L) (l_0 l_1 l_1 l_2 \cdots l_{L-1} l_L)}^{\nicefrac{1}{2L}} \! = 2L \! \pb{l_0 (l_1)^3 (l_2)^3 \cdots (l_{L-1})^3 (l_L)^2}^{\nicefrac{1}{2L}} \\
& \ge 2L \! \pb{l_0 (l_1)^2 (l_2)^2 \cdots (l_{L-1})^2}^{\nicefrac{1}{2L}}\!.
\end{split}
\end{equation}

\noindent
Hence, we obtain that
\begin{equation}
l_{L-1} l_{L-2} \cdots l_1 \sqrt{l_0} \le \! \pb{\frac{\paramANN(\Phi)}{2L}}^L \!.
\end{equation}

\noindent
Combining this and \cref{eqn:cor:on_Realization_of_DNN_2} shows that
\begin{equation}
\begin{split}
\norm{(\functionANN(\Phi))(x_0)}_{\infty} & \le l_{L-1}l_{L-2} \cdots l_1 \sqrt{l_0} \, \alpha^L (\norm{x_0}_2 + L) \\
& \le \! \pb{\frac{\paramANN(\Phi)\alpha}{2L}}^L \! (\norm{x_0}_2 + L) \\
& = \! \pb{\frac{\paramANN(\Phi) \max\!\pc{1, \norm{\vectorNN(\Phi)}_{\infty}}}{2 \lengthANN(\Phi)}}^{\lengthANN(\Phi)} \! (\norm{x_0}_2+\lengthANN(\Phi)).
\end{split}
\end{equation}

\noindent
This proves \cref{eqn:cor:on_Realization_of_DNN_0} in the case $\lengthANN(\Phi) >1 $.
\end{cproof}

\cfclear
\begin{lemma}\label{lemma:on_the_norm_of_realization}
Let $d \in \N \cap [4, \infty)$, $\beta, \expconst \in (0, \infty)$, $\Phi \in \ANNs$ satisfy $\inDimANN(\Phi)=d$ and $\outDimANN(\Phi)=1$ and let $\varphi \colon \R^d \to \R$ satisfy for all $x \in \R^d$ that $\varphi(x) = (\nicefrac{\expconst}{\pi})^{\nicefrac{d}{2}} \exp(- \expconst \norm{x}_2^2)$ \cfload. Then
\begin{equation}
\begin{split}
& \int_{\pc{y \in \R^d \colon \norm{y}_2 \ge \frac{\sqrt{d(1+\beta)}}{\sqrt{2 \expconst}}}} \pabs{(\functionANN(\Phi))(x)}^2 \! \varphi(x) \, dx \\
& \quad \quad \quad \le \abs{\lengthANN(\Phi)}^2 \! \pb{\frac{\paramANN(\Phi) \max\!\pc{1, \norm{\vectorNN(\Phi)}_{\infty}}}{2 \lengthANN(\Phi)}}^{2\lengthANN(\Phi)} \! \pb{\frac{1+\beta}{e^{\beta}}}^{\nicefrac{d}{3}}\! \pb{\frac{d^3(6+4\beta+\expconst)}{4 \expconst}}\!\ifnocf.
\end{split}
\end{equation}
\cfout[.]
\end{lemma}
\begin{cproof}{lemma:on_the_norm_of_realization}
Throughout this proof let $\scrR \in \R$ satisfy $\sqrt{2 \expconst} \scrR = \sqrt{d(1+\beta)}$. \Nobs that the fact that for all $\fa, \fb \in \R$ it holds that $(\fa+\fb)^2 \le 2(\fa^2+\fb^2)$ and \cref{cor:on_Realization_of_DNN} imply that for all $x \in \R^d$ it holds that
\begin{equation}\label{eqn:lemma:on_the_norm_of_realization}
\begin{split}
\abs{(\functionANN(\Phi))(x)}^2 & = \norm{(\functionANN(\Phi))(x)}_{\infty}^2 \le \! \pb{\frac{\paramANN(\Phi) \max\!\pc{1, \norm{\vectorNN(\Phi)}_{\infty}}}{2 \lengthANN(\Phi)}}^{2\lengthANN(\Phi)} \! (\norm{x}_2+\lengthANN(\Phi))^2 \\
& \le 2 \! \pb{\frac{\paramANN(\Phi) \max\!\pc{1, \norm{\vectorNN(\Phi)}_{\infty}}}{2 \lengthANN(\Phi)}}^{2\lengthANN(\Phi)} \! (\norm{x}_2^2+\abs{\lengthANN(\Phi)}^2)\ifnocf.
\end{split}
\end{equation}

\noindent
\cfload[.]\Nobs that \cref{item:lemma:GammaBeta1} in \cref{lemma:GammaBeta} ensures that $\Gamma\!\pa{\nicefrac{d}{2}+1} = [\nicefrac{d}{2}] \Gamma\!\pa{\nicefrac{d}{2}}$. Combining this with \cref{lemma:sharp_bounds_specific_integrals} (applied with $d \with d$, $\beta \with \beta$, $\expconst \with \expconst$, $k \with 0$, $k \with 2$ in the notation of \cref{lemma:sharp_bounds_specific_integrals}) assures that
\begin{equation}
\begin{split}
& \int_{\{y \in \R^d \colon \norm{y}_2 \ge \scrR\}} \big[\abs{\lengthANN(\Phi)}^2+\norm{x}_2^2\big] \varphi(x) \, dx \\
& = \abs{\lengthANN(\Phi)}^2 \int_{\{y \in \R^d \colon \norm{y}_2 \ge \scrR\}} \varphi(x) \, dx + \int_{\{y \in \R^d \colon \norm{y}_2 \ge \scrR \}} \norm{x}_2^2 \, \varphi(x) \, dx \\
& = \abs{\lengthANN(\Phi)}^2 \int_{\{y \in \R^d \colon \norm{y}_2 \ge \scrR \}} \! \pb{\frac{\expconst}{\pi}}^{\nicefrac{d}{2}} \! e^{- \expconst \norm{x}_2^2} \, dx + \int_{\{y \in \R^d \colon \norm{y}_2 \ge \scrR \}} \! \pb{\frac{\expconst}{\pi}}^{\nicefrac{d}{2}} \! \norm{x}_2^2 \, e^{- \expconst \norm{x}_2^2} \, dx\\
& \le \abs{\lengthANN(\Phi)}^2 \! \pb{d \! \pb{\frac{1+\beta}{e^{\beta}}}^{\nicefrac{d}{2}} + de^{-\frac{d(1+\beta)}{2}}} \! + \! \pb{\frac{d^3(1+\beta)}{2 \expconst} \! \pb{\frac{1+\beta}{e^{\beta}}}^{\nicefrac{d}{2}} \! + \frac{\Gamma\!\pa{\frac{d}{2}+1}}{\Gamma\!\pa{\frac{d}{2}}} \pb{\frac{d+2}{\expconst}}  e^{-\frac{d^2(1+\beta)}{2(d+2)}}} \\
& = \abs{\lengthANN(\Phi)}^2 \! \pb{d \! \pb{\frac{1+\beta}{e^{\beta}}}^{\nicefrac{d}{2}} + de^{-\frac{d(1+\beta)}{2}}} \! + \! \pb{\frac{d^3(1+\beta)}{2\expconst} \! \pb{\frac{1+\beta}{e^{\beta}}}^{\nicefrac{d}{2}}\! + \pb{\frac{d(d+2)}{2 \expconst}} e^{-\frac{d^2(1+\beta)}{2(d+2)}}} \\
& \le \abs{\lengthANN(\Phi)}^2 \! \pb{d \! \pb{\frac{1+\beta}{e^{\beta}}}^{\nicefrac{d}{2}} + de^{-\frac{d(1+\beta)}{2}} + \pb{\frac{d^3(1+\beta)}{2\expconst}} \pb{\frac{1+\beta}{e^{\beta}}}^{\nicefrac{d}{2}}\! + \pb{\frac{d(d+2)}{2 \expconst}}  e^{-\frac{d^2(1+\beta)}{2(d+2)}}}\!.
\end{split}
\end{equation}

\noindent
This, the fact that $d^3(1+\beta)+d(d+2)+4d \expconst \le d^3\big(\frac{3}{2}+\beta+\frac{\expconst}{4}\big)$, and the fact that
\begin{equation}
\max\!\pc{e^{\frac{-d(1+\beta)}{2}}, e^{- \frac{d^2(1+\beta)}{2(d+2)}}, \! \pb{\frac{1+\beta}{e^{\beta}}}^{\nicefrac{d}{2}}} \! \le \! \pb{\frac{1+\beta}{e^{\beta}}}^{\nicefrac{d}{3}}\!
\end{equation}

\noindent
imply that
\begin{equation}
\begin{split}
\int_{\{y \in \R^d \colon \norm{y}_2 \ge \scrR\}} \big[\abs{\lengthANN(\Phi)}^2+\norm{x}_2^2\big] \varphi(x) \, dx & \le \abs{\lengthANN(\Phi)}^2 \! \pb{\frac{1+\beta}{e^{\beta}}}^{\nicefrac{d}{3}}\! \pb{\frac{d^3(1+\beta)+d(d+2)+4d \expconst}{2 \expconst}} \\
& \le \abs{\lengthANN(\Phi)}^2 \! \pb{\frac{1+\beta}{e^{\beta}}}^{\nicefrac{d}{3}} \! \pb{\frac{d^3(6+4\beta+\expconst)}{8 \expconst}}\!.
\end{split}
\end{equation}

\noindent
Combining this with \cref{eqn:lemma:on_the_norm_of_realization} demonstrates that
\begin{equation}
\begin{split}
& \int_{\pc{y \in \R^d \colon \norm{y}_2 \ge \frac{\sqrt{d(1+\beta)}}{\sqrt{2 \expconst}}}} \abs{(\functionANN(\Phi))(x)}^2 \varphi(x) \, dx = \int_{\{y \in \R^d \colon \norm{y}_2 \ge \scrR\}} \abs{(\functionANN(\Phi))(x)}^2 \varphi(x) \, dx \\
& \le 2 \! \pb{\frac{\paramANN(\Phi) \max\!\pc{1, \norm{\vectorNN(\Phi)}_{\infty}}}{2 \lengthANN(\Phi)}}^{2\lengthANN(\Phi)} \! \pb{\int_{\{y \in \R^d \colon \norm{y}_2 \ge \scrR\}} \big[\abs{\lengthANN(\Phi)}^2+\norm{x}_2^2\big] \varphi(x) \, dx} \\
& \le \abs{\lengthANN(\Phi)}^2 \! \pb{\frac{\paramANN(\Phi) \max\!\pc{1, \norm{\vectorNN(\Phi)}_{\infty}}}{2 \lengthANN(\Phi)}}^{2\lengthANN(\Phi)} \! \pb{\frac{1+\beta}{e^{\beta}}}^{\nicefrac{d}{3}}\! \pb{\frac{d^3(6+4\beta+\expconst)}{4 \expconst}}\!.
\end{split}
\end{equation}
\end{cproof}



\subsection{Upper bounds for scalar products involving realizations of ANNs}
\label{subsec:upper_bounds_scalar_prod}


\cfclear
\begin{lemma}\label{lemma:upper_bounds_for_the_scalar_product}
 Let $d \in \N \cap [4, \infty)$, $\beta, \expconst \in (0, \infty)$, $\Phi \in \ANNs$ satisfy $\inDimANN(\Phi)=d$ and $\outDimANN(\Phi)=1$, let $\varphi \colon \R^d \to \R$, $\ff \colon \R^d \to \R$, and $\fg \colon \R^d \to \R$ be measurable, and assume for all $x \in \R^d$ that $\varphi(x) = (\nicefrac{\expconst}{\pi})^{\nicefrac{d}{2}} \exp(- \expconst \norm{x}_2^2)$, $\int_{\R^d} \abs{(\functionANN(\Phi))(y)} \, dy > 0$, $\int_{\R^d} \abs{\fg(y)}^2 \, dy = 1$, and
 \begin{equation}
\ff(x) = \textstyle \! \pb{\int_{\R^d} \abs{(\functionANN(\Phi))(y)}^2 \varphi(y) \, dy}^{\nicefrac{-1}{2}} \! \displaystyle (\functionANN(\Phi))(x)\pb{\varphi(x)}^{\nicefrac{1}{2}}
\end{equation}

\noindent
\cfload. Then
\begin{equation}
\begin{split}
& \int_{\R^d} \abs{\ff(x) \fg(x)} \, dx \le \!\pb{\int_{\pc{y \in \R^d \colon \norm{y}_2 \le \frac{\sqrt{d(1+\beta)}}{\sqrt{2 \expconst}}}} \abs{\fg(x)}^2 \, dx}^{\nicefrac{1}{2}} \\
&\quad + \lengthANN(\Phi) \! \pb{\frac{\paramANN(\Phi) \max\!\pc{1, \norm{\vectorNN(\Phi)}_{\infty}}}{2 \lengthANN(\Phi)}}^{\lengthANN(\Phi)} \! \pb{\frac{1+\beta}{e^{\beta}}}^{\nicefrac{d}{6}}\! \pb{\frac{d^{\nicefrac{3}{2}}(6+4\beta+\expconst)^{\nicefrac{1}{2}}}{2\sqrt{\expconst}\pb{\int_{\R^d} \abs{(\functionANN(\Phi))(y)}^2 \varphi(y) \, dy}^{\nicefrac{1}{2}}}}\!\ifnocf.
\end{split}
\end{equation}
\cfout[.]
\end{lemma}
\begin{cproof}{lemma:upper_bounds_for_the_scalar_product} Throughout this proof let $\Gamma \colon (0, \infty) \to (0, \infty)$ satisfy for all $x \in (0, \infty)$ that $\Gamma(x) = \int_0^{\infty} t^{x-1} e^{-t} \, dt$, let $\fa \in \R$ satisfy $\fa = \! \pb{\int_{\R^d} \abs{(\functionANN(\Phi))(y)}^2 \varphi(y) \, dy}^{\nicefrac{1}{2}}$, and let $\scrR \in \R$ satisfy $\sqrt{2 \expconst} \scrR = \sqrt{d(1+\beta)}$. \Nobs that $\fa \in (0, \infty)$ and
\begin{equation}
\int_{\R^d} \abs{\ff(x)}^2 \, dx = \! \pb{\int_{\R^d} \abs{(\functionANN(\Phi))(y)}^2 \varphi(y) \, dy}^{-1} \! \int_{\R^d} \abs{(\functionANN(\Phi))(x)}^2 \varphi(x) \, dx = 1.
\end{equation}

\noindent
Combining this with the H\"older inequality shows that
\begin{equation}\label{eqn:lemma:upper_bounds_Schwarz}
\begin{split}
\int_{\R^d} \abs{\ff(x) \fg(x)} \, dx & = \int_{\{y \in \R^d \colon \norm{y}_2 \le \scrR\}} \abs{\ff(x)\fg(x)} \, dx + \int_{\{y \in \R^d \colon \norm{y}_2 \ge \scrR\}} \abs{\ff(x)\fg(x)} \, dx \\
& \le \pb{\int_{\{y \in \R^d \colon \norm{y}_2 \le \scrR\}} \abs{\ff(x)}^2 \, dx}^{\nicefrac{1}{2}} \! \pb{\int_{\{y \in \R^d \colon \norm{y}_2 \le \scrR\}} \abs{\fg(x)}^2 \, dx}^{\nicefrac{1}{2}}  \\
& \quad + \pb{\int_{\{y \in \R^d \colon \norm{y}_2 \ge \scrR\}} \abs{\ff(x)}^2 \, dx}^{\nicefrac{1}{2}} \! \pb{\int_{\{y \in \R^d \colon \norm{y}_2 \ge \scrR\}} \abs{\fg(x)}^2 \, dx}^{\nicefrac{1}{2}} \\
& \le \pb{\int_{\R^d} \abs{\ff(x)}^2 \, dx}^{\nicefrac{1}{2}} \! \pb{\int_{\{y \in \R^d \colon \norm{y}_2 \le \scrR\}} \abs{\fg(x)}^2 \, dx}^{\nicefrac{1}{2}}  \\
& \quad + \pb{\int_{\{y \in \R^d \colon \norm{y}_2 \ge \scrR\}} \abs{\ff(x)}^2 \, dx}^{\nicefrac{1}{2}} \! \pb{\int_{\R^d} \abs{\fg(x)}^2 \, dx}^{\nicefrac{1}{2}} \\
& = \pb{\int_{\{y \in \R^d \colon \norm{y}_2 \le \scrR\}} \abs{\fg(x)}^2 \, dx}^{\nicefrac{1}{2}} \! + \pb{\int_{\{y \in \R^d \colon \norm{y}_2 \ge \scrR\}} \abs{\ff(x)}^2 \, dx}^{\nicefrac{1}{2}}\!.
\end{split}
\end{equation}

\noindent
Next we obtain that \cref{lemma:on_the_norm_of_realization} (applied with $d \with d$, $\beta \with \beta$, $\expconst \with \expconst$, $\Phi \with \Phi$, $\varphi \with \varphi$ in the notation of \cref{lemma:on_the_norm_of_realization}) implies that
\begin{equation}
\begin{split}
& \int_{\{y \in \R^d \colon \norm{y}_2 \ge \scrR\}}  \abs{\ff(x)}^2 \, dx = \fa^{-2} \int_{\{y \in \R^d \colon \norm{y}_2 \ge \scrR\}} \abs{(\functionANN(\Phi))(x)}^2 \varphi(x) \, dx \\
& \le \fa^{-2} \abs{\lengthANN(\Phi)}^2 \! \pb{\frac{\paramANN(\Phi) \max\!\pc{1, \norm{\vectorNN(\Phi)}_{\infty}}}{2 \lengthANN(\Phi)}}^{2\lengthANN(\Phi)} \! \pb{\frac{1+\beta}{e^{\beta}}}^{\nicefrac{d}{3}}\! \pb{\frac{d^3(6+4\beta+\expconst)}{4 \expconst}}\!.
\end{split}
\end{equation}

\noindent
This and \cref{eqn:lemma:upper_bounds_Schwarz} imply that
\begin{equation}
\begin{split}
& \int_{\R^d} \abs{\ff(x) \fg(x)} \, dx \le \pb{\int_{\{y \in \R^d \colon \norm{y}_2 \le \scrR\}} \abs{\fg(x)}^2 \, dx}^{\nicefrac{1}{2}} \!+ \! \pb{\int_{\{y \in \R^d \colon \norm{y}_2 \ge \scrR\}} \abs{f(x)}^2 \, dx}^{\nicefrac{1}{2}} \\
& \le \pb{\int_{\{y \in \R^d \colon \norm{y}_2 \le \scrR\}} \abs{\fg(x)}^2 \, dx}^{\nicefrac{1}{2}} \\
& \quad + \fa^{-1} \lengthANN(\Phi) \! \pb{\frac{\paramANN(\Phi) \max\!\pc{1, \norm{\vectorNN(\Phi)}_{\infty}}}{2 \lengthANN(\Phi)}}^{\lengthANN(\Phi)} \! \pb{\frac{1+\beta}{e^{\beta}}}^{\nicefrac{d}{6}}\! \pb{\frac{d^{\nicefrac{3}{2}}(6+4\beta+\expconst)^{\nicefrac{1}{2}}}{2\sqrt{\expconst}}}\!.
\end{split}
\end{equation}
\end{cproof}




\subsection{On the connection of distances and scalar products}
\label{subsec:connection_distance&scalar_prod}


\begin{lemma}\label{lemma:CauchySchwarz_like}
Let $d \in \N$, $\alpha \in \R$, let $\ff \colon \R^d \to \R$ and $\fg \colon \R^d \to \R$ be measurable, and assume $\int_{\R^d} \abs{\ff(x)}^2 \, dx=\int_{\R^d} \abs{\fg(x)}^2 \, dx=1$. Then
\begin{equation}
\int_{\R^d} \abs{\alpha \ff(x) - \fg(x)}^2 \, dx \ge 1 - \int_{\R^d} \abs{\ff(x) \fg(x)} \, dx.
\end{equation}
\end{lemma}
\begin{cproof}{lemma:CauchySchwarz_like} \Nobs that the H\"older inequality implies that
\begin{equation}\label{eqn:Schwarz_apply}
\int_{\R^d} \pabs{\ff(x) \fg(x)}\! \, dx \le \! \pb{\int_{\R^d} \pabs{\ff(x)}^2\! \, dx}^{\nicefrac{1}{2}}\! \pb{\int_{\R^d} \pabs{\fg(x)}^2 \! \, dx}^{\nicefrac{1}{2}} \! = 1.
\end{equation}

\noindent
Next \nobs that
\begin{equation}
\begin{split}
\int_{\R^d} \abs{\alpha \ff(x) - \fg(x)}^2 \, dx & = \alpha^2 + 1 - 2 \alpha \int_{\R^d} \ff(x) \fg(x) \, dx \\
& = \! \pb{\alpha - \int_{\R^d} \ff(x) \fg(x) \, dx}^2 \! + 1 - \!\pb{\int_{\R^d} \ff(x) \fg(x) \, dx}^2  \\
& \ge 1 - \!\pb{\int_{\R^d} \ff(x) \fg(x) \, dx}^2 \! \ge 1 - \!\pb{\int_{\R^d} \abs{\ff(x) \fg(x)} \, dx}^2 \!.
\end{split}
\end{equation}

\noindent
This and \cref{eqn:Schwarz_apply} ensure that
\begin{equation}
\int_{\R^d} \abs{\alpha \ff(x) - \fg(x)}^2 \, dx \ge 1 - \! \pb{\int_{\R^d} \abs{\ff(x) \fg(x)} \, dx}^2 \ge 1 - \int_{\R^d} \abs{\ff(x) \fg(x)} \, dx.
\end{equation}
\end{cproof}




\subsection{ANN approximations for a class of general high-dimensional functions}
\label{subsec:lower_bounds_anns_for_general_class}




\cfclear
\begin{theorem}\label{thm:main1}
Let $d \in \N \cap [4, \infty)$, $\beta, \expconst \in (0, \infty)$, $\Phi \in \ANNs$ satisfy $\inDimANN(\Phi) = d$ and $\outDimANN(\Phi) = 1$, let $\varphi\colon \R^d \to \R$, $g \colon \R^d \to \R$, and $\fg \colon \R^d\rightarrow \R$ be measurable, and assume for all $x \in \R^d$ that $\varphi(x)=(\nicefrac{\expconst}{\pi})^{\nicefrac{d}{2}}\exp\!\left(- \expconst \norm{x}_2^2\right)$, $\int_{\R^d}\abs{g(y)}^2 \varphi(y)\,dy \in (0, \infty)$, $\int_{\R^d}\abs{\pr{\functionANN\pr{\Phi}}(y)}\,dy>0$, and $\fg(x)=\big[\int_{\R^d}\abs{g(y)}^2\varphi(y)\,dy\big]^{\nicefrac{-1}{2}}g(x)$ \cfload. Then
\begin{equation}
\begin{split}
\lengthANN(\Phi) \! & \pb{\frac{\paramANN(\Phi) \max\!\pc{1, \norm{\vectorNN(\Phi)}_{\infty}}}{2 \lengthANN(\Phi)}}^{\lengthANN(\Phi)}\! \ge \! \pb{\frac{e^{\beta}}{1+\beta}}^{\nicefrac{d}{6}} \! \pb{\frac{2 \sqrt{\expconst}\pb{\int_{\R^d}\abs{\pr{\functionANN\pr{\Phi}}(x)}^2\varphi(x)\,dx}^{\nicefrac{1}{2}}}{d^{\nicefrac{3}{2}}(6+4 \beta + \expconst)^{\nicefrac{1}{2}}}} \\
\cdot & \!\left[1-\!\pb{\int_{\pc{y \in \R^d \colon \norm{y}_2 \le \frac{\sqrt{d(1+\beta)}}{\sqrt{2 \expconst}}}} \abs{\fg(x)}^2 \varphi(x) \, dx}^{\nicefrac{1}{2}}\!-\int_{\R^d}\abs{\pr{\functionANN\pr{\Phi}}(x)-\fg(x)}^2\varphi(x)\,dx \right]\!\ifnocf.
\end{split}
\end{equation}
\cfout[.]
\end{theorem}
\begin{cproof}{thm:main1}
Throughout this proof let $\mathbf{f} \colon \R^d \to \R$ and $\mathbf{g} \colon \R^d \to \R$ satisfy for all $x \in \R^d$ that $\mathbf{g} (x) = \fg(x) \! \pb{\varphi (x)}^{\nicefrac{1}{2}}$ and $\mathbf{f}(x) = \! \pb{\int_{\R^d}\abs{(\functionANN(\Phi))(y)}^2\varphi(y)\,dy}^{\nicefrac{-1}{2}} \!(\functionANN(\Phi))(x) \! \pb{\varphi(x)}^{\nicefrac{1}{2}}\!$ and let $\scrR \in \R$ satisfy $\sqrt{2 \expconst} \scrR = \sqrt{d(1+\beta)}$ \cfload. \Nobs that $\smallint_{\R^d} \abs{\mathbf{f}(x)}^2 \, dx = \int_{\R^d} \abs{\mathbf{g}(x)}^2 \, dx = 1$. \cref{lemma:CauchySchwarz_like} (applied with $d \with d$, $\alpha \with \! \pb{\int_{\R^d}\abs{(\functionANN(\Phi))(y)}^2\varphi(y)\,dy}^{\nicefrac{1}{2}}$, $\ff \with \mathbf{f}$, $\fg \with \mathbf{g}$ in the notation of \cref{lemma:CauchySchwarz_like}) hence ensures that
\begin{equation}
\begin{split}
\int_{\R^d}\abs{\pr{\functionANN\pr{\Phi}}(x)-\fg(x)}^2\varphi(x)\,dx & = \int_{\R^d} \! \pabs{\mathbf{f}(x) \! \pb{\int_{\R^d}\abs{(\functionANN(\Phi))(y)}^2\varphi(y)\,dy}^{\nicefrac{1}{2}}\! - \mathbf{g}(x)}^2 dx \\
& \ge 1-\int_{\R^d} \abs{\mathbf{f}(x)\mathbf{g}(x)} \, dx.
\end{split}
\end{equation}

\noindent
Combining this with \cref{lemma:upper_bounds_for_the_scalar_product} (applied with $d \with d$, $\beta \with \beta$, $\expconst \with \expconst$, $\Phi \with \Phi$, $\varphi \with \varphi$, $\ff \with \mathbf{f}$, $\fg \with \mathbf{g}$ in the notation of \cref{lemma:upper_bounds_for_the_scalar_product}) demonstrates that
\begin{equation}
\begin{split}
& \int_{\R^d}\abs{\pr{\functionANN\pr{\Phi}}(x)-\fg(x)}^2\varphi(x)\,dx \ge 1-\int_{\R^d} \abs{\mathbf{f}(x)\mathbf{g}(x)} \, dx \ge 1 - \pb{\int_{\{y \in \R^d \colon \norm{y}_2 \le \scrR\}} \abs{\mathbf{g}(x)}^2 \, dx}^{\nicefrac{1}{2}} \\
& - \lengthANN(\Phi) \! \pb{\frac{\paramANN(\Phi) \max\!\pc{1, \norm{\vectorNN(\Phi)}_{\infty}}}{2 \lengthANN(\Phi)}}^{\lengthANN(\Phi)} \! \pb{\frac{1+\beta}{e^{\beta}}}^{\nicefrac{d}{6}}\! \pb{\frac{d^{\nicefrac{3}{2}}(6+4\beta+\expconst)^{\nicefrac{1}{2}}}{2\sqrt{\expconst}\pb{\int_{\R^d} \abs{(\functionANN(\Phi))(y)}^2 \varphi(y) \, dy}^{\nicefrac{1}{2}}}} \\
& = 1 - \pb{\int_{\{y \in \R^d \colon \norm{y}_2 \le \scrR\}} \abs{\fg(x)}^2 \varphi(x) \, dx}^{\nicefrac{1}{2}}\\
& - \lengthANN(\Phi) \! \pb{\frac{\paramANN(\Phi) \max\!\pc{1, \norm{\vectorNN(\Phi)}_{\infty}}}{2 \lengthANN(\Phi)}}^{\lengthANN(\Phi)} \! \pb{\frac{1+\beta}{e^{\beta}}}^{\nicefrac{d}{6}}\! \pb{\frac{d^{\nicefrac{3}{2}}(6+4\beta+\expconst)^{\nicefrac{1}{2}}}{2\sqrt{\expconst}\pb{\int_{\R^d} \abs{(\functionANN(\Phi))(y)}^2 \varphi(y) \, dy}^{\nicefrac{1}{2}}}}\ifnocf.
\end{split}
\end{equation}

\noindent
\cfload[.]This implies that
\begin{multline}
\lengthANN(\Phi) \! \pb{\frac{\paramANN(\Phi) \max\!\pc{1, \norm{\vectorNN(\Phi)}_{\infty}}}{2 \lengthANN(\Phi)}}^{\lengthANN(\Phi)}\! \ge \! \pb{\frac{e^{\beta}}{1+\beta}}^{\nicefrac{d}{6}} \! \pb{\frac{2 \sqrt{\expconst}\pb{\int_{\R^d}\abs{\pr{\functionANN\pr{\Phi}}(x)}^2\varphi(x)\,dx}^{\nicefrac{1}{2}}}{d^{\nicefrac{3}{2}}(6+4 \beta + \expconst)^{\nicefrac{1}{2}}}} \\
\cdot \!\left[1-\!\pb{\int_{\{y \in \R^d \colon \norm{y}_2 \le \scrR\}} \abs{\fg(x)}^2 \varphi(x) \, dx}^{\nicefrac{1}{2}}\!-\int_{\R^d}\abs{\pr{\functionANN\pr{\Phi}}(x)-\fg(x)}^2\varphi(x)\,dx \right]\!.
\end{multline}
\end{cproof}


\subsection{ANN approximations for certain specific high-dimensional functions}
\label{subsec:lower_bounds_anns_for_specific_class}


\cfclear
\begin{corollary}\label{cor:main}
Let $d \in \N \cap [4, \infty)$, $\eps \in (0, \nicefrac{1}{4}]$, let $\varphi\colon \R^d \to \R$ and $g \colon \R^d \to \R$ satisfy for all $x = (x_1, x_2, \ldots, x_d) \in \R^d$ that $\varphi(x)=(2 \pi)^{-\nicefrac{d}{2}}\exp(- \frac{1}{2} \norm{x}_2^2)$ and $g(x) = \smallsum_{j=1}^d [\max\{\abs{x_j}-\sqrt{2d}, 0\}]^2$, let $\fg\colon \R^d\rightarrow \R$ satisfy for all $x \in \R^d$ that $\fg(x)=[\int_{\R^d}\abs{g(y)}^2\varphi(y)\,dy]^{\nicefrac{-1}{2}}g(x)$, and let $\Phi \in \ANNs$ satisfy $\inDimANN(\Phi) = d$, $\outDimANN(\Phi) = 1$, and $\int_{\R^d}\abs{\pr{\functionANN\pr{\Phi}}(x)-\fg(x)}^2\varphi(x)\,dx \le \eps$ \cfload. Then
\begin{equation}
\paramANN(\Phi) \max\!\pc{1, \norm{\vectorNN(\Phi)}_{\infty}} \! \ge \! \big[\tfrac{2}{7}\big] d^{\nicefrac{-3}{2}} \exp\!\big(\tfrac{d}{20 \lengthANN(\Phi)}\big)\!\ifnocf.
\end{equation}
\cfout[.]
\end{corollary}
\begin{cproof}{cor:main}
\Nobs that the triangle inequality ensures that
\begin{equation}\label{eqn:cor:main:norm_of_realization_lower_bound}
\begin{split}
& \pb{\int_{\R^d}\abs{(\functionANN(\Phi))(x)}^2\varphi(x)\,dx}^{\nicefrac{1}{2}}\! \\
& \ge \! \pb{\int_{\R^d}\abs{\fg(x)}^2\varphi(x)\,dx}^{\nicefrac{1}{2}}\! - \! \pb{\int_{\R^d}\abs{\pr{\functionANN\pr{\Phi}}(x)-\fg(x)}^2\varphi(x)\,dx}^{\nicefrac{1}{2}} \\
& = 1- \! \pb{\int_{\R^d}\abs{\pr{\functionANN\pr{\Phi}}(x)-\fg(x)}^2\varphi(x)\,dx}^{\nicefrac{1}{2}} \ge 1 - \eps^{\nicefrac{1}{2}} \\
& \ge 1-4^{\nicefrac{-1}{2}} = \frac{1}{2} > 0.
\end{split}
\end{equation}

\noindent
Hence, we obtain that $\int_{\R^d} \abs{(\functionANN(\Phi))(x)} \, dx >0$. Next \nobs that for all $x = (x_1, x_2, \ldots, x_d) \in \{y \in \R^d \colon \norm{y}_2 \le \sqrt{2 d}\}$, $j \in \{1, 2, \ldots, d\}$ it holds that $\abs{x_j} \le \norm{x}_2 \le \sqrt{2 d}$. This ensures that for all $x = (x_1, x_2, \ldots, x_d) \in \{y \in \R^d \colon \norm{y}_2 \le \sqrt{2 d}\}$ it holds that $\fg(x)=g(x) = 0$. Combining \cref{thm:main1} (applied with $d \with d$, $\beta \with 1$, $\expconst \with \nicefrac{1}{2}$, $\Phi \with \Phi$, $\varphi \with \varphi$, $g \with g$, $\fg \with \fg$ in the notation of \cref{thm:main1}), the fact that $\nicefrac{e}{2} \ge e^{\nicefrac{3}{10}}$, the fact that $\int_{\R^d} \abs{(\functionANN(\Phi))(x)} \, dx >0$, and \cref{eqn:cor:main:norm_of_realization_lower_bound} therefore  implies that
\begin{equation}
\begin{split}
& \lengthANN(\Phi) \! \pb{\frac{\paramANN(\Phi) \max\!\pc{1, \norm{\vectorNN(\Phi)}_{\infty}}}{2 \lengthANN(\Phi)}}^{\lengthANN(\Phi)}\! \ge \! \pb{\frac{e}{2}}^{\nicefrac{d}{6}} \! \pb{\frac{\sqrt{2} \pb{\int_{\R^d}\abs{\pr{\functionANN\pr{\Phi}}(x)}^2\varphi(x)\,dx}^{\nicefrac{1}{2}}}{d^{\nicefrac{3}{2}}(6+4+\nicefrac{1}{2})^{\nicefrac{1}{2}}}} \\
& \quad \cdot \!\left[1-\!\pb{\int_{\{y \in \R^d \colon \norm{y}_2 \le \sqrt{2 d}\}} \abs{\fg(x)}^2 \varphi(x) \, dx}^{\nicefrac{1}{2}}\!-\int_{\R^d}\abs{\pr{\functionANN\pr{\Phi}}(x)-\fg(x)}^2\varphi(x)\,dx \right] \\
& = \! \pb{\frac{e}{2}}^{\nicefrac{d}{6}} \! \pb{\frac{\sqrt{2} \pb{\int_{\R^d}\abs{\pr{\functionANN\pr{\Phi}}(x)}^2\varphi(x)\,dx}^{\nicefrac{1}{2}}}{d^{\nicefrac{3}{2}}(6+4+\nicefrac{1}{2})^{\nicefrac{1}{2}}}} \!\left[1-\int_{\R^d}\abs{\pr{\functionANN\pr{\Phi}}(x)-\fg(x)}^2\varphi(x)\,dx \right] \\
& \ge \, [(21)^{\nicefrac{-1}{2}}] \! \pb{\frac{e}{2}}^{\nicefrac{d}{6}}\! d^{\nicefrac{-3}{2}} (1-\eps) \ge \frac{e^{\nicefrac{d}{20}}}{7 d^{\nicefrac{3}{2}}}\ifnocf.
\end{split}
\end{equation}

\noindent
\cfload[.]Hence, we obtain that
\begin{equation}
\paramANN(\Phi) \max\!\pc{1, \norm{\vectorNN(\Phi)}_{\infty}} \ge 2\lengthANN(\Phi) \! \pb{\frac{e^{\nicefrac{d}{20}}}{7 d^{\nicefrac{3}{2}}\lengthANN(\Phi)}}^{\nicefrac{1}{\lengthANN(\Phi)}} \! \ge \!\pb{\frac{2}{7}} \! d^{\nicefrac{-3}{2}} \exp\!\pa{\frac{d}{20 \lengthANN(\Phi)}}\!.
\end{equation}
\end{cproof}

\cfclear
\begin{corollary}\label{cor:main_2}
Let $\varphi_d \colon \R^d \to \R$, $d \in \N$, and $g_d \colon \R^d \to \R$, $d \in \N$, satisfy for all $d \in \N$, $x = (x_1, x_2, \ldots, x_d) \in \R^d$ that $\varphi_d(x)=(2 \pi)^{\nicefrac{-d}{2}}\exp(- \frac{1}{2} (\smallsum_{j=1}^d \abs{x_j}^2))$ and $g_d(x) = \smallsum_{j=1}^d [\max\{\abs{x_j}-\sqrt{2 d}, 0\}]^2$, let $\fg_d\colon \R^d\rightarrow \R$, $d \in \N$, satisfy for all $d \in \N$, $x \in \R^d$ that $\fg_d(x)=[\int_{\R^d}\abs{g_d(y)}^2 \allowbreak \varphi_d(y)\,dy]^{-1/2}g_d(x)$, and let $\delta \in (0, 1]$, $\constantfrakC \in \! [100 (\delta \ln (1.03))^{-2}, \infty)$ satisfy $2 {\constantfrakC}^{\nicefrac{5}{\delta}} \le \! \allowbreak (1.03)^{\sqrt{\constantfrakC}}$ \cfload. Then it holds for all $\constantfrakc \in [\constantfrakC, \infty)$, $d \in \N$, $\eps \in (0, \nicefrac{1}{2}]$, $\Phi \in \ANNs$ with $\inDimANN(\Phi) = d$, $\outDimANN(\Phi) = 1$, $\hiddenLength(\Phi) \le \constantfrakc d^{1-\delta}$, $\norm{\vectorNN(\Phi)}_{\infty} \le \constantfrakc d^{\constantfrakc}$, and $\pb{\int_{\R^d}\abs{\pr{\functionANN\pr{\Phi}}(x)-\fg_d(x)}^2\varphi_d(x)\,dx}^{\nicefrac{1}{2}} \le \eps$ that $\paramANN(\Phi) \ge (1+{\constantfrakc}^{-3})^{(d^{\delta})}$ \cfout.
\end{corollary}
\begin{cproof}{cor:main_2}
\Nobs that the assumption that $\constantfrakC \in \! [100 (\delta \ln (1.03))^{-2}, \allowbreak \infty)$ and the chain rule ensure that for all $x \in [\constantfrakC, \infty)$ it holds that
\begin{equation}
\begin{split}
\big[2^{-1} (1.03)^{\sqrt{x}} x^{\nicefrac{-5}{\delta}}\big]' \! & = (1.03)^{\sqrt{x}} \ln (1.03) \! \pb{\frac{1}{4\sqrt{x}}} \! x^{\nicefrac{-5}{\delta}} - (1.03)^{\sqrt{x}} \! \pb{\frac{5}{2 \delta x}} \! x^{\nicefrac{-5}{\delta}} \\
& = (1.03)^{\sqrt{x}} \! \pb{\frac{x^{\nicefrac{-5}{\delta}}}{4  x}} \! \ln (1.03) \big[\sqrt{x} - 10 (\delta \ln (1.03))^{-1}\big] \\
& \ge (1.03)^{\sqrt{x}} \! \pb{\frac{x^{\nicefrac{-5}{\delta}}}{4  x}} \! \ln (1.03) \big[\sqrt{\constantfrakC} - 10 (\delta \ln (1.03))^{-1}\big] \ge 0.
\end{split}
\end{equation}

\noindent
This implies that the function $[\constantfrakC, \infty) \ni x \mapsto 2^{-1} (1.03)^{\sqrt{x}} x^{\nicefrac{-5}{\delta}} \in \R$ is non-decreasing. The assumption that $\constantfrakC \in \! [100 (\delta \ln (1.03))^{-2}, \allowbreak \infty)$ and the assumption that $2 {\constantfrakC}^{\nicefrac{5}{\delta}} \le \! \allowbreak (1.03)^{\sqrt{\constantfrakC}}$ therefore assure that for all $\constantfrakc \in [\constantfrakC, \infty)$ it holds that $\constantfrakc \ge 100 (\delta \ln (1.03))^{-2}$ and
\begin{equation}\label{eqn:cor:main_2_on_constantc}
2^{-1} (1.03)^{\sqrt{\constantfrakc}} {\constantfrakc}^{\nicefrac{-5}{\delta}} \ge 2^{-1} (1.03)^{\sqrt{\constantfrakC}} \, {\constantfrakC}^{\nicefrac{-5}{\delta}} \ge 1.
\end{equation}

\noindent
The fact that for all $x \in (0, \infty)$ it holds that $(1+x^{-1})^x \le e$ hence ensures that for all $\constantfrakc \in [\constantfrakC, \infty)$, $d \in \N$, $\Phi \in \ANNs$ with $d \le {\constantfrakc}^{5/(2 \delta)}$ it holds that
\begin{equation}\label{eqn:cor:main_2_small_d}
(1+{\constantfrakc}^{-3})^{(d^{\delta})} \le (1+{\constantfrakc}^{-3})^{({\constantfrakc}^{\nicefrac{5}{2}})} = \! \big[(1+{\constantfrakc}^{-3})^{({\constantfrakc}^3)}\big]^{\frac{1}{\sqrt{\constantfrakc}}} \le e^{\frac{1}{\sqrt{\constantfrakc}}} \le 2 \le \paramANN(\Phi)\ifnocf.
\end{equation}

\noindent
\cfload[.]Moreover, \nobs that the chain rule and \cref{eqn:cor:main_2_on_constantc} show that for all $\constantfrakc \in [\constantfrakC, \infty)$, $x \in [{\constantfrakc}^{5/(2 \delta)}, \infty)$ it holds that
\begin{equation}
\begin{split}
\big[(1.03)^{\nicefrac{(x^{\delta})}{\constantfrakc}} x^{-2 \constantfrakc} \big]' & = (1.03)^{\nicefrac{(x^{\delta})}{\constantfrakc}} \ln (1.03) \! \pb{\frac{\delta}{\constantfrakc}} \! x^{-2 \constantfrakc - 1 + \delta} - 2 \constantfrakc (1.03)^{\nicefrac{(x^{\delta})}{\constantfrakc}} x^{-2 \constantfrakc - 1} \\
& = (1.03)^{\nicefrac{(x^{\delta})}{\constantfrakc}} \, x^{-2 \constantfrakc - 1} \! \pb{\frac{\delta}{\constantfrakc}} \! \ln (1.03) [x^{\delta} - 2 {\constantfrakc}^2 (\delta \ln (1.03))^{-1}] \\
& \ge (1.03)^{\nicefrac{(x^{\delta})}{\constantfrakc}} \, x^{-2 \constantfrakc - 1} \! \pb{\frac{\delta}{\constantfrakc}} \! \ln (1.03) [{\constantfrakc}^{\nicefrac{5}{2}} - 2 {\constantfrakc}^2 (\delta \ln (1.03))^{-1}] \\
& = (1.03)^{\nicefrac{(x^{\delta})}{\constantfrakc}} \, x^{-2 \constantfrakc - 1} \delta \constantfrakc \ln (1.03) [\sqrt{\constantfrakc} - 2 (\delta \ln (1.03))^{-1}] \\
& \ge (1.03)^{\nicefrac{(x^{\delta})}{\constantfrakc}} \, x^{-2 \constantfrakc - 1} \, 8 \constantfrakc  > 0.
\end{split}
\end{equation}

\noindent
This implies for all $\constantfrakc \in [\constantfrakC, \infty)$ that the function $[{\constantfrakc}^{5/(2 \delta)}, \infty) \ni x \mapsto (1.03)^{\nicefrac{(x^{\delta})}{\constantfrakc}} x^{-2 \constantfrakc} \in \R$ is strictly increasing. The fact that $e^{\nicefrac{1}{30}} \ge 1.03$, \cref{eqn:cor:main_2_on_constantc}, and the fact that for all $\constantfrakc \in [\constantfrakC, \infty)$ it holds that $2^{\constantfrakc + 1} \ge 7 \constantfrakc$ therefore demonstrate that for all $\constantfrakc \in [\constantfrakC, \infty)$, $d \in \N$ with $d \ge {\constantfrakc}^{5/(2 \delta)}$ it holds that
\begin{equation}
e^{(d^{\delta})/(30 \constantfrakc)} d^{-2 \constantfrakc} \ge (1.03)^{\nicefrac{(d^{\delta})}{\constantfrakc}} d^{-2 \constantfrakc} \ge (1.03)^{({\constantfrakc}^{\nicefrac{3}{2}})} {\constantfrakc}^{\nicefrac{-5 \constantfrakc}{\delta}} = \! \big[ (1.03)^{\sqrt{\constantfrakc}} {\constantfrakc}^{\nicefrac{-5}{\delta}}\big]^{\constantfrakc} \! \ge 2^{\constantfrakc} \ge \! \pb{\frac{7}{2}} \! \constantfrakc.
\end{equation}

\noindent
The fact that for all $\constantfrakc \in [\constantfrakC, \infty)$ it holds that $(25 \constantfrakc)^{-1} \ge (30 \constantfrakc)^{-1} + {\constantfrakc}^{-3}$, the fact that for all $x \in \R$ it holds that $e^x \ge 1+x$, and \cref{cor:main} hence ensure that for all $\constantfrakc \in [\constantfrakC, \infty)$, $d \in \N$, $\eps \in (0, \nicefrac{1}{2}]$, $\Phi \in \ANNs$ with $d \ge {\constantfrakc^{5/(2 \delta)}}$, $\inDimANN(\Phi) = d$, $\outDimANN(\Phi) = 1$, $\hiddenLength(\Phi) \le \constantfrakc d^{1-\delta}$, $\norm{\vectorNN(\Phi)}_{\infty} \le \constantfrakc d^{\constantfrakc}$, and $\pb{\int_{\R^d}\abs{\pr{\functionANN\pr{\Phi}}(x)-\fg_d(x)}^2\varphi_d(x)\,dx}^{\nicefrac{1}{2}} \le \eps$ it holds that
\begin{equation}
\begin{split}
\paramANN(\Phi) & \ge (\max\{1, \norm{\vectorNN(\Phi)}_{\infty}\})^{-1} \! \pb{\tfrac{2}{7}} \! d^{\nicefrac{-3}{2}} \exp\!\big(\tfrac{d}{20 \lengthANN(\Phi)}\big) \ge \! \pb{\tfrac{2}{7}} \! \exp\!\big(\tfrac{d^{\delta}}{25 \constantfrakc}\big) d^{- 2 \constantfrakc} {\constantfrakc}^{-1} \\
& \ge \! \pb{\tfrac{2}{7}} \! \exp\!\big(\tfrac{d^{\delta}}{30 \constantfrakc}\big) d^{- 2 \constantfrakc} {\constantfrakc}^{-1} \exp\!\big(\tfrac{d^{\delta}}{{\constantfrakc}^3}\big) \ge \exp\!\big(\tfrac{d^{\delta}}{{\constantfrakc}^3}\big) \ge (1+{\constantfrakc}^{-3})^{(d^{\delta})}\ifnocf.
\end{split}
\end{equation}
\cfload[.]Combining this with \cref{eqn:cor:main_2_small_d} assures that for all $\constantfrakc \in [\constantfrakC, \infty)$, $d \in \N$, $\eps \in (0, \nicefrac{1}{2}]$, $\Phi \in \ANNs$ with $\inDimANN(\Phi) = d$, $\outDimANN(\Phi) = 1$, $\hiddenLength(\Phi) \le \constantfrakc d^{1-\delta}$, $\norm{\vectorNN(\Phi)}_{\infty} \le \constantfrakc d^{\constantfrakc}$, and $\pb{\int_{\R^d}\abs{\pr{\functionANN\pr{\Phi}}(x)-\fg_d(x)}^2\varphi_d(x)\,dx}^{\nicefrac{1}{2}} \le \eps$ it holds that $\paramANN(\Phi) \ge (1+{\constantfrakc}^{-3})^{(d^{\delta})}$.
\end{cproof}



\section{Upper bounds for the number of ANN parameters in the approximation of high-dimensional functions}\label{sec:upper_bounds_for_number_of_parameters_in_ANN_approximations}

In this section we establish in \cref{cor:main_3} in Subsection~\ref{subsec:upper_bounds_anns_for_specific_class} below appropriate upper bounds for the number of parameters of suitable ANNs that approximate certain high-dimensional target functions. \cref{cor:main_3} is a consequence of the ANN approximation result in \cref{thm:main2} in Subsection~\ref{subsec:upper_bounds_anns_for_specific_class} below. Our proof of \cref{thm:main2} employs (i) the elementary ANN representation result for multiplications with powers of real numbers which we establish in \cref{lemma:on_Realization_of_DNN_2} in Subsection~\ref{subsec:ann_representation_powers_of_reals} below, (ii) the lower and upper bounds for appropriate Gaussian integrals which we present in \cref{lemma:bounds_for_the_norm_of_specific_function} in Subsection~\ref{subsec:bounds_for_the_norm_of_specific_function_lower} below, and (iii) the ANN approximation result for appropriate shifted squared rectifier functions in \cref{cor:approximation_of_rectifier_square} in Subsection~\ref{subsec:ann_approximation_for_shifted_squared_rectifier} below.

Our proof of \cref{lemma:bounds_for_the_norm_of_specific_function} employs the well-known Gaussian tail estimates in \cref{lemma:lower_bound_gaussian_tail_from_Klenke,lemma:lower_bound_gaussian_tail} in Subsection~\ref{subsec:bounds_for_the_norm_of_specific_function_lower} below. \cref{lemma:lower_bound_gaussian_tail_from_Klenke} is, e.g., proved as Lemma~22.2 in Klenke~\cite{k08b} and only for completeness we include in Subsection~\ref{subsec:bounds_for_the_norm_of_specific_function_lower} also the detailed proofs for \cref{lemma:lower_bound_gaussian_tail_from_Klenke,lemma:lower_bound_gaussian_tail}. 
Our proof of \cref{cor:approximation_of_rectifier_square} uses the elementary ANN representation result for compositions with shifted absolute value functions which we present in \cref{lemma:realization_absolute-c} in Subsection~\ref{subsec:ann_approximation_for_shifted_squared_rectifier} below as well as the ANN approximation result for the squared rectifier function in \cref{lemma:approximation_of_rectifier_square} in Subsection~\ref{subsec:ann_approximation_squared_rectifier} below.

Our proof of \cref{lemma:approximation_of_rectifier_square} uses the well-known ANN approximation result for the square function in \cref{cor:approximation_of_square_on_segment} in Subsection~\ref{subsec:ann_approximation_square_function} below. The proof of \cref{cor:approximation_of_square_on_segment}, in turn, employs the well-known ANN representation result in \cref{lemma:recurrence_on_functions_g_f_r} in Subsection~\ref{subsec:ann_approximation_square_function}. \cref{lemma:recurrence_on_functions_g_f_r,cor:approximation_of_square_on_segment} and their proofs are stongly based on Yarotsky~\cite[Proposition~2]{Yarotsky2017Errorbounds}. In the current form \cref{lemma:recurrence_on_functions_g_f_r,cor:approximation_of_square_on_segment} and their proofs are slight extensions of, e.g., the statement and the proof of Proposition~3.3 in Grohs et al.~\cite{Zimmermann2019spacetime} (cf., e.g., also Elbr\"{a}chter et al.~\cite[Lemma 6.1]{ElbraechterSchwab2018}). Only for completeness we include in Subsection~\ref{subsec:ann_approximation_square_function} also the detailed proofs for \cref{lemma:recurrence_on_functions_g_f_r,cor:approximation_of_square_on_segment}.

\subsection{ANN approximations for the square function}
\label{subsec:ann_approximation_square_function}


\cfclear
\begin{lemma}\label{lemma:recurrence_on_functions_g_f_r}
Let $(A_k)_{k\in\N} \subseteq \R^{4 \times 4}$, $\mathbb{B} \in \R^{4 \times 1}$, $(c_k)_{k \in \N} \subseteq \R$ satisfy for all $k\in \N$ that
\begin{equation}\label{eqn:lemma:recurrence_on_functions_g_f_r:matrices}
A_k = \! 
\begin{pmatrix}
2 & -4 & 2 & 0\\
2 & -4 & 2 & 0\\
2 & -4 & 2 & 0\\
-c_k & 2c_k & -c_k & 1
\end{pmatrix}\!,
\qquad
\mathbb{B} = \! 
\begin{pmatrix}
0\\
-\frac{1}{2}\\
-1\\
0
\end{pmatrix}\!, \qquad \text{and} \qquad c_k = 2^{1-2k},
\end{equation}

\noindent
let $g_n \colon \R \to [0,1]$, $n \in \N$, satisfy for all $n \in \N$, $x \in \R$ that
\begin{equation}\label{eqn:lemma:recurrence_on_functions_g_f_r:g_1}
g_1(x)=
\begin{cases}
2x & \colon x \in [0,\frac{1}{2}) \\[1ex]
2-2x & \colon x \in [\frac{1}{2},1] \\[1ex]
0 & \colon x \in \R \backslash [0,1] \\
\end{cases}
\end{equation}

\noindent
and $g_{n+1}(x)=g_1(g_{n}(x))$, let $f_n \colon [0,1] \to [0,1]$, $n \in \N_0$, satisfy for all $n \in \N_0$, $k \in \{0, 1, \dots, 2^n - 1\}$, $x \in \! \big[\tfrac{k}{2^n}, \tfrac{k+1}{2^n}\big)$ that $f_n(1) = 1$ and
\begin{equation}\label{eqn:lemma:recurrence_on_functions_g_f_r:f_n}
f_n(x) = \! \big[\tfrac{2k + 1}{2^n}\big] x - \tfrac{(k^2 + k)}{2^{2n}},
\end{equation}

\noindent
and let $r_k=(r_{k,1}, r_{k,2}, r_{k,3}, r_{k,4}) \colon \R \to \R^4$, $k \in \N$, satisfy for all $k \in \N$, $x \in \R$ that $r_1(x) = \Rect(x, x-\tfrac{1}{2}, x-1, x)$ and $r_{k+1}(x) = \Rect(A_{k} r_{k}(x) + \mathbb{B})$ \cfload. Then
\begin{enumerate}[label=(\roman *)]
\item
\label{item1:lemma:recurrence_on_functions_g_f_r} it holds for all $k \in \N$, $x \in \R$ that 
\begin{equation}\label{eqn:lemma:recurrence_on_functions_g_f_r:rec_r&g}
2r_{k,1}(x)-4r_{k,2}(x)+2r_{k,3}(x)=g_k(x)
\end{equation}

and
\item
\label{item2:lemma:recurrence_on_functions_g_f_r} it holds for all $k \in \N$, $x\in \R$ that
\begin{equation}\label{eqn:lemma:recurrence_on_functions_g_f_r:rec_r&f}
r_{k,4}(x)
=
\begin{cases}
f_{k-1}(x) & \colon x \in [0,1] \\
\max\{x, 0\} & \colon x \in \R \backslash [0,1]. \\
\end{cases}
\end{equation}
\end{enumerate}
\cfout[.]
\end{lemma}
\begin{cproof}{lemma:recurrence_on_functions_g_f_r}
We prove \cref{eqn:lemma:recurrence_on_functions_g_f_r:rec_r&g,eqn:lemma:recurrence_on_functions_g_f_r:rec_r&f} by induction on $k \in \N$. \Nobs that \cref{eqn:lemma:recurrence_on_functions_g_f_r:g_1} and the assumption that for all $x \in \R$ it holds that $r_1(x) = \Rect\big(x, x-\tfrac{1}{2}, x-1, x\big)$ show that for all $x \in \R$ it holds that
\begin{multline}\label{eqn:cor:approximation_of_square_on_segment:function_r_initial_eq}
2r_{1,1}(x) - 4 r_{1,2}(x) + 2 r_{1,3}(x) = 2 \Rect(x) - 4 \Rect(x-\tfrac{1}{2}) + 2 \Rect(x-1) \\
= 2 \max\{x, 0\} - 4 \max\{x-\tfrac{1}{2}, 0\} + 2 \max\{x-1, 0\}
= g_1(x).
\end{multline}

\noindent
Furthermore, \nobs that the assumption that for all $x \in \R$ it holds that $r_1(x) = \Rect(x, x-\tfrac{1}{2}, x-1, x)$ and the fact that for all $x \in [0,1]$ it holds that $f_0(x) = x = \max\{x, 0\}$ imply that for all $x \in \R$ it holds that 
\begin{equation}\label{eqn:cor:approximation_of_square_on_segment:function_r_initial_eqb}
r_{1,4}(x) = \max\{x, 0\} =
\begin{cases}
f_{0}(x) & \colon x \in [0,1] \\
\max\{x, 0\} & \colon x \in \R \backslash [0,1]. \\
\end{cases}
\end{equation}

\noindent 
Combining this with \cref{eqn:cor:approximation_of_square_on_segment:function_r_initial_eq} proves \cref{eqn:lemma:recurrence_on_functions_g_f_r:rec_r&g,eqn:lemma:recurrence_on_functions_g_f_r:rec_r&f} in the base case $k=1$. For the induction step let $k \in \N$ satisfy for all $x \in \R$ that
\begin{equation}\label{eqn:cor:approximation_of_square_on_segment:function_r_induction_initial:InductionHypothesis}
2 r_{k, 1}(x) - 4 r_{k, 2}(x) + 2 r_{k, 3}(x) = g_k(x)
\end{equation}

\noindent
and
\begin{equation}\label{eqn:cor:approximation_of_square_on_segment:function_r_induction_b:InductionHypothesis}
r_{k, 4}(x) =
\begin{cases}
f_{k-1}(x) & \colon x \in [0,1] \\
\max\{x, 0\} & \colon x \in \R \backslash [0,1]. \\
\end{cases}
\end{equation}

\noindent
\Nobs that \cref{eqn:lemma:recurrence_on_functions_g_f_r:matrices}, \cref{eqn:cor:approximation_of_square_on_segment:function_r_initial_eq}, \cref{eqn:cor:approximation_of_square_on_segment:function_r_induction_initial:InductionHypothesis}, and the assumption that for all $n \in \N$, $x \in \R$ it holds that $r_{n+1}(x) = \Rect(A_{n} r_{n}(x) + \mathbb{B})$ ensure that for all $x \in \R$ it holds that
\begin{equation}\label{eqn:cor:approximation_of_square_on_segment:function_r_induction_initial:InductionCalculation}
\begin{split}
g_{k+1}(x) & = g_1(g_{k}(x)) = g_1(2 r_{k,1}(x) - 4 r_{k,2}(x) + 2 r_{k,3}(x)) \\
& = 2 \Rect\big(2 r_{k,1}(x) - 4 r_{k,2}(x) + 2 r_{k,3}(x)\big) \\ 
& \quad - 4 \Rect\big(2 r_{k,1}(x) - 4 r_{k,2}(x) + 2 r_{k,3}(x) -\tfrac{1}{2}\big) \\
& \quad + 2 \Rect\big(2 r_{k,1}(x) - 4 r_{k,2}(x) + 2 r_{k,3}(x) - 1\big) \\
& = 2 r_{k+1,1}(x) - 4 r_{k+1,2}(x) + 2 r_{k+1,3}(x).
\end{split}
\end{equation}

\noindent
In addition, \nobs that \cref{eqn:lemma:recurrence_on_functions_g_f_r:matrices}, \cref{eqn:cor:approximation_of_square_on_segment:function_r_induction_initial:InductionHypothesis}, and the assumption that for all $n \in \N$, $x \in \R$ it holds that $r_{n+1}(x) = \Rect(A_{n} r_{n}(x) + \mathbb{B})$ demonstrate that for all $x \in \R$ it holds that 
\begin{equation}\label{eqn:cor:approximation_of_square_on_segment:Induction_FourthComponentOne}
\begin{split}
r_{k+1,4}(x) & = \Rect(- c_k r_{k,1}(x) + 2 c_k r_{k,2}(x) - c_k r_{k,3}(x) + r_{k,4}(x)) \\
& = \Rect(- [2^{1-2k}] r_{k,1}(x) + [2^{2-2k}] r_{k,2}(x) - [2^{1-2k}] r_{k,3}(x) + r_{k,4}(x)) \\
& = \Rect(-[2^{-2k}][2 r_{k,1}(x) - 4 r_{k,2}(x) + 2 r_{k,3}(x)] + r_{k,4}(x)) \\
& = \Rect(-[2^{-2k}] g_k(x) + r_{k,4}(x)).
\end{split}
\end{equation}

\noindent
Combining this with \cref{eqn:cor:approximation_of_square_on_segment:function_r_induction_b:InductionHypothesis}, \cite[Lemma 3.2]{Zimmermann2019spacetime}, and the fact that for all $x \in [0,1]$ it holds that $f_k(x) \ge 0$ shows that for all $x \in [0,1]$ it holds that
\begin{equation}\label{eqn:cor:approximation_of_square_on_segment:Induction_CalculationFourthComponentTwo}
\begin{split}
r_{k+1,4}(x) & = \Rect(-[2^{-2k}] g_k(x) + r_{k, 4}(x)) = \Rect(-[2^{-2k} g_k(x)] + f_{k-1}(x)) \\
& = \Rect\big(-[2^{-2k} g_k(x)] + x - \! \big[\smallsum\nolimits_{j=1}^{k-1} [2^{-2j} g_j(x)]\big]\big) \\
& = \Rect\big(x - \!\big[\smallsum\nolimits_{j=1}^{k} [2^{-2j} g_j(x)]\big]\big) \! = \Rect(f_k(x)) = f_k(x).
\end{split}
\end{equation}

\noindent
Next \nobs that \cref{eqn:cor:approximation_of_square_on_segment:function_r_induction_b:InductionHypothesis}, \cref{eqn:cor:approximation_of_square_on_segment:Induction_FourthComponentOne}, and the fact that for all $x \in \R \backslash [0, 1]$ it holds that $g_k(x) = 0$ prove that for all $x \in \R \backslash [0,1]$ it holds that 
\begin{equation}\label{eqn:cor:approximation_of_square_on_segment:Induction_CalculationFourthComponentThree}
r_{k+1,4}(x) = \Rect(-[2^{-2k}] g_k(x) + r_{k,4}(x)) = \Rect(r_{k, 4}(x)) = \Rect(\max\{x, 0\}) = \max\{x, 0\}.
\end{equation}

\noindent
Combining \cref{eqn:cor:approximation_of_square_on_segment:function_r_induction_initial:InductionCalculation,eqn:cor:approximation_of_square_on_segment:Induction_CalculationFourthComponentTwo} hence proves \cref{eqn:lemma:recurrence_on_functions_g_f_r:rec_r&g,eqn:lemma:recurrence_on_functions_g_f_r:rec_r&f} in the case $k+1$. Induction thus establishes \cref{item1:lemma:recurrence_on_functions_g_f_r,item2:lemma:recurrence_on_functions_g_f_r}.
\end{cproof}

\cfclear
\begin{lemma}\label{cor:approximation_of_square_on_segment}
Let $M \in \N$, $(A_k)_{k\in\N} \subseteq \R^{4 \times 4}$, $\mathbb{A}, \mathbb{B} \in \R^{4 \times 1}$, $(C_k)_{k\in\N} \subseteq \R^{1 \times 4}$, $(c_k)_{k \in \N} \subseteq \R$ satisfy for all $k\in \N$ that
\begin{equation}\label{eqn:cor:approximation_of_square_on_segment:matrices}
A_k = \! \begin{pmatrix}
2 & -4 & 2 & 0\\
2 & -4 & 2 & 0\\
2 & -4 & 2 & 0\\
-c_k & 2c_k & -c_k & 1
\end{pmatrix}\!,
\quad
\mathbb{A} = \! \begin{pmatrix}
1\\
1\\
1\\
1
\end{pmatrix}\!,
\quad
\mathbb{B} = \! \begin{pmatrix}
0\\
-\frac{1}{2}\\
-1\\
0
\end{pmatrix}\!, \quad
C_k = \!\begin{pmatrix}
-c_k & 2c_k & -c_k & 1
\end{pmatrix}\!,
\end{equation}

\noindent
and $c_k = 2^{1-2k}$ and let $\Phi \in \ANNs$ satisfy \begin{equation}\label{eqn:cor:approximation_of_square_on_segment:Phi_intro}
\Phi = \begin{cases}
\pa{(\mathbb{A}, \mathbb{B}), (C_1, 0)} & \colon M = 1\\
\pa{(\mathbb{A}, \mathbb{B}), (A_1, \mathbb{B}), (A_2, \mathbb{B}), \dots, (A_{M-1}, \mathbb{B}),
(C_M, 0)} & \colon M > 1\ifnocf.
\end{cases}
\end{equation} \cfload. Then
\begin{enumerate}[label=(\roman *)]
\item
\label{item1:cor:approximation_of_square_on_segment}
it holds that $\functionANN(\Phi) \in C(\R,\R)$,

\item
\label{item2:cor:approximation_of_square_on_segment}
it holds for all $x \in [0,1]$ that $\abs{x^2 - (\functionANN(\Phi))(x)} \le 4^{-M-1}$,

\item
\label{item3:cor:approximation_of_square_on_segment}
it holds for all $x \in \R \backslash [0,1]$ that $(\functionANN(\Phi))(x) = \Rect(x)$,

\item
\label{item4:cor:approximation_of_square_on_segment}
it holds that $\dims(\Phi) = (1,4,4,\dots,4,1) \in \N^{M+2}$,

\item
\label{item5:cor:approximation_of_square_on_segment}
it holds that $\infnorm{\vectorNN(\Phi)} \le 4$,

\item
\label{item6:cor:approximation_of_square_on_segment}
it holds that $\hiddenLength(\Phi) = M$, and

\item
\label{item7:cor:approximation_of_square_on_segment}
it holds that $\paramANN(\Phi) = 20M - 7$\ifnocf.
\end{enumerate}
\cfout[.]
\end{lemma}
\begin{cproof}{cor:approximation_of_square_on_segment}
Throughout this proof let $g_n \colon \R \to [0,1]$, $n \in \N$, satisfy for all $n \in \N$, $x \in \R$ that
\begin{equation}\label{eqn:cor:approximation_of_square_on_segment:function_g_1}
g_1(x)=
\begin{cases}
2x & \colon x \in [0,\frac{1}{2}) \\
2-2x & \colon x \in [\frac{1}{2},1] \\
0 & \colon x \in \R \backslash [0,1] \\
\end{cases}
\end{equation}

\noindent
and $g_{n+1}(x)=g_1(g_{n}(x))$, let $f_n \colon [0,1] \to [0,1]$, $n \in \N_0$, satisfy for all $n \in \N_0$, $k \in \{0, 1, \dots, 2^n - 1\}$, $x \in \! \big[\tfrac{k}{2^n}, \tfrac{k+1}{2^n}\big)$ that $f_n(1) = 1$ and
\begin{equation}\label{eqn:cor:approximation_of_square_on_segment:function_f_k}
f_n(x) = \! \big[\tfrac{2k + 1}{2^n}\big] x - \tfrac{(k^2 + k)}{2^{2n}},
\end{equation}

\noindent
and let $r_k=(r_{k,1}, r_{k,2}, r_{k,3}, r_{k,4}) \colon \R \to \R^4$, $k \in \N$, satisfy for all $k \in \N$, $x \in \R$ that 
\begin{equation}\label{eqn:cor:approximation_of_square_on_segment:function_r_1}
r_1(x) = \Rect\big(x, x-\tfrac{1}{2}, x-1, x\big)
\end{equation}  

\noindent
and 
\begin{equation}\label{eqn:cor:approximation_of_square_on_segment:function_r_k+1}
r_{k+1}(x) = \Rect\big(A_{k} r_{k}(x) + \mathbb{B}\big)\ifnocf.
\end{equation}

\noindent
\cfload[.]\Nobs that \cref{item1:lemma:recurrence_on_functions_g_f_r} in \cref{lemma:recurrence_on_functions_g_f_r} (applied with $(A_k)_{k \in \N} \with (A_k)_{k \in \N}$, $\mathbb{B} \with \mathbb{B}$, $(c_k)_{k \in \N} \with (c_k)_{k \in \N}$, $(g_n)_{n \in \N} \with (g_n)_{n \in \N}$, $(f_n)_{n \in \N} \with (f_n)_{n \in \N}$, $(r_k)_{k \in \N} \with (r_k)_{k \in \N}$ in the notation of \cref{lemma:recurrence_on_functions_g_f_r}), \cref{eqn:cor:approximation_of_square_on_segment:matrices,eqn:cor:approximation_of_square_on_segment:Phi_intro,eqn:cor:approximation_of_square_on_segment:function_r_1,eqn:cor:approximation_of_square_on_segment:function_r_k+1} assure that for all $x \in \R$ it holds that
\begin{equation}\label{NNsquare:CalculationANNOne}
\begin{split}
(\functionANN(\Phi))(x) & = - c_M r_{M, 1}(x) + 2 c_M r_{M, 2}(x) - c_M r_{M, 3}(x) + r_{M, 4}(x) \\
& = -[2^{1-2M}] r_{M, 1}(x) + [2^{2-2M}] r_{M, 2}(x) - [2^{1-2M}] r_{M, 3}(x) + r_{M, 4}(x) \\
& = -[2^{-2M}][2 r_{M, 1}(x) - 4 r_{M, 2}(x) + 2 r_{M, 3}(x)] + r_{M, 4}(x) \\
& = -[2^{-2M}] g_{M}(x) + r_{M, 4}(x)\ifnocf.
\end{split}
\end{equation} 

\noindent
\cfload[.]This establishes \cref{item1:cor:approximation_of_square_on_segment}. Moreover, \nobs that \cref{NNsquare:CalculationANNOne}, \cite[Lemma 3.2]{Zimmermann2019spacetime}, and \cref{item2:lemma:recurrence_on_functions_g_f_r} in \cref{lemma:recurrence_on_functions_g_f_r} (applied with $(A_k)_{k \in \N} \with (A_k)_{k \in \N}$, $\mathbb{B} \with \mathbb{B}$, $(c_k)_{k \in \N} \with (c_k)_{k \in \N}$, $(g_n)_{n \in \N} \with (g_n)_{n \in \N}$, $(f_n)_{n \in \N} \with (f_n)_{n \in \N}$, $(r_k)_{k \in \N} \with (r_k)_{k \in \N}$ in the notation of \cref{lemma:recurrence_on_functions_g_f_r}) show that for all $x \in [0,1]$ it holds that
\begin{equation}\label{eqn:cor:approximation_of_square_on_segment:CalculationANNTwo}
\begin{split}
(\functionANN(\Phi))(x) & = -[2^{-2M}] g_{M}(x) + r_{M, 4}(x) = -[2^{-2M} g_{M}(x)] + f_{M-1}(x) \\
& = - [2^{-2M} g_{M}(x)] + x - \! \big[\smallsum\nolimits_{j=1}^{M-1} [2^{-2j} g_j(x)]\big] \\
& = x - \! \big[\smallsum\nolimits_{j=1}^{M} [2^{-2j} g_j(x)]\big] \! = f_{M}(x).
\end{split}
\end{equation}

\noindent
This and \cite[Lemma~3.2]{Zimmermann2019spacetime} imply that for all $x \in [0,1]$ it holds that 
\begin{equation}\label{eqn:cor:approximation_of_square_on_segment:approx_error}
\abs{x^2 - (\functionANN(\Phi))(x)} = \abs{x^2 - f_M(x)} \leq 2^{-2M-2} = 4^{-M-1}.
\end{equation}

\noindent
This establishes \cref{item2:cor:approximation_of_square_on_segment}.
Furthermore, \nobs that \cref{NNsquare:CalculationANNOne}, the fact that for all $x \in \R \backslash [0, 1]$ it holds that $g_M(x) = 0$, and \cref{item2:lemma:recurrence_on_functions_g_f_r} in \cref{lemma:recurrence_on_functions_g_f_r} (applied with $(A_k)_{k \in \N} \with (A_k)_{k \in \N}$, $\mathbb{B} \with \mathbb{B}$, $(c_k)_{k \in \N} \with (c_k)_{k \in \N}$, $(g_n)_{n \in \N} \with (g_n)_{n \in \N}$, $(f_n)_{n \in \N} \with (f_n)_{n \in \N}$, $(r_k)_{k \in \N} \with (r_k)_{k \in \N}$ in the notation of \cref{lemma:recurrence_on_functions_g_f_r}) ensure that for all $x \in \R \backslash [0,1]$ it holds that
\begin{equation}\label{eqn:cor:approximation_of_square_on_segment:ANNoutside}
(\functionANN(\Phi))(x) = - [2^{-2M}] g_{M}(x) + r_{M,4}(x) = r_{M,4}(x) = \max\{x, 0\} = \Rect(x).
\end{equation}

\noindent
This establishes \cref{item3:cor:approximation_of_square_on_segment}. In addition, \nobs that \cref{eqn:cor:approximation_of_square_on_segment:matrices,eqn:cor:approximation_of_square_on_segment:Phi_intro} imply that $\dims(\Phi) = (1,4,4,\dots,4,1) \in \N^{M+2}$, $\norm{\vectorNN(\Phi)}_{\infty} \le 4$, $\hiddenLength(\Phi) = M$, and 
\begin{equation}\label{NNsLMphi_2}
\paramANN(\Phi) = 4(1+1) + \! \big[\smallsum\nolimits_{j=2}^{M}4(4+1)\big] \! + (4+1) = 8 + 20(M - 1) + 5 = 20M - 7.
\end{equation}

\noindent
This establishes \cref{item4:cor:approximation_of_square_on_segment,item5:cor:approximation_of_square_on_segment,item6:cor:approximation_of_square_on_segment,item7:cor:approximation_of_square_on_segment}.
\end{cproof}




\subsection{ANN approximations for the squared rectifier function}
\label{subsec:ann_approximation_squared_rectifier}


\cfclear
\begin{corollary}\label{lemma:approximation_of_rectifier_square}
Let $M \in \N$, $\constantR \in [1, \infty)$, $q \in (2, \infty)$, $(A_k)_{k \in \N} \subseteq \R^{4 \times 4}$, $\mathbb{A}, \mathbb{B} \in \R^{4 \times 1}$, $(C_k)_{k \in \N} \subseteq \R^{1 \times 4}$, $(c_k)_{k \in \N} \subseteq \R$ satisfy for all $k \in \N$ that
\begin{equation}\label{eqn:lemma:approximation_of_rectifier_square_weights_biases}
A_k = \! \begin{pmatrix}
2 & -4 & 2 & 0\\
2 & -4 & 2 & 0\\
2 & -4 & 2 & 0\\
-c_k & 2c_k & -c_k & 1
\end{pmatrix}\!,
\quad
\mathbb{A} = \! \begin{pmatrix}
1\\
1\\
1\\
1
\end{pmatrix}\!,
\quad
\mathbb{B} = \! \begin{pmatrix}
0\\
-\frac{1}{2}\\
-1\\
0
\end{pmatrix}\!,
\quad
C_k = \!\begin{pmatrix}
-c_k & 2c_k & -c_k & 1
\end{pmatrix}\!,
\end{equation}

\noindent
and $c_k = 2^{1-2k}$, and let $\Psi, \Phi \in \ANNs$ satisfy
\begin{equation}\label{eqn:lemma:approximation_of_rectifier_square_Psi}
\Psi = \begin{cases}
\pa{(\mathbb{A}, \mathbb{B}), (C_1, 0)} & \colon M = 1\\
\pa{(\mathbb{A}, \mathbb{B}), (A_1, \mathbb{B}), (A_2, \mathbb{B}), \dots, (A_{M-1}, \mathbb{B}),
(C_M, 0)} & \colon M > 1
\end{cases}
\end{equation}

\noindent
and $\Phi = \AffineANN_{{\constantR}^2, 0} \compANNbullet \Psi \compANNbullet \AffineANN_{{\constantR}^{-1}, 0}$ \cfload. Then
\begin{enumerate}[label=(\roman *)]
\item
\label{item1:lemma:approximation_of_rectifier_square} it holds that $\functionANN(\Phi)\in C(\R,\R)$,

\item
\label{item2:lemma:approximation_of_rectifier_square} it holds for all $x \in (-\infty, 0]$ that $\pabs{[\Rect(x)]^2 - (\functionANN(\Phi))(x)}=0$,

\item
\label{item3:lemma:approximation_of_rectifier_square} it holds for all $x \in [0, \constantR]$ that $\pabs{[\Rect(x)]^2 - (\functionANN(\Phi))(x)} \le 4^{-M-1} {\constantR}^2$,

\item
\label{item4:lemma:approximation_of_rectifier_square} it holds for all $x \in [\constantR, \infty)$ that $\pabs{[\Rect(x)]^2 - (\functionANN(\Phi))(x)} \le \abs{\Rect(x)}^q {\constantR}^{2-q}$,

\item
\label{item5:lemma:approximation_of_rectifier_square} it holds that $\dims(\Phi) = (1, 4, \ldots, 4, 1) \in \N^{M+2}$,

\item
\label{item6:lemma:approximation_of_rectifier_square}
it holds that $\hiddenLength(\Phi) = M$,

\item
\label{item7:lemma:approximation_of_rectifier_square} it holds that $\paramANN(\Phi) = 20M-7$, and

\item
\label{item8:lemma:approximation_of_rectifier_square} it holds that $\infnorm{\vectorNN(\Phi)} \le \max\{4, {\constantR}^2\}$\ifnocf.
\end{enumerate}
\cfout[.]
\end{corollary}
\begin{cproof}{lemma:approximation_of_rectifier_square}
\Nobs that \cref{cor:approximation_of_square_on_segment} (applied with $M \with M$, $(A_k)_{k \in \N} \with (A_k)_{k \in \N}$, $\mathbb{A} \with \mathbb{A}$, $\mathbb{B} \with \mathbb{B}$, $(C_k)_{k \in \N} \with (C_k)_{k \in \N}$, $(c_k)_{k \in \N} \with (c_k)_{k \in \N}$, $\Phi \with \Psi$ in the notation of \cref{cor:approximation_of_square_on_segment}) assures that
\begin{enumerate}[label=(\Roman *)]
\item
\label{item1:lemma:approximation_of_rectifier_square_local} it holds that $\functionANN(\Psi) \in C(\R, \R)$,

\item
\label{item2:lemma:approximation_of_rectifier_square_local} it holds for all $x \in \R \backslash [0,1]$ that $(\functionANN(\Psi))(x) = \Rect(x)$, and

\item
\label{item3:lemma:approximation_of_rectifier_square_local} it holds for all $x \in [0,1]$ that $\abs{x^2 - (\functionANN(\Psi))(x)} \le 4^{-M-1}$\ifnocf.
\end{enumerate}

\noindent
\cfload[.]Next \nobs that \cref{prop:ANNcomposition_elementary_properties} and \cref{lem:ANN:affine} imply that for all $x \in \R$ it holds that $\functionANN(\Phi) \in C(\R, \R)$ and
\begin{equation}\label{eqn:approximation_of_rectifier_square_local_realization}
\begin{split}
(\functionANN(\Phi))(x) & = (\functionANN(\AffineANN_{{\constantR}^2, 0} \compANNbullet \Psi \compANNbullet \AffineANN_{{\constantR}^{-1}, 0}))(x) = \big(\functionANN(\AffineANN_{{\constantR}^2, 0})\big)\!\Big(\big(\functionANN(\Psi)\big)\!\big((\functionANN(\AffineANN_{{\constantR}^{-1}, 0}))(x)\big)\Big) \\
& = (\functionANN(\AffineANN_{{\constantR}^2, 0}))\!\pa{(\functionANN(\Psi))\!\pa{{\constantR}^{-1} x}} \! = {\constantR}^2 \!\pb{(\functionANN(\Psi))\!\pa{{\constantR}^{-1} x}}\!\ifnocf.
\end{split}
\end{equation}

\noindent
\cfload[.]This establishes \cref{item1:lemma:approximation_of_rectifier_square}. Moreover, \nobs that \cref{eqn:approximation_of_rectifier_square_local_realization}, \cref{item1:lemma:approximation_of_rectifier_square_local}, \cref{item2:lemma:approximation_of_rectifier_square_local}, and the fact that for all $x \in (-\infty, 0]$ it holds that ${\constantR}^{-1} x \in (-\infty, 0]$ ensure that for all $x \in (-\infty, 0]$ it holds that
\begin{equation}\label{eqn:lemma:approximation_of_rectifier_square_left}
\begin{split}
\pabs{[\Rect(x)]^2 - (\functionANN(\Phi))(x)} \! & = \! \pabs{[\Rect(x)]^2 - {\constantR}^2 \!\pb{(\functionANN(\Psi))\!\pa{{\constantR}^{-1}x}}} \\
& = \! \pabs{[\Rect(x)]^2 - {\constantR}^2 \Rect\!\pa{{\constantR}^{-1} x}} \! = 0.
\end{split}
\end{equation}

\noindent
This establishes \cref{item2:lemma:approximation_of_rectifier_square}. In the next step we \nobs that \cref{item2:lemma:approximation_of_rectifier_square_local}, \cref{eqn:approximation_of_rectifier_square_local_realization}, and the fact that for all $x \in [\constantR, \infty)$ it holds that ${\constantR}^{-1} x \in [1, \infty)$ demonstrate that for all $x \in [\constantR, \infty)$ it holds that
\begin{equation}
0 \le (\functionANN(\Phi))(x) = {\constantR}^2 \!\pb{(\functionANN(\Psi))\!\pa{{\constantR}^{-1} x}} = {\constantR}^2 \Rect\!\pa{{\constantR}^{-1} x} = \constantR x \le x^2 = \abs{\Rect(x)}^2.
\end{equation}

\noindent
The triangle inequality and the assumption that $q \in (2, \infty)$ therefore ensure that for all $x \in [\constantR, \infty)$ it holds that
\begin{equation}\label{eqn:lemma:approximation_of_rectifier_square_right}
\begin{split}
\pabs{[\Rect(x)]^2 - (\functionANN(\Phi))(x)} \! & = \abs{\Rect(x)}^2 - (\functionANN(\Phi))(x) \le \abs{\Rect(x)}^2 \\
& = \abs{x}^2 = \abs{x}^q \abs{x}^{2-q} \le \abs{x}^q {\constantR}^{2-q} = \abs{\Rect(x)}^q {\constantR}^{2-q}.
\end{split}
\end{equation}

\noindent
This establishes \cref{item4:lemma:approximation_of_rectifier_square}. Next \nobs that \cref{item3:lemma:approximation_of_rectifier_square_local}, \cref{eqn:approximation_of_rectifier_square_local_realization}, and the fact that for all $x \in [0, \constantR]$ it holds that ${\constantR}^{-1} x \in [0, 1]$ demonstrate that for all $x \in [0, \constantR]$ it holds that
\begin{equation}
\begin{split}
\pabs{[\Rect(x)]^2 - (\functionANN(\Phi))(x)} \! & = \! \pabs{x^2 - {\constantR}^2 \!\pb{(\functionANN(\Psi))\!\pa{{\constantR}^{-1} x}}} \\
& = {\constantR}^2 \! \pabs{[{\constantR}^{-1} x]^2 - (\functionANN(\Psi))\!\pa{{\constantR}^{-1} x}} \! \le 4^{-M-1} {\constantR}^2.
\end{split}
\end{equation}

\noindent
This establishes \cref{item3:lemma:approximation_of_rectifier_square}. Next \nobs that \cref{eqn:lemma:approximation_of_rectifier_square_weights_biases,eqn:lemma:approximation_of_rectifier_square_Psi} show that
\begin{equation}
{\constantR}^2 C_M = \! \pa{-2^{1-2M} {\constantR}^2 \quad \,\,\,\,\,\, 2^{2-2M} {\constantR}^2 \quad \,\,\,\,\,\, -2^{1-2M} {\constantR}^2 \quad \,\,\,\,\,\, {\constantR}^2 \, } \! \in \R^{1 \times 4}
\end{equation}

\noindent
and
\begin{equation}
\begin{split}
\Phi & = \AffineANN_{{\constantR}^2, 0} \compANNbullet \Psi \compANNbullet \AffineANN_{{\constantR}^{-1}, 0} = \begin{cases}
\pa{({\constantR}^{-1} \mathbb{A}, \mathbb{B}), ({\constantR}^2 C_1, 0)} & \colon M = 1\\
\pa{({\constantR}^{-1} \mathbb{A}, \mathbb{B}), (A_1, \mathbb{B}),\dots,(A_{M-1}, \mathbb{B}),
({\constantR}^2 C_M, 0)} & \colon M > 1.
\end{cases}
\end{split}
\end{equation}

\noindent
Combining this with \cref{eqn:lemma:approximation_of_rectifier_square_weights_biases} implies that $\dims(\Phi) = (1, 4, \ldots, 4, 1) \in \N^{M+2}$, $\hiddenLength(\Phi) = M$, $\paramANN(\Phi) = 20M - 7$, and $\infnorm{\vectorNN(\Phi)} \le \max\{4, {\constantR}^2\}$ \cfload. This establishes \cref{item5:lemma:approximation_of_rectifier_square,item6:lemma:approximation_of_rectifier_square,item7:lemma:approximation_of_rectifier_square,item8:lemma:approximation_of_rectifier_square}.
\end{cproof}

\subsection{ANN approximations for shifted squared rectifier functions}
\label{subsec:ann_approximation_for_shifted_squared_rectifier}

\cfclear
\begin{lemma}\label{lemma:realization_absolute-c}
Let $\fa \in \R$, $\shiftedabsnetwork, \Phi, \Psi \in \ANNs$ satisfy
\begin{equation}
\shiftedabsnetwork = \pa{\pa{\begin{pmatrix}
1 \\
-1
\end{pmatrix}\!, \begin{pmatrix}
0 \\
0
\end{pmatrix}}\!, \begin{pmatrix}
\begin{pmatrix}
1 & 1
\end{pmatrix}\!, (-\fa)
\end{pmatrix}} \in \pa{\pa{\R^{2 \times 1} \times \R^2} \times \pa{\R^{1 \times 2} \times \R^1}}\!,
\end{equation}
$\inDimANN(\Phi) = 1$, and $\Psi = \Phi \compANNbullet \shiftedabsnetwork$ \cfload. Then 
\begin{enumerate}[label=(\roman *)]
\item
\label{item1:lemma:realization_absolute-c} it holds that $\dims(\Psi) = \!\pa{1, 2, \dimANNlevel_1(\Phi), \ldots,\dimANNlevel_{\lengthANN(\Phi)}(\Phi)} \! \in \N^{\lengthANN(\Phi)+2}$,

\item
\label{item2:lemma:realization_absolute-c} it holds for all $x \in \R$ that $(\functionANN(\shiftedabsnetwork))(x) = \abs{x} - \fa$,

\item
\label{item3:lemma:realization_absolute-c} it holds for all $x \in \R$ that $(\functionANN(\Psi))(x) = (\functionANN(\Phi))(\abs{x} - \fa)$, and

\item
\label{item4:lemma:realization_absolute-c} it holds that $\norm{\vectorNN(\Psi)}_{\infty} \le (\abs{\fa}+1) \max\!\pc{1, \norm{\vectorNN(\Phi)}_{\infty}}$\ifnocf.
\end{enumerate}
\cfout[.]
\end{lemma}
\begin{cproof}{lemma:realization_absolute-c}
Throughout this proof let $L \in \N$, $l_0, l_1, \ldots, l_L \in \N$ satisfy $(l_0, l_1, \ldots, l_L) = \dims(\Phi)$ and let $W_k \in \R^{l_k \times l_{k-1}}$, $k \in \{1, 2, \ldots, L\}$, and $B_k \in \R^{l_k}$, $k \in \{1, 2, \ldots, L\}$, satisfy $\Phi = ((W_1, B_1), \allowbreak (W_2, B_2), \ldots, (W_L, B_L))$. \Nobs that $\dims(\shiftedabsnetwork) = (1, 2, 1) \in \N^3$. \cref{prop:ANNcomposition_elementary_properties} therefore ensures that $\dims(\Psi) = \dims(\Phi \compANNbullet \shiftedabsnetwork) = (1, 2, \dimANNlevel_1(\Phi), \allowbreak \ldots, \allowbreak \dimANNlevel_{\lengthANN(\Phi)}(\Phi)) \! \in \N^{\lengthANN(\Phi)+2}$. This establishes \cref{item1:lemma:realization_absolute-c}. Next \nobs that for all $x \in \R$ it holds that
\begin{equation}\label{eqn:lemma:realization_absolute-c_shiftedabs}
(\functionANN(\shiftedabsnetwork))(x) = \begin{pmatrix}
1 & 1
\end{pmatrix} \!\begin{pmatrix}
\Rect(x+0) \\
\Rect(-x+0)
\end{pmatrix} \! -\fa = \Rect(x)+\Rect(-x)-\fa = \abs{x}-\fa\ifnocf.
\end{equation}

\noindent
\cfload[.]This establishes \cref{item2:lemma:realization_absolute-c}. Moreover, \nobs that \cref{eqn:lemma:realization_absolute-c_shiftedabs} and \cref{prop:ANNcomposition_elementary_properties} assure that for all $x \in \R$ it holds that
\begin{equation}
(\functionANN(\Psi))(x) = (\functionANN(\Phi \compANNbullet \shiftedabsnetwork))(x) = (\functionANN(\Phi))((\functionANN(\shiftedabsnetwork))(x)) = (\functionANN(\Phi))(\abs{x}-\fa).
\end{equation}

\noindent
This establishes \cref{item3:lemma:realization_absolute-c}. In addition, \nobs that
\begin{equation}
\Psi = \Phi \compANNbullet \shiftedabsnetwork = \! \pa{\pa{\begin{pmatrix}
1 \\
-1
\end{pmatrix}\!, \begin{pmatrix}
0 \\
0
\end{pmatrix}}\!, \begin{pmatrix}
W_1 \begin{pmatrix}
1 & 1
\end{pmatrix}\!, W_1 (-\fa) +B_1
\end{pmatrix}, (W_2, B_2), \ldots, (W_L, B_L)}\!.
\end{equation}

\noindent
The fact that for all $\fW=(w_{i})_{i \in\{1,2,\dots,l_1\}}\in\R^{l_1 \times 1}$, $\fB = (b_1,b_2,\dots,b_{l_1}) \in \R^{l_1}$ it holds that
\begin{equation}
\fW \begin{pmatrix}
1 & 1
\end{pmatrix}\! = \!\pmat{
w_{1} & w_{1}\\
w_{2} & w_{2}\\
\vdots & \vdots\\
w_{l_1} & w_{l_1}}\!\in\R^{l_1\times 2} \qquad \text{and} \qquad \fW \, (-\fa) + \fB = \!\pmat{-\fa w_1+b_1\\-\fa w_2+b_2\\\vdots\\-\fa w_{l_1} + b_{l_1}} \! \in \R^{l_1}
\end{equation}

\noindent
hence demonstrates that
\begin{equation}
\begin{split}
\norm{\vectorNN(\Psi)}_{\infty} & \le \max\!\pc{1, \norm{\vectorNN(\Phi)}_{\infty}, (\abs{\fa}+1) \norm{\vectorNN(\Phi)}_{\infty}} = \max\!\pc{1, (\abs{\fa}+1) \norm{\vectorNN(\Phi)}_{\infty}} \\
& \le \max\!\pc{(\abs{\fa}+1), (\abs{\fa}+1) \norm{\vectorNN(\Phi)}_{\infty}} = (\abs{\fa}+1) \max\!\pc{1, \norm{\vectorNN(\Phi)}_{\infty}}\!\ifnocf.
\end{split}
\end{equation}
\cfload. This establishes \cref{item4:lemma:realization_absolute-c}.
\end{cproof}


\cfclear
\begin{corollary}\label{cor:approximation_of_rectifier_square}
Let $\fa \in [0, \infty)$, $M \in \N \cap [2, \infty)$, $\constantR \in [1, \infty)$, $q \in (2, \infty)$, $(A_k)_{k \in \N} \subseteq \R^{4 \times 4}$, $\mathbb{A}, \mathbb{B} \in \R^{4 \times 1}$, $(C_k)_{k \in \N} \subseteq \R^{1 \times 4}$, $(c_k)_{k \in \N} \subseteq \R$ satisfy for all $k \in \N$ that
\begin{equation}
A_k = \! \begin{pmatrix}
2 & -4 & 2 & 0\\
2 & -4 & 2 & 0\\
2 & -4 & 2 & 0\\
-c_k & 2c_k & -c_k & 1
\end{pmatrix}\!,
\quad
\mathbb{A} = \! \begin{pmatrix}
1 \\
1 \\
1 \\
1
\end{pmatrix}\!,
\quad
\mathbb{B} = \! \begin{pmatrix}
0\\
-\frac{1}{2}\\
-1\\
0
\end{pmatrix}\!,
\quad
C_k = \!\begin{pmatrix}
-c_k & 2c_k & -c_k & 1
\end{pmatrix}\!,
\end{equation}

\noindent
and $c_k = 2^{1-2k}$ and let $\Phi, \shiftedabsnetwork, \Psi \in \ANNs$ satisfy
\begin{equation}
\Phi = \pa{({\constantR}^{-1} \mathbb{A}, \mathbb{B}), (A_1, \mathbb{B}), \dots, (A_{M-1}, \mathbb{B}), ({\constantR}^2 C_M, 0)}\!,
\end{equation}
\begin{equation}
\shiftedabsnetwork = \pa{\pa{\begin{pmatrix}
1 \\
-1
\end{pmatrix}\!, \begin{pmatrix}
0 \\
0
\end{pmatrix}}\!, \begin{pmatrix}
\begin{pmatrix}
1 & 1
\end{pmatrix}\!, (-\fa)
\end{pmatrix}}\!,
\end{equation}

\noindent
and $\Psi = \Phi \compANNbullet \shiftedabsnetwork$ \cfload. Then
\begin{enumerate}[label=(\roman *)]
\item
\label{item1:cor:approximation_of_rectifier_square} it holds that $\functionANN(\Psi)\in C(\R,\R)$,

\item
\label{item2:cor:approximation_of_rectifier_square} it holds that $\dims(\Psi) = (1, 2, 4, 4, \ldots, 4, 1) \in \N^{M+3}$,

\item
\label{item3:cor:approximation_of_rectifier_square} it holds that $\hiddenLength(\Psi) = M + 1$,

\item
\label{item4:cor:approximation_of_rectifier_square} it holds that $\paramANN(\Psi) = 20M+1$,

\item
\label{item5:cor:approximation_of_rectifier_square} it holds that $\infnorm{\vectorNN(\Psi)} \le (\abs{\fa}+1) \max\{4, {\constantR}^2\}$,

\item
\label{item6:cor:approximation_of_rectifier_square} it holds for all $x \in \R$ that $(\functionANN(\Psi))(x) = (\functionANN(\Psi))(-x)$,

\item
\label{item7:cor:approximation_of_rectifier_square} it holds for all $x \in \R$ with $\abs{x} \le \fa$ that $\pabs{[\Rect(\abs{x}-\fa)]^2 - (\functionANN(\Psi))(x)} = 0$,

\item
\label{item8:cor:approximation_of_rectifier_square} it holds for all $x \in \R$ with $\fa \le \abs{x} \le \constantR + \fa$ that $\pabs{[\Rect(\abs{x}-\fa)]^2 - (\functionANN(\Psi))(x)} \le 4^{-M-1} {\constantR}^2$, and

\item
\label{item9:cor:approximation_of_rectifier_square} it holds for all $x \in \R$ with $\abs{x} \ge \constantR + \fa$ that $\pabs{[\Rect(\abs{x}-\fa)]^2 - (\functionANN(\Psi))(x)} \le  [\abs{x}-\fa]^q {\constantR}^{2-q}$\ifnocf.
\end{enumerate}
\cfout[.]
\end{corollary}
\begin{cproof}{cor:approximation_of_rectifier_square}
\Nobs that \cref{lemma:approximation_of_rectifier_square} (applied with $M \with M$, $\constantR \with \constantR$, $q \with q$, $(A_k)_{k \in \N} \with (A_k)_{k \in \N}$, $\mathbb{A} \with \mathbb{A}$, $\mathbb{B} \with \mathbb{B}$, $(C_k)_{k \in \N} \with (C_k)_{k \in \N}$, $(c_k)_{k \in \N} \with (c_k)_{k \in \N}$, $\Phi \with \Phi$ in the notation of \cref{lemma:approximation_of_rectifier_square}) implies that
\begin{enumerate}[label=(\Roman *)]
\item
\label{item1:cor:approximation_of_rectifier_square_local} it holds that $\functionANN(\Phi)\in C(\R,\R)$,

\item
\label{item2:cor:approximation_of_rectifier_square_local} it holds for all $x \in (-\infty, 0]$ that $\pabs{[\Rect(x)]^2 - (\functionANN(\Phi))(x)}=0$,

\item
\label{item3:cor:approximation_of_rectifier_square_local} it holds for all $x \in [0, \constantR]$ that $\pabs{[\Rect(x)]^2 - (\functionANN(\Phi))(x)} \le 4^{-M-1} {\constantR}^2$,

\item
\label{item4:cor:approximation_of_rectifier_square_local} it holds for all $x \in [\constantR, \infty)$ that $\pabs{[\Rect(x)]^2 - (\functionANN(\Phi))(x)} \le \abs{\Rect(x)}^q {\constantR}^{2-q}$,

\item
\label{item5:cor:approximation_of_rectifier_square_local} it holds that $\dims(\Phi) = (1, 4, \ldots, 4, 1) \in \N^{M+2}$, and

\item
\label{item6:cor:approximation_of_rectifier_square_local} it holds that $\infnorm{\vectorNN(\Phi)} \le \max\{4, {\constantR}^2\}$\ifnocf.
\end{enumerate}
\cfload[.]Next \nobs that \cref{lemma:realization_absolute-c} (applied with $\fa \with \fa$, $\shiftedabsnetwork \with \shiftedabsnetwork$, $\Phi \with \Phi$, $\Psi \with \Psi$ in the notation of \cref{lemma:realization_absolute-c}), \cref{item5:cor:approximation_of_rectifier_square_local}, and \cref{item6:cor:approximation_of_rectifier_square_local} ensure that 
\begin{equation}\label{eqn:cor:approximation_of_rectifier_square_dims_length}
\dims(\Psi) = \!\pa{1, 2, \dimANNlevel_1(\Phi), \ldots,\dimANNlevel_{\lengthANN(\Phi)}(\Phi)} \! = (1, 2, \underbrace{4, \ldots, 4}_{M}, 1) \in \N^{M+3}, \qquad \hiddenLength(\Psi) = M + 1,
\end{equation}
\begin{equation}\label{eqn:cor:approximation_of_rectifier_square_param}
\paramANN(\Psi) = 2(1+1) + 4(2+1) + \underbrace{4(4+1) +\ldots + 4(4+1)}_{M-1} + 1(4+1) = 20M+1,
\end{equation}

\noindent
and
\begin{equation}
\norm{\vectorNN(\Psi)}_{\infty} \le (\abs{\fa}+1) \max\!\pc{1, \norm{\vectorNN(\Phi)}_{\infty}}\! \le (\abs{\fa}+1) \max\{4, {\constantR}^2\}\ifnocf.
\end{equation}

\noindent
\cfload[.]This establishes \cref{item1:cor:approximation_of_rectifier_square,item2:cor:approximation_of_rectifier_square,item3:cor:approximation_of_rectifier_square,item4:cor:approximation_of_rectifier_square,item5:cor:approximation_of_rectifier_square}. Next \nobs that \cref{lemma:realization_absolute-c} (applied with $\fa \with \fa$, $\shiftedabsnetwork \with \shiftedabsnetwork$, $\Phi \with \Phi$, $\Psi \with \Psi$ in the notation of \cref{lemma:realization_absolute-c}) assures that for all $x \in \R$ it holds that 
\begin{equation}\label{eqn:cor:approximation_of_rectifier_square_shiftedabsnetwork}
(\functionANN(\shiftedabsnetwork))(x) = \abs{x} - \fa
\end{equation}

\noindent
and
\begin{equation}
(\functionANN(\Psi))(x) = (\functionANN(\Phi))(\abs{x} - \fa) = (\functionANN(\Phi))(\abs{-x} - \fa) = (\functionANN(\Psi))(-x).
\end{equation}

\noindent
This establishes \cref{item6:cor:approximation_of_rectifier_square}. Furthermore, \nobs that \cref{eqn:cor:approximation_of_rectifier_square_shiftedabsnetwork} shows that for all $x \in [-\fa, \fa]$ it holds that $(\functionANN(\shiftedabsnetwork))(x) = \abs{x} - \fa \le 0$. Combining this with \cref{item2:cor:approximation_of_rectifier_square_local} proves that for all $x \in [-\fa, \fa]$ it holds that
\begin{equation}
\pabs{[\Rect(\abs{x}-\fa)]^2-(\functionANN(\Psi))(x)}\! = \! \pabs{[\Rect\!\pa{(\functionANN(\shiftedabsnetwork))(x)}]^2 - (\functionANN(\Phi))((\functionANN(\shiftedabsnetwork))(x))} \! = 0.
\end{equation}

\noindent
This establishes \cref{item7:cor:approximation_of_rectifier_square}. Moreover, \nobs that \cref{eqn:cor:approximation_of_rectifier_square_shiftedabsnetwork} demonstrates that for all $x \in \R$ with $\fa \le \abs{x} \le \constantR + \fa$ it holds that $(\functionANN(\shiftedabsnetwork))(x) = \abs{x} - \fa \in [0, \constantR]$. This and \cref{item3:cor:approximation_of_rectifier_square_local} ensure that for all $x \in \R$ with $\fa \le \abs{x} \le \constantR + \fa$ it holds that
\begin{equation}
\pabs{[\Rect(\abs{x}-\fa)]^2 - (\functionANN(\Psi))(x)} = \! \pabs{[\Rect((\functionANN(\shiftedabsnetwork))(x))]^2 - (\functionANN(\Phi))((\functionANN(\shiftedabsnetwork))(x))} \le 4^{-M-1} {\constantR}^2.
\end{equation}

\noindent
This establishes \cref{item8:cor:approximation_of_rectifier_square}. In addition, \nobs that \cref{eqn:cor:approximation_of_rectifier_square_shiftedabsnetwork} proves that for all $x \in \R$ with $\abs{x} \ge \allowbreak \constantR + \fa$ it holds that $(\functionANN(\shiftedabsnetwork))(x) = \abs{x} - \fa \in [\constantR, \infty)$. Item~\ref{item4:cor:approximation_of_rectifier_square_local} hence shows that for all $x \in \R$ with $\abs{x} \ge \allowbreak \constantR + \fa$ it holds that
\begin{equation}
\begin{split}
\pabs{[\Rect(\abs{x}-\fa)]^2 - (\functionANN(\Psi))(x)} & = \! \pabs{[\Rect\!\pa{(\functionANN(\shiftedabsnetwork))(x)}]^2 - (\functionANN(\Phi))((\functionANN(\shiftedabsnetwork))(x))} \\
& \le \abs{\Rect\!\pa{(\functionANN(\shiftedabsnetwork))(x)}}^q {\constantR}^{2-q} = \abs{\Rect\!\pa{\abs{x}-\fa}}^q {\constantR}^{2-q}.
\end{split}
\end{equation}

\noindent
This establishes \cref{item9:cor:approximation_of_rectifier_square}.
\end{cproof}



\subsection{Lower and upper bounds for integrals of certain specific high-dimensional functions}
\label{subsec:bounds_for_the_norm_of_specific_function_lower}


\cfclear
\begin{lemma}\label{lemma:lower_bound_gaussian_tail_from_Klenke}
Let $s \in (0, \infty)$ \cfload. Then
\begin{equation}
\int_s^{\infty} e^{- \frac{1}{2} x^2} \, dx \ge \frac{e^{- \frac{1}{2} s^2}}{s + s^{-1}}\ifnocf.
\end{equation}
\cfout[.]
\end{lemma}
\begin{cproof}{lemma:lower_bound_gaussian_tail_from_Klenke}
\Nobs that the integration by parts formula ensures that
\begin{equation}
\begin{split}
\int_s^{\infty} e^{- \frac{1}{2} x^2} \, dx & = \int_s^{\infty} -x^{-1} \big[e^{- \frac{1}{2} x^2}\big]' \, dx = \lim_{T \to \infty} \pa{\! \pb{- x^{-1} e^{- \frac{1}{2} x^2}}_{x = s}^{x = T}} - \int_{s}^{\infty} \! \pb{x^{-2} e^{- \frac{1}{2} x^2}} \! dx \\
& = s^{-1} e^{- \frac{1}{2} s^2} - \int_{s}^{\infty} \! \pb{x^{-2} e^{- \frac{1}{2} x^2}} \! dx \ge s^{-1} e^{-  \frac{1}{2} s^2} - s^{-2} \int_{s}^{\infty} e^{- \frac{1}{2} x^2} \, dx.
\end{split}
\end{equation}

\noindent
Hence, we obtain that
\begin{equation}
\begin{split}
\int_s^{\infty} e^{- \frac{1}{2} x^2} \, dx & = \! \pb{\frac{s^2}{1 + s^2}} \! \pb{1+\frac{1}{s^2}} \! \pb{\int_s^{\infty} e^{- \frac{1}{2} x^2} \, dx} \\
& = \! \pb{\frac{s^2}{1 + s^2}} \! \pb{\int_s^{\infty} e^{- \frac{1}{2} x^2} \, dx + \frac{1}{s^2} \int_{s}^{\infty} e^{- \frac{1}{2} x^2} \, dx} \\
& \ge \! \pb{\frac{s^2}{1 + s^2}} \! \pb{\frac{e^{- \frac{1}{2} s^2}}{s}} \! = \frac{e^{- \frac{1}{2} s^2}}{s + s^{-1}}\ifnocf.
\end{split}
\end{equation}
\end{cproof}

\cfclear
\begin{lemma}\label{lemma:lower_bound_gaussian_tail}
Let $\expconst, s \in (0, \infty)$ \cfload. Then
\begin{equation}
\int_s^{\infty} e^{- \expconst x^2} \, dx \ge \frac{e^{- \expconst s^2}}{s^{-1} + 2 \expconst s}\ifnocf.
\end{equation}
\cfout[.]
\end{lemma}
\begin{cproof}{lemma:lower_bound_gaussian_tail}
\Nobs that the integral transformation theorem and \cref{lemma:lower_bound_gaussian_tail_from_Klenke} (applied with $s \with s \sqrt{2 \sigma}$ in the notation of \cref{lemma:lower_bound_gaussian_tail_from_Klenke}) ensure that
\begin{equation}
\int_s^{\infty} e^{- \expconst x^2} \, dx = \frac{1}{\sqrt{2 \sigma}} \int_{s \sqrt{2 \sigma}}^{\infty} e^{- \frac{1}{2} x^2} \, dx \ge \frac{1}{\sqrt{2 \sigma}} \! \pb{\frac{e^{-\frac{1}{2} (s\sqrt{2 \sigma})^2}}{s \sqrt{2 \sigma} + (s \sqrt{2 \sigma})^{-1}}} \! = \frac{e^{- \expconst s^2}}{s^{-1} + 2 \expconst s}\ifnocf.
\end{equation}
\cfload[.]
\end{cproof}

\cfclear
\begin{lemma}\label{lemma:small_facts_sec_5}
Let $d \in \N$. Then
\begin{equation}
\frac{\sqrt{2d}(2d+1)}{4d^2(4d^2+6d+1)} \! \pb{\frac{2}{\pi}}^{\nicefrac{1}{2}} \! e^{-1-\frac{1}{4d}} \ge 50^{-1} d^{\nicefrac{-5}{2}}.
\end{equation}
\end{lemma}
\begin{cproof}{lemma:small_facts_sec_5}
\Nobs that $48 d^2 -28 d \ge 20d^2 \ge 13$. This implies that $4d^2 + 6d + 1 \le (\nicefrac{25}{13})(4d^2 + 2d) = (\nicefrac{50}{13})d(2d + 1)$. The fact that $13 \ge 2\sqrt{\pi} e^{\nicefrac{5}{4}}$ and the fact that $-1-\frac{1}{4d} \ge -\frac{5}{4}$ hence ensure that
\begin{equation}
\frac{\sqrt{2d}(2d+1)}{4d^2(4d^2+6d+1)} \! \pb{\frac{2}{\pi}}^{\nicefrac{1}{2}} \! e^{-1-\frac{1}{4d}} \ge \frac{\sqrt{2d}}{4d^2} \! \pb{\frac{13}{50 d}} \! \pb{\frac{2}{\pi}}^{\nicefrac{1}{2}} \! e^{\nicefrac{-5}{4}} \ge \frac{\sqrt{2d}}{4d^2} \! \pb{\frac{2\sqrt{2}}{50 d}} \! = 50^{-1} d^{\nicefrac{-5}{2}}.
\end{equation}
\end{cproof}

\cfclear
\begin{lemma}\label{lemma:bounds_for_the_norm_of_specific_function}
Let $d \in \N$ and let $\varphi \colon \R^d \to \R$ and $g \colon \R^d \to \R$ satisfy for all $x = (x_1, x_2, \ldots, x_d) \in \R^d$ that $\varphi(x) = (2 \pi)^{\nicefrac{-d}{2}} \exp(-\frac{1}{2}(\smallsum_{j=1}^d \abs{x_j}^2))$ and $g(x) = \smallsum_{j=1}^d [\max\{\abs{x_j}-\sqrt{2 d}, 0\}]^2$ \cfload. Then
\begin{equation}
(50)^{-1} d^{\nicefrac{-3}{2}} e^{- d} \le \int_{\R^d} \abs{g(x)}^2 \varphi(x) \, dx \le 3 d^2 e^{- d}\ifnocf.
\end{equation}
\cfout[.]
\end{lemma}
\begin{cproof}{lemma:bounds_for_the_norm_of_specific_function}
Throughout this proof let $\Gamma \colon (0, \infty) \to (0, \infty)$ satisfy for all $x \in (0, \infty)$ that $\Gamma(x) = \int_{0}^{\infty} t^{x-1} e^{-t} \, dt$.  \Nobs that the fact that for all $k \in \N$, $\fa_1, \fa_2, \ldots, \fa_k \in \R$ it holds that
\begin{equation}
 \abs{\fa_1}^2+\abs{\fa_2}^2+\ldots+\abs{\fa_k}^2 \le (\abs{\fa_1}+\abs{\fa_2}+\ldots+\abs{\fa_k})^2 \le k(\abs{\fa_1}^2+\abs{\fa_2}^2+\ldots+\abs{\fa_k}^2)
\end{equation}

\noindent
ensures that for all $(x_1, x_2, \ldots, x_d) \in \R^d$ it holds that
\begin{equation}
\sum_{j=1}^d \!\big[\Rect\big(\abs{x_j}-\sqrt{2 d}\big)\big]^4\! \le \pb{\sum_{j=1}^d \!\big[\Rect\big(\abs{x_j}-\sqrt{2 d}\big)\big]^2}^2 \! \le d \! \pb{\sum_{j=1}^d \!\big[\Rect\big(\abs{x_j}-\sqrt{2 d}\big)\big]^4}\!\ifnocf.
\end{equation}

\noindent
\cfload[.]The fact that for all $k \in \N$ it holds that $\int_{\R^k} (2 \pi)^{\nicefrac{-k}{2}} e^{- \frac{1}{2} \norm{x}_2^2} \, dx = 1$ therefore demonstrates that
\begin{equation}
\begin{split}
& d \! \int_{\R} \!\big[\Rect\big(\abs{x}-\sqrt{2 d}\big)\big]^4 (2 \pi)^{\nicefrac{-1}{2}} e^{- \frac{1}{2} x^2} dx \\
& = d \! \int_{\R} \int_{\R} \ldots \int_{\R} \!\big[\Rect\big(\abs{x_1}-\sqrt{2 d}\big)\big]^4 (2 \pi)^{\nicefrac{-d}{2}}e^{- \frac{1}{2} [\sum_{j=1}^d \abs{x_j}^2]} \, dx_d \ldots dx_2 \, dx_1 \\
& = \sum_{j=1}^d \int_{\R} \int_{\R} \ldots \int_{\R} \!\big[\Rect\big(\abs{x_j}-\sqrt{2 d}\big)\big]^4 (2 \pi)^{\nicefrac{-d}{2}}e^{- \frac{1}{2} [\sum_{j=1}^d \abs{x_j}^2]} \, dx_d \ldots dx_2 \, dx_1 \\
& = \int_{\R} \int_{\R} \ldots \int_{\R} \! \pb{\sum_{j=1}^d \!\big[\Rect\big(\abs{x_j}-\sqrt{2 d}\big)\big]^4} \! (2 \pi)^{\nicefrac{-d}{2}}e^{- \frac{1}{2} [\sum_{j=1}^d \abs{x_j}^2]} \, dx_d \ldots dx_2 \, dx_1 \\
& \le \int_{\R} \int_{\R} \ldots \int_{\R} \! \pb{\sum_{j=1}^d \!\big[\Rect\big(\abs{x_j}-\sqrt{2 d}\big)\big]^2}^2 \! (2 \pi)^{\nicefrac{-d}{2}} e^{- \frac{1}{2} [\sum_{j=1}^d \abs{x_j}^2]} \, dx_d \ldots dx_2 \, dx_1 \\
& = \int_{\R^d} \abs{g(x)}^2 \varphi(x) \, dx
\end{split}
\end{equation}

\noindent
and
\begin{equation}
\begin{split}
& d^2 \! \int_{\R} \!\big[\Rect\big(\abs{x}-\sqrt{2 d}\big)\big]^4 (2 \pi)^{\nicefrac{-1}{2}} e^{- \frac{1}{2} x^2} dx \\
& = d^2 \! \int_{\R} \int_{\R} \ldots \int_{\R} \!\big[\Rect\big(\abs{x_1}-\sqrt{2 d}\big)\big]^4 (2 \pi)^{\nicefrac{-d}{2}} e^{- \frac{1}{2} [\sum_{j=1}^d \abs{x_j}^2]} \, dx_d \ldots dx_2 \, dx_1 \\
& = d \! \int_{\R} \int_{\R} \ldots \int_{\R} \! \pb{\sum_{j=1}^d \!\big[\Rect\big(\abs{x_j}-\sqrt{2 d}\big)\big]^4} \! (2 \pi)^{\nicefrac{-d}{2}} e^{- \frac{1}{2} [\sum_{j=1}^d \abs{x_j}^2]} \, dx_d \ldots dx_2 \, dx_1 \\
& \ge \int_{\R} \int_{\R} \ldots \int_{\R} \! \pb{\sum_{j=1}^d \!\big[\Rect\big(\abs{x_j}-\sqrt{2 d}\big)\big]^2}^2 \! (2 \pi)^{\nicefrac{-d}{2}} e^{- \frac{1}{2} [\sum_{j=1}^d \abs{x_j}^2]} \, dx_d \ldots dx_2 \, dx_1 \\
& = \int_{\R^d} \abs{g(x)}^2 \varphi(x) \, dx\ifnocf.
\end{split}
\end{equation}

\noindent
\cfload[.]Hence, we obtain that
\begin{equation}
\label{eqn:lemma:bounds_for_the_norm_of_specific_function_lower_upper}
\begin{split}
d \! \int_{\R} \!\big[\Rect\big(\abs{x}-\sqrt{2 d}\big)\big]^4  (2 \pi)^{\nicefrac{-1}{2}} e^{- \frac{1}{2} x^2} dx & \le \int_{\R^d} \abs{g(x)}^2 \varphi(x) \, dx \\
& \le d^2 \! \int_{\R} \!\big[\Rect\big(\abs{x}-\sqrt{2 d}\big)\big]^4 (2\pi)^{\nicefrac{-1}{2}} e^{- \frac{1}{2} x^2} dx.
\end{split}
\end{equation}

\noindent
Next \nobs that \cref{lemma:lower_bound_gaussian_tail} (applied with $\expconst \with \nicefrac{1}{2}$, $s \with (2d)^{\nicefrac{1}{2}}+(2d)^{\nicefrac{-1}{2}}$ in the notation of \cref{lemma:lower_bound_gaussian_tail}) and \cref{lemma:small_facts_sec_5} (applied with $d \with d$ in the notation of \cref{lemma:small_facts_sec_5}) ensure that
\begin{equation}\label{eqn:lemma:bounds_for_the_norm_of_specific_function_lower}
\begin{split}
& \int_{\R} \!\big[\Rect\big(\abs{x}-\sqrt{2 d}\big)\big]^4 (2\pi)^{\nicefrac{-1}{2}} e^{- \frac{1}{2} x^2} dx \\
& = 2 \! \pb{\int_{0}^{\infty} \!\big[\Rect\big(\abs{x}-\sqrt{2 d}\big)\big]^4 (2\pi)^{\nicefrac{-1}{2}} e^{- \frac{1}{2} x^2} dx} \! = \! \pb{\frac{2}{\pi}}^{\nicefrac{1}{2}} \! \pb{\int_{\sqrt{2 d}}^{\infty} \big[x-\sqrt{2 d}\big]^4 e^{- \frac{1}{2} x^2} dx} \\
& \ge \frac{1}{4 d^2} \! \pb{\frac{2}{\pi}}^{\nicefrac{1}{2}} \! \pb{\int_{(2d)^{\nicefrac{1}{2}}+(2d)^{\nicefrac{-1}{2}}}^{\infty}e^{- \frac{1}{2} x^2} dx} \! \ge \frac{1}{4 d^2} \! \pb{\frac{2}{\pi}}^{\nicefrac{1}{2}} \! \pb{\frac{[(2d)^{\nicefrac{1}{2}}+(2d)^{\nicefrac{-1}{2}}] e^{-\frac{1}{2}(2d + 2 + (2d)^{-1})}}{1 + (2d + 2 +(2d)^{-1})}} \\
& = e^{-d} \! \pb{\frac{\sqrt{2d}(2d+1)}{4d^2(4d^2+6d+1)} \! \pb{\frac{2}{\pi}}^{\nicefrac{1}{2}} \! e^{-1-\frac{1}{4d}}} \ge 50^{-1} d^{\nicefrac{-5}{2}} e^{-d}.
\end{split}
\end{equation}

\noindent
Moreover, \nobs that the integral transformation theorem and \cref{lemma:GammaBeta} demonstrate that
\begin{equation}
\begin{split}
& \int_{\R} \!\big[\Rect\big(\abs{x}-\sqrt{2 d}\big)\big]^4 (2 \pi)^{\nicefrac{-1}{2}} e^{- \frac{1}{2} x^2} dx \\
& = 2 \! \pb{\int_{0}^{\infty} \!\big[\Rect\big(\abs{x}-\sqrt{2 d}\big)\big]^4 (2 \pi)^{\nicefrac{-1}{2}} e^{- \frac{1}{2} x^2} dx} \! = \! \pb{\frac{2}{\pi}}^{\nicefrac{1}{2}} \! \pb{\int_{\sqrt{2 d}}^{\infty} \big[x-\sqrt{2 d}\big]^4\! e^{- \frac{1}{2} x^2} dx} \\
& = \! \pb{\frac{2}{\pi}}^{\nicefrac{1}{2}} \! \pb{\int_{0}^{\infty} x^4 \, e^{- \frac{1}{2} (x+\sqrt{2d})^2} \, dx} \! = \! \pb{\frac{2}{\pi}}^{\nicefrac{1}{2}} \! e^{-d} \! \pb{\int_{0}^{\infty} x^4 \, e^{- \frac{1}{2} (x^2 + 2\sqrt{2d} x)} \, dx} \\
& \le \! \pb{\frac{2}{\pi}}^{\nicefrac{1}{2}} \! e^{-d} \! \pb{\int_{0}^{\infty} x^4 \, e^{- \frac{1}{2} x^2} \, dx} \! = \! \pb{\frac{4}{\sqrt{\pi}}} \! e^{-d} \! \pb{\int_{0}^{\infty} x^{\nicefrac{3}{2}} e^{-x} \, dx} \! = \! \pb{\frac{4}{\sqrt{\pi}}} \! e^{-d} \, \Gamma\!\pa{\frac{5}{2}} = 3 e^{-d}.
\end{split}
\end{equation}

\noindent
Combining this with \cref{eqn:lemma:bounds_for_the_norm_of_specific_function_lower_upper,eqn:lemma:bounds_for_the_norm_of_specific_function_lower} demonstrates that
\begin{equation}
50^{-1} d^{\nicefrac{-3}{2}} e^{- d} \le \int_{\R^d} \abs{g(x)}^2 \varphi(x) \, dx \le 3 d^2 e^{- d}.
\end{equation}
\end{cproof}




\subsection{ANN representations for multiplications with powers of real numbers}
\label{subsec:ann_representation_powers_of_reals}

\cfclear
\begin{lemma}\label{lemma:on_Realization_of_DNN_2}
Let $n \in \N$, $\lambda \in \R$, $\Phi, \idPowerrr, \Psi\in \ANNs$ satisfy $\idPowerrr = (\scalarMultANN{\lambda}{\idRelu_{\outDimANN(\Phi)}}) \compANNbullet \AffineANN_{\lambda \idMatrix_{\outDimANN(\Phi)}, 0}$ and $\Psi = (\powANN{\idPowerrr}{n}) \compANNbullet \Phi$ \cfload. Then
\begin{enumerate}[label=(\roman *)]
\item
\label{item1:lemma:on_Realization_of_DNN_2} it holds that $\inDimANN(\Psi) = \inDimANN(\Phi)$,

\item
\label{item2:lemma:on_Realization_of_DNN_2} it holds that $\hiddenLength(\Psi) = \hiddenLength(\Phi) + n$,

\item
\label{item3:lemma:on_Realization_of_DNN_2} it holds that $\paramANN(\Psi) \le 2 \paramANN(\Phi) + 6n \abs{\outDimANN(\Phi)}^2$,

\item
\label{item4:lemma:on_Realization_of_DNN_2} it holds that $\norm{\vectorNN(\Psi)}_{\infty} \le  \max\!\pc{1, \abs{\lambda}} \max\!\pc{ \abs{\lambda}, \norm{\vectorNN(\Phi)}_{\infty}}$, and

\item
\label{item5:lemma:on_Realization_of_DNN_2} it holds for all $x \in \R^{\inDimANN(\Psi)}$ that $(\functionANN(\Psi))(x) = \lambda^{2 n}(\functionANN(\Phi))(x)$\ifnocf.
\end{enumerate}
\cfout[.]
\end{lemma}
\begin{cproof}{lemma:on_Realization_of_DNN_2}
Throughout this proof let $d, l_0, l_1, l_2 \in \N$  satisfy $l_0 = l_2 = d = \outDimANN(\Phi)$ and $l_1 = 2d$, let $O_{\fn} \in \R^{\fn}$, $\fn \in \N$, satisfy for all $\fn \in \N$ that $O_{\fn} = 0$, and let $W_{k} \in \R^{l_k \times l_{k-1}}$, $k \in \{1, 2\}$, satisfy $\idRelu_d = ((W_{1}, O_{2d}), (W_{2}, O_d))$ (cf.\ \cref{lem:Relu:identity})\cfload. \Nobs that \cref{Lemma:PropertiesOfANNenlargementGeometry}, \cref{prop:ANNcomposition_elementary_properties}, and \cref{lem:ANNscalar} show that
\begin{equation}\label{eqn:length_dims}
\begin{split}
& \dims(\powANN{\idPowerrr}{n})=(d, 2d, 2d, \dots, 2d, d)\in\N^{n+2}, \qquad \hiddenLength(\powANN{\idPowerrr}{n}) = n, \qquad \hiddenLength(\Psi) = \hiddenLength(\Phi) + n, \\
& \text{and} \qquad \dims(\Psi) = (\dimANNlevel_0(\Phi), \dimANNlevel_1(\Phi), \ldots, \dimANNlevel_{\hiddenLength(\Phi)}(\Phi), \underbrace{2d, 2d, \ldots, 2d}_{n}, d) \in \N^{\lengthANN(\Phi)+n+1}\ifnocf.
\end{split}
\end{equation}

\noindent
Therefore, we obtain that
\begin{equation}\label{eqn:lemma:on_Realization_of_DNN_2_param}
\begin{split}
\paramANN(\Psi) & = \paramANN(\Phi) + \dimANNlevel_{\lengthANN(\Phi)}(\Phi)(\dimANNlevel_{\hiddenLength(\Phi)}(\Phi) + 1) + \underbrace{2d(2d+1) + \ldots + 2d(2d+1)}_{n-1} + d(2d+1) \\
& = \paramANN(\Phi) + \dimANNlevel_{\lengthANN(\Phi)}(\Phi)(\dimANNlevel_{\hiddenLength(\Phi)}(\Phi) + 1) + (n-1) (4d^2+2d) + (2d^2+d) \\
& \le 2 \paramANN(\Phi) + 6 d^2 n = 2 \paramANN(\Phi) + 6n \abs{\outDimANN(\Phi)}^2\ifnocf.
\end{split}
\end{equation}

\noindent
\cfload[.]Moreover, \nobs that \cref{ANNoperations:Composition} and the fact that for all $\alpha \in \R$, $\phi \in \ANNs$ it holds that $\scalarMultANN{\alpha}{\phi} = \AffineANN_{\alpha \idMatrix_{\outDimANN(\phi)}, 0} \allowbreak \compANNbullet \allowbreak \phi$ ensure that $\idPowerrr = ((\lambda W_1, O_{2d}), (\lambda W_2, O_d))$. Therefore, we obtain that
\begin{equation}\label{eqn:lemma:on_Realization_of_DNN_2:idpowerrr_n}
\powANN{\idPowerrr}{n} = ((\lambda W_1, O_{2d}), \underbrace{(\lambda^2 W_1 W_2, O_{2d}), \dots, (\lambda^2 W_1 W_2, O_{2d})}_{n-1}, (\lambda W_2, O_d)).
\end{equation}

\noindent
Next \nobs that the fact that $\idRelu_d = ((W_{1}, O_{2d}), (W_{2}, O_d))$, \cref{parallelisationSameLengthDef,eq:def:id:1,eq:def:id:2} demonstrate that $\infnorm{\vectorNN(((W_1, O_{2d})))} = \infnorm{\vectorNN(((W_2, O_d)))} = \infnorm{\vectorNN(((W_1 W_2, O_{2d})))} = 1$ \cfload. Combining this with \cref{eqn:lemma:on_Realization_of_DNN_2:idpowerrr_n} establishes that
\begin{equation}\label{eq:extension3}
\asinfnorm{\vectorNN\bigl(\powANN{\idPowerrr}{n}\bigr)} = \begin{cases}
\abs{\lambda} &\colon n=1 \\
\abs{\lambda} \max\!\pc{1, \abs{\lambda}} & \colon n>1.
\end{cases}
\end{equation}

\noindent
Furthermore, \nobs that the fact that $\idRelu_d = ((W_{1}, O_{2d}), (W_{2}, O_d))$, \cref{parallelisationSameLengthDef,eq:def:id:1,eq:def:id:2} show that
for all $k\in\N$, $\fW \in \R^{d\times k}$, $\fB \in \R^{d}$ it holds that
\begin{equation}
\infnorm{\vectorNN(((\lambda W_1 \fW, \lambda W_1 \fB + O_{2d})))} = \abs{\lambda} \infnorm{\vectorNN(((\fW, \fB)))}.
\end{equation}

\noindent
This, \cref{lem:composition_infnorm}, \cref{eqn:lemma:on_Realization_of_DNN_2:idpowerrr_n}, and \cref{eq:extension3} establish that
\begin{equation}\label{eqn:lemma:on_Realization_of_DNN_2_size_of_param}
\begin{split}
\infnorm{\vectorNN(\Psi)} & = \asinfnorm{\vectorNN\bigl((\powANN{\idPowerrr}{n}) \compANNbullet \Phi \bigr)} \leq \max\!\pc{\norm{\vectorNN(\powANN{\idPowerrr}{n})}_{\infty}, \norm{\vectorNN(\Phi)}_{\infty}, \abs{\lambda} \norm{\vectorNN(\Phi)}_{\infty}} \\
& \leq \max\bigl\{ \abs{\lambda} \max\!\pc{1, \abs{\lambda}}, \norm{\vectorNN(\Phi)}_{\infty}, \abs{\lambda} \norm{\vectorNN(\Phi)}_{\infty} \bigr\} \\
& = \max\!\pc{1, \abs{\lambda}} \max\!\pc{ \abs{\lambda}, \norm{\vectorNN(\Phi)}_{\infty}}\!\ifnocf.
\end{split}
\end{equation}

\noindent
\cfload[.]In addition, \nobs that \cref{prop:ANNcomposition_elementary_properties}, \cref{lem:Relu:identity},  \cref{lem:ANN:affine}, and \cref{lem:ANNscalar} demonstrate that for all $x \in \R^d$ it holds that
\begin{equation}
\begin{split}
(\functionANN(\idPowerrr))(x) & = (\functionANN((\scalarMultANN{\lambda}{\idRelu_d})\compANNbullet \AffineANN_{\lambda \idMatrix_d, 0}))(x) = (\functionANN(\scalarMultANN{\lambda}{\idRelu_d}))((\functionANN(\AffineANN_{\lambda \idMatrix_d, 0}))(x)) \\
& = (\functionANN(\scalarMultANN{\lambda}{\idRelu_d}))(\lambda x) = \lambda [(\functionANN(\idRelu_d))(\lambda x)] = \lambda [\lambda x] = \lambda^2 x.
\end{split}
\end{equation}

\noindent
Induction therefore shows that for all $x \in \R^d$ it holds that $(\functionANN(\powANN{\idPowerrr}{n}))(x) = \lambda^{2n} x$. Hence, we obtain that for all $x \in \R^{\inDimANN(\Psi)}$ it holds that
\begin{equation}\label{eqn:lemma:on_Realization_of_DNN_2_realization}
\functionANN(\Psi)(x) = (\functionANN((\powANN{\idPowerrr}{n}) \compANNbullet \Phi))(x) = (\functionANN(\powANN{\idPowerrr}{n}))((\functionANN(\Phi))(x)) = \lambda^{2n} (\functionANN(\Phi))(x).
\end{equation}

\noindent
Combining this with \cref{eqn:lemma:on_Realization_of_DNN_2_param,eqn:length_dims,eqn:lemma:on_Realization_of_DNN_2_size_of_param} establishes \cref{item1:lemma:on_Realization_of_DNN_2,item2:lemma:on_Realization_of_DNN_2,item3:lemma:on_Realization_of_DNN_2,item4:lemma:on_Realization_of_DNN_2,item5:lemma:on_Realization_of_DNN_2}.
\end{cproof}

\subsection{ANN approximations for certain specific high-dimensional functions}
\label{subsec:upper_bounds_anns_for_specific_class}




\cfclear
\begin{theorem}\label{thm:main2}
Let $d \in \N$, $M \in \N \cap [2, \infty)$, $\constantR \in [1, \infty)$, let $\varphi\colon \R^d \to \R$ and $g \colon \R^d \to \R$ satisfy for all $x = (x_1, x_2, \ldots, x_d) \allowbreak \in \R^d$ that $\varphi(x)=(2\pi)^{\nicefrac{-d}{2}}\exp(- \frac{1}{2} (\smallsum_{j=1}^d \abs{x_j}^2))$ and $g(x) = \sum_{j=1}^d [\max\{\abs{x_j}-\sqrt{2 d}, 0\}]^2\!$, and let $\fg \colon \R^d\rightarrow \R$ satisfy for all $x \in \R^d$ that $\fg(x)=[\int_{\R^d}\abs{g(y)}^2\varphi(y)\,dy]^{-1/2}g(x)$ \cfload. Then there exists $\Phi \in \ANNs$ such that
\begin{enumerate}[label=(\roman *)]
\item
\label{item1:thm2} it holds that $\functionANN(\Phi) \in C(\R^d, \R)$,

\item
\label{item2:thm2} it holds that $\hiddenLength(\Phi) = d + M + 1$,

\item
\label{item3:thm2} it holds that $\paramANN(\Phi) \le 42 d^2 M + 6d$,

\item
\label{item4:thm2} it holds that $\norm{\vectorNN(\Phi)}_{\infty} \le 12 d^{\nicefrac{3}{2}} \! \max\{4, {\constantR}^2\}$, and

\item
\label{item5:thm2} it holds that $\int_{\R^d} \abs{(\functionANN(\Phi))(x) - \fg(x)}^2 \varphi(x) \, dx \le 50 d^{\nicefrac{7}{2}} \! \pb{16^{-M-1} {\constantR}^4 + 105 {\constantR}^{-4}}\!$\ifnocf.
\end{enumerate}
\cfout[.]
\end{theorem}
\begin{cproof}{thm:main2}
Throughout this proof let $\Gamma \colon (0, \infty) \to (0, \infty)$ satisfy for all $x \in (0, \infty)$ that $\Gamma (x) = \int_0^{\infty} t^{x-1} e^{-t} \, dt$, let $\psi \in \ANNs$ satisfy that
\begin{enumerate}[label=(\Roman *)]
\item
\label{item1:thm2:approximation_of_rectifier_square_local} it holds that $\functionANN(\psi)\in C(\R,\R)$,

\item
\label{item2:thm2:approximation_of_rectifier_square_local} it holds that $\dims(\psi) = (1, 2, \underbrace{4, \ldots, 4}_{M}, 1) \in \N^{M+3}$,

\item
\label{item3:thm2:approximation_of_rectifier_square_local} it holds that $\norm{\vectorNN(\psi)}_{\infty} \le (\sqrt{2d}+1)\max\{4, {\constantR}^2\}$,

\item
\label{item4:thm2:approximation_of_rectifier_square_local} it holds for all $x \in \R$ that $(\functionANN(\psi))(x) = (\functionANN(\psi))(-x)$,

\item
\label{item5:thm2:approximation_of_rectifier_square_local} it holds for all $x \in \R$ with $\abs{x} \le \sqrt{2d}$ that $\Abs{[\Rect(\abs{x}-\sqrt{2 d})]^2 - (\functionANN(\psi))(x)} \! = 0$,

\item
\label{item6:thm2:approximation_of_rectifier_square_local} it holds for all $x \in \R$ with $\sqrt{2d} \le \abs{x} \le \constantR + \sqrt{2d}$ that 
\begin{equation}
\Abs{[\Rect(\abs{x}-\sqrt{2 d})]^2 - (\functionANN(\psi))(x)} \! \le 4^{-M-1} {\constantR}^2,
\end{equation}

and
\item
\label{item7:thm2:approximation_of_rectifier_square_local} it holds for all $x \in \R$ with $\abs{x} \ge \constantR + \sqrt{2d}$ that 
\begin{equation}
\Abs{[\Rect(\abs{x}-\sqrt{2 d})]^2 - (\functionANN(\psi))(x)} \! \le \! \big[\abs{x}-\sqrt{2 d}\big]^4 {\constantR}^{-2}\ifnocf.
\end{equation}
\end{enumerate}
(cf. \cref{cor:approximation_of_rectifier_square}), let $\lambda \in \R$ satisfy $\lambda = [\int_{\R^d} \abs{g(y)}^2 \varphi(y) \, dy]^{-1/(4d)}$, and let $\idPowerrr, \Psi, \Phi \in \ANNs$ satisfy $\idPowerrr = (\scalarMultANN{\lambda}{\idRelu_1})\compANNbullet \AffineANN_{\lambda, 0}$, $\Psi = \sumANN_{1, d} \compANNbullet \parallelizationSpecial_d(\psi, \psi, \ldots, \psi)$, and $\Phi = (\powANN{\idPowerrr}{d}) \compANNbullet \Psi$ \cfload. \Nobs that \cref{lemma:bounds_for_the_norm_of_specific_function} (applied with $d \with d$, $\varphi \with \varphi$, $g \with g$ in the notation of \cref{lemma:bounds_for_the_norm_of_specific_function}) implies that
\begin{equation}\label{eqn:thm2:lambda}
\begin{split}
0 < \lambda & = \! \pb{\int_{\R^d} \abs{g(y)}^2 \varphi(y) \, dy}^{-\frac{1}{4d}} \! \le \! \pb{50^{-1}d^{\nicefrac{-3}{2}}e^{- d}}^{-\frac{1}{4d}} \! = \! \pb{50d^{\nicefrac{3}{2}} e^{d}}^{\frac{1}{4d}} \\
& \le \! \pb{64 d^2 4^{d}}^{\frac{1}{4d}} \! = \! \pb{8 d \, 2^{d}}^{\frac{1}{2d}} \! \le \! \pb{8^d \, 2^{d}}^{\frac{1}{2d}} \! = \! \pb{16^{d}}^{\frac{1}{2d}} \! = \! \pb{4^{2d}}^{\frac{1}{2d}} \! = 4.
\end{split}
\end{equation}

\noindent
This and \cref{lemma:on_Realization_of_DNN_2} (applied with $n \with d$, $\lambda \with \lambda$, $\Phi \with \Psi$, $\idPowerrr \with \idPowerrr$, $\Psi \with \Phi$ in the notation of \cref{lemma:on_Realization_of_DNN_2}) ensure that for all $x \in \R^{\inDimANN(\Phi)}$ it holds that
\begin{equation}\label{eqn:thm2:Phi_Psi}
(\functionANN(\Phi))(x) = \lambda^{2d} (\functionANN(\Psi))(x) = \!\pb{\int_{\R^d} \abs{g(y)}^2 \varphi(y) \, dy}^{\nicefrac{-1}{2}}\! (\functionANN(\Psi))(x).
\end{equation}

\noindent
Next \nobs that \cref{item1:thm2:approximation_of_rectifier_square_local}, \cref{lem:def:ANNsum} (applied with $m \with 1$, $n \with d$ in the notation of \cref{lem:def:ANNsum}), and \cref{Lemma:PropertiesOfParallelizationEqualLength} (applied with $n \with d$, $(\Phi_1, \Phi_2, \ldots, \Phi_n) \with (\psi, \psi, \ldots, \psi)$ in the notation of \cref{Lemma:PropertiesOfParallelizationEqualLength}) assure that for all $x = (x_1, x_2, \ldots, x_d) \in \R^d$ it holds that $\functionANN(\Psi) \in C(\R^d, \R)$ and
\begin{equation}\label{eqn:thm2:Psi_psi}
(\functionANN(\Psi))(x) = \sum_{j=1}^d (\functionANN(\psi))(x_j).
\end{equation}

\noindent
Combining this with \cref{eqn:thm2:Phi_Psi} establishes \cref{item1:thm2}. In the next step \nobs that \cref{item2:thm2:approximation_of_rectifier_square_local}, \cref{lem:def:ANNsum} (applied with $m \with 1$, $n \with d$ in the notation of \cref{lem:def:ANNsum}),  \cref{Lemma:PropertiesOfParallelizationEqualLengthDims} (applied with $n \with d$, $(\Phi_1, \Phi_2, \ldots, \Phi_n) \with (\psi, \psi, \ldots, \psi)$ in the notation of \cref{Lemma:PropertiesOfParallelizationEqualLengthDims}), and \cref{prop:ANNcomposition_elementary_properties} (applied with $\Phi_1 \with \sumANN_{1,d}$, $\Phi_2 \with \parallelizationSpecial_d(\psi, \psi, \ldots, \psi)$ in the notation of \cref{prop:ANNcomposition_elementary_properties}) show that
\begin{equation}
\dims(\parallelizationSpecial_d(\psi, \ldots, \psi)) = (d, 2d, \underbrace{4d, \ldots, 4d}_{M}, d) \in \N^{M+3}
\end{equation}

\noindent
and
\begin{equation}
\dims(\Psi) = (d, 2d, \underbrace{4d, \ldots, 4d}_{M}, 1) \in \N^{M+3}.
\end{equation}

\noindent
Therefore, we obtain that $\hiddenLength(\Psi) = M + 1$ and
\begin{equation}
\begin{split}
\paramANN(\Psi) & = 2d(d+1) + 4d(2d+1) + \underbrace{4d(4d+1) +\ldots + 4d(4d+1)}_{M-1} + 1(4d+1) \\
& = 10d^2+10d+1+(M-1)(16d^2+4d) \le 21d^2 M.
\end{split}
\end{equation}

\noindent
Combining this with \cref{lemma:on_Realization_of_DNN_2} (applied with $n \with d$, $\lambda \with \lambda$, $\Phi \with \Psi$, $\idPowerrr \with \idPowerrr$, $\Psi \with \Phi$ in the notation of \cref{lemma:on_Realization_of_DNN_2}) ensures that $\hiddenLength(\Phi) = \hiddenLength(\Psi) + d = d + M + 1$ and $\paramANN(\Phi) \le 2 \paramANN(\Psi) + 6d \abs{\outDimANN(\Psi)}^2 \le 42 d^2 M + 6d$. This establishes \cref{item2:thm2,item3:thm2}. Next \nobs that for all $\fW=(w_{i,j})_{(i,j)\in\{1,2,\dots,d\}\times\{1,2,\dots,4d\}}\in\R^{d\times 4d}$, $\fB = (b_1,b_2,\dots,b_d) \in \R^{d}$ it holds that
\begin{equation}
\underbrace{\pmat{1 & 1 & \cdots & 1}}_{\in \R^{1 \times d}} \, \fW =  \! \pmat{\big[\smallsum_{i=1}^d w_{i, 1}\big], \big[\smallsum_{i=1}^d w_{i, 2}\big], \ldots, \big[\smallsum_{i=1}^d w_{i, 4d}\big]} \! \in \R^{1 \times 4d}
\end{equation}

\noindent
and
\begin{equation}
\underbrace{\pmat{1 & 1 & \cdots & 1}}_{\in \R^{1 \times d}} \, \fB + 0 = \big[\smallsum_{i = 1}^d b_i\big] \in \R.
\end{equation}

\noindent
The fact that $\norm{\vectorNN(\parallelizationSpecial_d(\psi, \psi, \ldots, \psi))}_{\infty} = \norm{\vectorNN(\psi)}_{\infty}$ therefore implies that
\begin{equation}
\norm{\vectorNN(\Psi)}_{\infty} = \norm{\vectorNN(\sumANN_{1, d} \compANNbullet \parallelizationSpecial_d(\psi, \psi, \ldots, \psi))}_{\infty} \le d \, \norm{\vectorNN(\parallelizationSpecial_d(\psi, \psi, \ldots, \psi))}_{\infty} = d \, \norm{\vectorNN(\psi)}_{\infty}.
\end{equation}

\noindent
Combining this with \cref{item3:thm2:approximation_of_rectifier_square_local} assures that
\begin{equation}
\norm{\vectorNN(\Psi)}_{\infty} \le d \, \norm{\vectorNN(\psi)}_{\infty} \le d\big(\sqrt{2 d}+1\big) \! \max\{4, {\constantR}^2\} \le 3d^{\nicefrac{3}{2}} \! \max\{4, {\constantR}^2\}.
\end{equation}

\noindent
\cref{lemma:on_Realization_of_DNN_2} (applied with $n \with d$, $\lambda \with \lambda$, $\Phi \with \Psi$, $\idPowerrr \with \idPowerrr$, $\Psi \with \Phi$ in the notation of \cref{lemma:on_Realization_of_DNN_2}) and \cref{eqn:thm2:lambda} hence demonstrate that
\begin{equation}
\begin{split}
\norm{\vectorNN(\Phi)}_{\infty} & \le \max\!\pc{1, \abs{\lambda}} \max\!\pc{\abs{\lambda}, \norm{\vectorNN(\Psi)}_{\infty}} \! \\
& \le 4 \max\!\pc{4, 3d^{\nicefrac{3}{2}} \! \max\{4, {\constantR}^2\}} = 12 d^{\nicefrac{3}{2}} \! \max\{4, {\constantR}^2\}.
\end{split}
\end{equation}

\noindent
This establishes \cref{item4:thm2}. Moreover, \nobs that the fact that for all $\fa_1, \fa_2, \ldots, \fa_d \allowbreak \in \R$ it holds that $(\fa_1+\fa_2+\ldots+\fa_d)^2 \le d (\abs{\fa_1}^2+\abs{\fa_2}^2+\ldots+\abs{\fa_d}^2)$ and \cref{eqn:thm2:Psi_psi} ensure that for all $x = (x_1, x_2, \ldots, x_d) \in \R^d$ it holds that
\begin{equation}
\begin{split}
\abs{(\functionANN(\Psi))(x) - g(x)}^2 & = \! \pabs{\sum_{j=1}^d \!\pb{(\functionANN(\psi))(x_j) - \!\big[\Rect\big(\abs{x_j}-\sqrt{2 d}\big)\big]^2}}^2 \\
& \le d \sum_{j=1}^d \!\pb{(\functionANN(\psi))(x_j) - \!\big[\Rect\big(\abs{x_j}-\sqrt{2 d}\big)\big]^2}^2\!.
\end{split}
\end{equation}

\noindent
Combining this with the fact that for all $k \in \N$ it holds that $\int_{\R^k} (2 \pi)^{\nicefrac{-k}{2}} e^{- \frac{1}{2} \norm{x}_2^2} \, dx = 1$, \cref{lemma:bounds_for_the_norm_of_specific_function} (applied with $d \with d$, $\varphi \with \varphi$, $g \with g$ in the notation of \cref{lemma:bounds_for_the_norm_of_specific_function}), and \cref{eqn:thm2:Phi_Psi} implies that
\begin{equation}
\begin{split}
& \int_{\R^d} \abs{(\functionANN(\Phi))(x) - \fg(x)}^2 \varphi(x) \, dx \\
& = \!\pb{\int_{\R^d}\abs{g(y)}^2 \varphi(y) \, dy}^{-1}\! \int_{\R^d} \abs{(\functionANN(\Psi))(x) - g(x)}^2 \varphi(x) \, dx \\
& \le 50 d^{\nicefrac{3}{2}} e^{d} \int_{\R^d} \abs{(\functionANN(\Psi))(x) - g(x)}^2 \varphi(x) \, dx \\
& \le 50 d^{\nicefrac{5}{2}} e^{d} \int_{\R^d} \! \pb{\sum_{j=1}^d \!\pb{(\functionANN(\psi))(x_j) - \!\big[\Rect\big(\abs{x_j}-\sqrt{2 d}\big)\big]^2}^2} \! \varphi(x_1, x_2, \ldots, x_d) \, d(x_1, x_2, \ldots, x_d) \\
& = 50 d^{\nicefrac{5}{2}} e^{d} \! \pb{\sum_{j=1}^d \int_{\R^d} \!\pb{(\functionANN(\psi))(x_j) - \!\big[\Rect\big(\abs{x_j}-\sqrt{2 d}\big)\big]^2}^2\! \varphi(x_1, x_2, \ldots, x_d) \, d(x_1, x_2, \ldots, x_d)} \\
& = 50 d^{\nicefrac{7}{2}} e^{d} \int_{\R^d} \!\pb{(\functionANN(\psi))(x_1) - \!\big[\Rect\big(\abs{x_1}-\sqrt{2 d}\big)\big]^2}^2\! \varphi(x_1, x_2, \ldots, x_d) \, d(x_1, x_2, \ldots, x_d) \\
& = 25 \sqrt{\frac{2}{\pi}} \, d^{\nicefrac{7}{2}} e^{d}  \int_{\R} \!\pb{(\functionANN(\psi))(x) - \!\big[\Rect\big(\abs{x}-\sqrt{2 d}\big)\big]^2}^2\! e^{- \frac{1}{2} x^2} \, dx.
\end{split}
\end{equation}

\noindent
The integral transformation theorem and \cref{item4:thm2:approximation_of_rectifier_square_local,item5:thm2:approximation_of_rectifier_square_local,item6:thm2:approximation_of_rectifier_square_local,item7:thm2:approximation_of_rectifier_square_local} therefore demonstrate that
\begin{equation}\label{eqn:thm:main2:error}
\begin{split}
& \int_{\R^d} \abs{(\functionANN(\Phi))(x) - \fg(x)}^2 \varphi(x) \, dx \\
& \le  25 \sqrt{\frac{2}{\pi}} \, d^{\nicefrac{7}{2}} e^{d} \int_{\R} \! \pb{(\functionANN(\psi))(x) - \!\big[\Rect\big(\abs{x}-\sqrt{2 d}\big)\big]^2}^2 \! e^{- \frac{1}{2} x^2} \, dx \\
& = 50 \sqrt{\frac{2}{\pi}} \, d^{\nicefrac{7}{2}} e^{d}  \int_{\sqrt{2 d}}^{\infty} \! \pb{(\functionANN(\psi))(x) - \!\big[\Rect\big(\abs{x}-\sqrt{2 d}\big)\big]^2}^2 \! e^{- \frac{1}{2} x^2} \, dx \\
& = 50 \sqrt{\frac{2}{\pi}} \, d^{\nicefrac{7}{2}} e^{d}  \int_{\sqrt{2 d}}^{\constantR + \sqrt{2d}} \! \pb{(\functionANN(\psi))(x) - \!\big[\Rect\big(\abs{x}-\sqrt{2 d}\big)\big]^2}^2 \! e^{- \frac{1}{2} x^2} \, dx \\
& \quad + 50 \sqrt{\frac{2}{\pi}} \, d^{\nicefrac{7}{2}} e^{d}  \int_{\constantR + \sqrt{2 d}}^{\infty} \! \pb{(\functionANN(\psi))(x) - \!\big[\Rect\big(\abs{x}-\sqrt{2 d}\big)\big]^2}^2 \! e^{- \frac{1}{2} x^2} \, dx \\
& \le 50 \sqrt{\frac{2}{\pi}} \, d^{\nicefrac{7}{2}} e^{d} \! \pb{4^{-2M-2} {\constantR}^4 \int_{\sqrt{2 d}}^{\constantR + \sqrt{2d}} e^{- \frac{1}{2} x^2} \, dx + {\constantR}^{-4}  \int_{\constantR + \sqrt{2 d}}^{\infty} \! \big[x - \sqrt{2d}\big]^8  e^{- \frac{1}{2} x^2} \, dx} \\
& = 50 \sqrt{\frac{2}{\pi}} \, d^{\nicefrac{7}{2}} e^{d} \! \pb{16^{-M-1} {\constantR}^4 \int_{0}^{\constantR} e^{- \frac{1}{2} (x^2 + 2x \sqrt{2d} + 2d)} \, dx + {\constantR}^{-4}  \int_{\constantR}^{\infty} x^8 e^{- \frac{1}{2} (x^2 + 2x \sqrt{2d} + 2d)} \, dx} \\
& = 50 \sqrt{\frac{2}{\pi}} \, d^{\nicefrac{7}{2}} \! \pb{16^{-M-1} {\constantR}^4 \int_{0}^{\constantR} e^{- \frac{1}{2} (x^2 + 2x \sqrt{2d})} \, dx + {\constantR}^{-4}  \int_{\constantR}^{\infty} x^8 e^{- \frac{1}{2} (x^2 + 2x \sqrt{2d})} \, dx} \\
& \le 50 \sqrt{\frac{2}{\pi}} \, d^{\nicefrac{7}{2}} \! \pb{16^{-M-1} {\constantR}^4 \int_{0}^{\infty} e^{- \frac{1}{2} x^2} \, dx + {\constantR}^{-4}  \int_{0}^{\infty} x^8 e^{- \frac{1}{2} x^2} \, dx}\!.
\end{split}
\end{equation}

\noindent
Next \nobs that the integral transformation theorem and \cref{lemma:GammaBeta} ensure that
\begin{equation}
\begin{split}
\int_{0}^{\infty} x^8 e^{- \frac{1}{2} x^2} \, dx & = 8\sqrt{2} \int_{0}^{\infty} x^{\nicefrac{7}{2}} e^{- x} \, dx = 8\sqrt{2} \, \Gamma\!\pa{\frac{9}{2}} \\
& = 8\sqrt{2} \! \pb{\frac{7}{2}} \! \pb{\frac{5}{2}} \! \pb{\frac{3}{2}} \! \pb{\frac{1}{2}} \! \Gamma\!\pa{\frac{1}{2}} = \frac{105 \sqrt{\pi}}{\sqrt{2}}.
\end{split}
\end{equation}

\noindent
The fact that $\int_{0}^{\infty} e^{-\frac{1}{2} x^2} dx = \sqrt{\frac{\pi}{2}} \int_{\R} (2 \pi)^{\nicefrac{-1}{2}} e^{-\frac{1}{2} x^2} dx = \sqrt{\frac{\pi}{2}}$ and \cref{eqn:thm:main2:error} therefore assure that
\begin{equation}
\begin{split}
& \int_{\R^d} \abs{(\functionANN(\Phi))(x) - \fg(x)}^2 \varphi(x) \, dx  \\
& \le 50 \sqrt{\frac{2}{\pi}} \, d^{\nicefrac{7}{2}} \! \pb{16^{-M-1} {\constantR}^4 \int_{0}^{\infty} e^{- \frac{1}{2} x^2} \, dx + {\constantR}^{-4}  \int_{0}^{\infty} x^8 e^{- \frac{1}{2} x^2} \, dx} \\
& = 50 d^{\nicefrac{7}{2}} \! \pb{16^{-M-1} {\constantR}^4 + 105 {\constantR}^{-4}}\!.
\end{split}
\end{equation}

\noindent
This establishes \cref{item5:thm2}.
\end{cproof}

\cfclear
\begin{corollary}\label{cor:main_3}
Let $\eps \in (0, 1]$, $\constantfrakC \in [1000 \eps^{-1}, \infty)$, $\constantfrakc \in [\constantfrakC, \infty)$, $d \in \N$, let $\varphi \colon \R^d \to \R$ and $g \colon \R^d \to \R$ satisfy for all $x = (x_1, x_2, \ldots, x_d) \in \R^d$ that $\varphi(x)=(2 \pi)^{\nicefrac{-d}{2}}\exp(- \frac{1}{2} (\smallsum_{j=1}^d \abs{x_j}^2))$ and $g(x) = \smallsum_{j=1}^d [\max\{\abs{x_j}-\sqrt{2 d}, 0\}]^2$, and let $\fg \colon \R^d \to \R$ satisfy for all $x \in \R^d$ that $\fg(x)=[\int_{\R^d}\abs{g(y)}^2 \allowbreak \varphi(y)\,dy]^{-1/2}g(x)$\cfload. Then there exists $\Phi \in \ANNs$ such that $\inDimANN(\Phi) = d$, $\outDimANN(\Phi) = 1$, $d \le \hiddenLength(\Phi) \le \constantfrakc d$, $\norm{\vectorNN(\Phi)}_{\infty} \le \constantfrakc d^{\constantfrakc}$, $\paramANN(\Phi) \le \constantfrakc d^3$, and $[\int_{\R^d}\abs{\pr{\functionANN\pr{\Phi}}(x)-\fg(x)}^2\varphi(x)\, dx]^{\nicefrac{1}{2}} \le \eps$ \cfout.
\end{corollary}
\begin{cproof}{cor:main_3}
Throughout this proof let $M \in \N \cap [2, \infty)$, $\constantR \in [1, \infty)$ satisfy $M = \max((-\infty, \allowbreak \constantR] \allowbreak \cap \N)$ and $\constantR = 9d \eps^{\nicefrac{-1}{2}}$. \Nobs that \cref{thm:main2} (applied with $d \with d$, $M \with M$, $\constantR \with \constantR$, $\varphi \with \varphi$, $g \with g$, $\fg \with \fg$ in the notation of \cref{thm:main2}) ensures that there exists $\Phi \in \ANNs$ which satisfies that
\begin{enumerate}[label=(\Roman *)]
\item
\label{item1:cor:main_3_local} it holds that $\functionANN(\Phi) \in C(\R^d, \R)$,

\item
\label{item2:cor:main_3_local} it holds that $\hiddenLength(\Phi) = d + M + 1$,

\item
\label{item3:cor:main_3_local} it holds that $\paramANN(\Phi) \le 42 d^2 M + 6d$,

\item
\label{item4:cor:main_3_local} it holds that $\norm{\vectorNN(\Phi)}_{\infty} \le 12 d^{\nicefrac{3}{2}} \! \max\{4, {\constantR}^2\}$, and

\item
\label{item5:cor:main_3_local} it holds that $\int_{\R^d} \abs{(\functionANN(\Phi))(x) - \fg(x)}^2 \varphi(x) \, dx \le 50 d^{\nicefrac{7}{2}} \! \pb{16^{-M-1} {\constantR}^4 + 105 {\constantR}^{-4}}\!$\ifnocf.
\end{enumerate}

\noindent
\cfload[.]Therefore, we obtain that $\inDimANN(\Phi) = d$, $\outDimANN(\Phi) = 1$, $d \le \hiddenLength(\Phi) = d + M + 1 \le d + \constantR + 1 = d + 9d \eps^{\nicefrac{-1}{2}} + 1 \le 11 d \eps^{\nicefrac{-1}{2}} \le \constantfrakC d \le \constantfrakc d$, $\norm{\vectorNN(\Phi)}_{\infty} \le 12 d^{\nicefrac{3}{2}} \! \max\{4, {\constantR}^2\} = 972 d^{\nicefrac{7}{2}} \eps^{-1} \le \constantfrakC d^{\constantfrakC} \le \constantfrakc d^{\constantfrakc}$, and $\paramANN(\Phi) \le 42 d^2 M + 6d \le 42 d^2 \constantR + 6d = 378 d^3 \eps^{\nicefrac{-1}{2}} +6d \le 384 d^3 \eps^{\nicefrac{-1}{2}} \le \constantfrakC d^3 \le \constantfrakc d^3$. Moreover, \nobs that the fact that for all $x \in [4, \infty)$ it holds that $x^2 \le 2^{x}$, the assumption that $M = \max((-\infty, \allowbreak \constantR] \allowbreak \cap \N)$, and the assumption that $\constantR = 9d \eps^{\nicefrac{-1}{2}}$ show that $16^{-M-1} {\constantR}^4 + 105 {\constantR}^{-4} \le 106 {\constantR}^{-4} = 106 \eps^2 (9d)^{-4} \le (50 d^{\nicefrac{7}{2}})^{-1} \eps^2$. Combining this with \cref{item5:cor:main_3_local} implies that
$[\int_{\R^d} \abs{(\functionANN(\Phi))(x) - \fg(x)}^2 \varphi(x) \, dx]^{\nicefrac{1}{2}} \allowbreak \le [50 d^{\nicefrac{7}{2}} [16^{-M-1} {\constantR}^4 + 105 {\constantR}^{-4}]]^{\nicefrac{1}{2}} \le \eps$.
\end{cproof}


\section{Lower and upper bounds for the number of ANN parameters in the approximation of high-dimensional \allowbreak functions}\label{sec:lower_and_upper_bounds_for_number_of_parameters_in_NN_approximations}

In \cref{sec:lower_and_upper_bounds_for_number_of_parameters_in_NN_approximations} we combine the lower bounds for the number of parameters of certain ANNs from \cref{sec:lower_bounds_for_number_of_parameters_in_ANN_approximations} with the upper bounds for the number of parameters of certain ANNs from \cref{sec:upper_bounds_for_number_of_parameters_in_ANN_approximations} to establish in \cref{thm:main_result} in Subsection~\ref{subsec:ann_approximation_specifying_the_target_function} below the main ANN approximation result of this article. \cref{thm:introduction} in the introduction is a direct consequence of \cref{cor:main_result} in Subsection~\ref{subsec:ann_approximation_without_specifying_the_target_function} below. The proof of \cref{cor:main_result}, in turn, is based on an application of \cref{thm:main_result}.

\subsection{ANN approximations with specifying the target functions}
\label{subsec:ann_approximation_specifying_the_target_function}

\cfclear
\begin{theorem}\label{thm:main_result}
Let $\varphi_d\colon \R^d \to \R$, $d \in \N$, and $f_d \colon \R^d \to \R$, $d \in \N$, satisfy for all $d \in \N$, $x = (x_1, x_2, \ldots, x_d) \in \R^d$ that $\varphi_d(x)=(2 \pi)^{\nicefrac{-d}{2}}\exp(- \frac{1}{2} (\smallsum_{j=1}^d\abs{x_j}^2))$ and $f_d(x) = \smallsum_{j=1}^d [\max\{\abs{x_j}-\sqrt{2 d}, 0\}]^2$, let $\ff_d \colon \R^d \to \R$, $d \in \N$, satisfy for all $d \in \N$, $x \in \R^d$ that $\ff_d(x) \allowbreak = \allowbreak [\int_{\R^d} \abs{f_d(y)}^2 \allowbreak \varphi_d(y) \, dy]^{-1/2} f_d(x)$, and let $\delta \in (0, 1]$, $\eps \in (0, \nicefrac{1}{2}]$ \cfload. Then there exists $\constantfrakC \in (0, \infty)$ such that
\begin{enumerate}[label=(\roman *)]
\item
\label{item1:thm:main_result} it holds for all $\constantfrakc \in [\constantfrakC, \infty)$, $d \in \N$ that
\begin{equation}
\min \! \pc{p \in \N \colon \!
\pb{
\begin{gathered}
\exists \, \Phi \in \ANNs \colon \, p=\paramANN(\Phi), \, \inDimANN(\Phi) = d, \, \outDimANN(\Phi) = 1, \\
d \le \hiddenLength(\Phi) \le \constantfrakc d, \, \norm{\vectorNN(\Phi)}_{\infty} \le \constantfrakc d^{\constantfrakc}\!, \\
[\smallint\nolimits_{\R^d}\abs{\pr{\functionANN\pr{\Phi}}(x)- \ff_d(x)}^2\varphi_d(x)\,dx]^{\nicefrac{1}{2}} \le \eps
\end{gathered}}}\! \le \constantfrakc d^{3}
\end{equation}

and
\item
\label{item2:thm:main_result} it holds for all $\constantfrakc \in [\constantfrakC, \infty)$, $d \in \N$ that
\begin{equation}
\min \! \pc{p \in \N \colon \!
\pb{
\begin{gathered}
\exists \, \Phi \in \, \ANNs \colon \, p=\paramANN(\Phi), \, \inDimANN(\Phi) = d, \, \outDimANN(\Phi)=1, \\
\hiddenLength(\Phi)\leq \constantfrakc d^{1-\delta}, \, \norm{\vectorNN(\Phi)}_{\infty} \le \constantfrakc d^{\constantfrakc}\!, \\
[\smallint\nolimits_{\R^d} \abs{\pr{\functionANN\pr{\Phi}}(x)- \ff_d(x)}^2\varphi_d(x)\,dx]^{\nicefrac{1}{2}} \le \eps
\end{gathered}}}\! \ge (1+{\constantfrakc}^{-3})^{(d^{\delta})}\!\ifnocf.
\end{equation}
\end{enumerate}
\cfout[.]
\end{theorem}
\begin{cproof}{thm:main_result}
Throughout this proof let $\constantfrakC \in [100 (\delta \ln (1.03))^{-2}, \infty) \cap [1000 \eps^{-1}, \infty)$, $\constantfrakc \in [\constantfrakC, \infty)$, $d \in \N$ satisfy $2 \constantfrakC^{\nicefrac{5}{\delta}} \le (1.03)^{\sqrt{\constantfrakC}}$. \Nobs that \cref{cor:main_3} (applied with $\eps \with \eps$, $\constantfrakC \with \constantfrakC$, $\constantfrakc \with \constantfrakc$, $d \with d$, $\varphi \with \varphi_d$, $g \with f_d$, $\fg \with \ff_d$ in the notation of \cref{cor:main_3}) assures that there exists $\Phi \in \ANNs$ such that $\inDimANN(\Phi) = d$, $\outDimANN(\Phi) = 1$, $d \le \hiddenLength(\Phi) \le \constantfrakc d$, $\norm{\vectorNN(\Phi)}_{\infty} \le \constantfrakc d^{\constantfrakc}$, $\paramANN(\Phi) \le \constantfrakc d^3$, and $[\int_{\R^d}\abs{\pr{\functionANN\pr{\Phi}}(x)-\fg_d(x)}^2\varphi_d(x)\,dx]^{\nicefrac{1}{2}} \le \eps$ \cfload. This establishes \cref{item1:thm:main_result}. Moreover, \nobs that \cref{cor:main_2} (applied with $\varphi_d \with \varphi_d$, $g_d \with f_d$, $\fg_d \with \ff_d$, $\delta \with \delta$, $\constantfrakC \with \constantfrakC$ in the notation of \cref{cor:main_2}) ensures that for all $\Phi \in \ANNs$ with $\inDimANN(\Phi) = d$, $\outDimANN(\Phi)=1$, $\hiddenLength(\Phi) \le \constantfrakc d^{1-\delta}$, $\norm{\vectorNN(\Phi)}_{\infty} \le \constantfrakc d^{\constantfrakc}$, and $[\int_{\R^{d}} \abs{(\functionANN(\Phi))(x) - \ff_{d}(x)}^2 \varphi_{d}(x)\, dx]^{\nicefrac{1}{2}} \le \eps$ it holds that $\paramANN(\Phi) \ge (1+{\constantfrakc}^{-3})^{d^{\delta}}$. This establishes \cref{item2:thm:main_result}.
\end{cproof}

\subsection{ANN approximations without specifying the target functions}
\label{subsec:ann_approximation_without_specifying_the_target_function}

\cfclear
\begin{corollary}\label{cor:main_result}
Let $\varphi_d\colon \R^d \to \R$, $d \in \N$, satisfy for all $d \in \N$, $x \in \R^d$ that $\varphi_d(x)=(2 \pi)^{\nicefrac{-d}{2}}\exp(- \frac{1}{2} \norm{x}_2^2)$ \cfload. Then there exist continuously differentiable $\ff_d \colon \R^d \allowbreak \to \R$, $d \in \N$, such that for all $\delta \in (0,1]$, $\eps \in (0, \nicefrac{1}{2}]$ there exists $\constantfrakC \in (0, \infty)$ such that
\begin{enumerate}[label=(\roman *)]
\item
\label{item1:cor:main_result} it holds for all $\constantfrakc \in [\constantfrakC, \infty)$, $d \in \N$ that
\begin{equation}
\min \! \pc{p \in \N \colon
\pb{
\begin{gathered}
\exists \, \Phi \in \ANNs \colon \, p=\paramANN(\Phi), \, \inDimANN(\Phi) = d, \, \outDimANN(\Phi) = 1, \\
d \le \hiddenLength(\Phi) \le \constantfrakc d, \, \norm{\vectorNN(\Phi)}_{\infty} \le \constantfrakc d^{\constantfrakc}, \\
[\smallint\nolimits_{\R^d}\abs{\pr{\functionANN\pr{\Phi}}(x)- \ff_d(x)}^2\varphi_d(x)\,dx]^{\nicefrac{1}{2}} \le \eps
\end{gathered}}} \le \constantfrakc d^3
\end{equation}

and
\item
\label{item2:cor:main_result} it holds for all $\constantfrakc \in [\constantfrakC, \infty)$, $d \in \N$ that
\begin{equation}
\min \! \pc{p \in \N \colon
\pb{
\begin{gathered}
\exists \, \Phi \in \ANNs \colon \, p=\paramANN(\Phi), \, \inDimANN(\Phi) = d, \, \outDimANN(\Phi)=1, \\
\hiddenLength(\Phi)\leq \constantfrakc d^{1-\delta}, \, \norm{\vectorNN(\Phi)}_{\infty} \le \constantfrakc d^{\constantfrakc},\\
[\smallint\nolimits_{\R^d}\abs{\pr{\functionANN\pr{\Phi}}(x)- \ff_d(x)}^2\varphi_d(x)\,dx]^{\nicefrac{1}{2}} \le \eps
\end{gathered}}} \ge (1 + {\constantfrakc}^{-3})^{(d^{\delta})}\ifnocf.
\end{equation}
\end{enumerate}
\cfout[.]
\end{corollary}
\begin{cproof}{cor:main_result}
Throughout this proof let $f_d \colon \R^d \to \R$, $d \in \N$, satisfy for all $d \in \N$, $x = (x_1, x_2, \ldots, x_d) \in \R^d$ that $f_d(x) = \smallsum_{j=1}^d [\max\{\abs{x_j}-\sqrt{2 d}, 0\}]^2$, let $\ff_d \colon \R^d \to \R$, $d \in \N$, satisfy for all $d \in \N$, $x \in \R^d$ that $\ff_d(x) \allowbreak = \allowbreak f_d(x) \allowbreak [\int_{\R^d} \abs{f_d(y)}^2 \allowbreak \varphi_d(y) \, dy]^{\nicefrac{-1}{2}}$, and let $\delta \in (0, 1]$, $\eps \in (0, \nicefrac{1}{2}]$ \cfload. \Nobs that \cref{thm:main_result} (applied with $(\varphi_d)_{d \in \N} \with (\varphi_d)_{d \in \N}$, $(f_d)_{d \in \N} \with (f_d)_{d \in \N}$, $(\ff_d)_{d \in \N} \with (\ff_d)_{d \in \N}$, $\delta \with \delta$, $\eps \with \eps$ in the notation of \cref{thm:main_result}) establishes \cref{item1:cor:main_result,item2:cor:main_result}.
\end{cproof}

\subsubsection*{Acknowledgements}
Benno Kuckuck and Philippe von Wurstemberger are gratefully acknowledged for their helpful assistance regarding \cref{lem:composition_infnorm}. Joshua Lee Padgett is gratefully acknowledged for his helpful assistance regarding \cref{lemma:recurrence_on_functions_g_f_r,cor:approximation_of_square_on_segment}. The third author acknowledges funding by the Deutsche Forschungsgemeinschaft (DFG, German Research
Foundation) under Germany’s Excellence Strategy EXC 2044-390685587, Mathematics Münster:
Dynamics-Geometry-Structure. The fourth author acknowledges funding by the Austrian Science Fund (FWF) 
through the projects P 30148 and I 3403.


\bibliographystyle{acm}
\bibliography{bibfile}

\end{document}